\documentclass[a4paper, 10pt, final]{amsart}

\usepackage{lmodern}
\usepackage{xcolor}
\usepackage[latin1]{inputenc}
\usepackage[english]{babel}
\usepackage{amsmath}
\usepackage{amsfonts}
\usepackage{amssymb}
\usepackage{graphicx}
\usepackage{amscd}
\usepackage{mathrsfs}
\usepackage{fancybox}

\usepackage{lmodern}
\usepackage[extra]{tipa}

\usepackage{amsmath, amsthm, amsfonts}

\DeclareMathOperator{\reg}{reg}
\DeclareMathOperator{\Lie}{Lie}

\DeclareMathOperator{\weight}{weight}
\DeclareMathOperator{\id}{id}

\DeclareMathOperator{\conv}{conv}

\DeclareMathOperator{\har}{har}
\DeclareMathOperator{\Ad}{Ad}

\DeclareMathOperator{\depth}{depth}

\DeclareMathOperator{\un}{un}
\DeclareMathOperator{\Har}{Har}
\DeclareMathOperator{\DR}{DR}
\DeclareMathOperator{\Wd}{Wd}

\DeclareMathOperator{\Spec}{Spec}

\DeclareMathOperator{\Reg}{Reg}
\DeclareMathOperator{\Cf}{Cf}

\DeclareMathOperator{\LA}{LA}

\DeclareMathOperator{\Li}{Li}

\DeclareMathOperator{\Vect}{Vect}

\DeclareMathOperator{\KZ}{KZ}

\DeclareMathOperator{\dec}{dec}
\DeclareMathOperator{\an}{an}

\theoremstyle{plain}

\newtheorem{Theoreme}{Theoreme}[subsection]
\newtheorem{Proposition}[Theoreme]{Proposition}

\newtheorem{Proposition-Definition}[Theoreme]{Proposition-Definition}
\newtheorem{Lemma-Notation}[Theoreme]{Lemma-Notation}
\newtheorem{Lemma-Definition}[Theoreme]{Lemma-Definition}

\newtheorem{Nota Bene}[Theoreme]{Nota Bene}
\newtheorem{Corollary}[Theoreme]{Corollary}
\theoremstyle{definition}
\newtheorem{Conjecture}[Theoreme]{Conjecture}

\newtheorem{Lemma}[Theoreme]{Lemma}
\newtheorem{Definition}[Theoreme]{Definition}
\newtheorem{Remark}[Theoreme]{Remark}

\newtheorem{Example}[Theoreme]{Example}

\newtheorem{Notation}[Theoreme]{Notation}

\input cyracc.def
\DeclareFontFamily{U}{russian}{}
\DeclareFontShape{U}{russian}{m}{n}
        { <5><6> wncyr5
        <7><8><9> wncyr7
        <10><10.95><12><14.4><17.28><20.74><24.88> wncyr10 }{}
\DeclareSymbolFont{Russian}{U}{russian}{m}{n}
\DeclareSymbolFontAlphabet{\mathcyr}{Russian}
\makeatletter
\let\@math@cyr\mathcyr
\renewcommand{\mathcyr}[1]{\@math@cyr{\cyracc #1}}
\makeatother
\newcommand{\sh}{\mathcyr{sh}} 

\newcommand{\simlra}{\buildrel \sim \over \longrightarrow}

\newcommand{\xdownarrow}[1]{%
{\left\downarrow\vbox to #1{} \right.\kern-\nulldelimiterspace} }

\makeatletter
\newcounter{subsubsubsection}[subsubsection]
\renewcommand\thesubsubsubsection{\thesubsubsection .\@alph\c@subsubsubsection}
\newcommand\subsubsubsection{\@startsection{subsubsubsection}{4}{\z@}%
                                {-3.25ex\@plus -1ex \@minus -.2ex}%
                                     {1.5ex \@plus .2ex}%
                                     {\normalfont\normalsize\bfseries}}
\newcommand*\l@subsubsubsection{\@dottedtocline{3}{10.0em}{4.1em}}
\newcommand*{\subsubsubsectionmark}[1]{}
\makeatother

\usepackage[top=2.5cm, bottom=2.5cm, left=2.5cm, right=2.5cm]{geometry}

\title{}

\author{David Jarossay}

\address{Universit\'{e} de Gen\`{e}ve, Section de mathématiques, 2-4 rue du Li\`{e}vre,
	Case postale 64,
	1211 Gen\`{e}ve, Switzerland}
\email{david.jarossay@unige.ch}

\begin{document}

\title{A bound on the norm of overconvergent $p$-adic multiple polylogarithms}

\maketitle

\begin{abstract} We generalize the definition of overconvergent $p$-adic multiple polylogarithms and of $p$-adic cyclotomic multiple zeta values and we prove a bound on their norm. A byproduct of the proof is a characterization of these objects in terms of certain regularized $p$-adic iterated integrals. The generalization of the definition consists in replacing the underlying Frobenius structure by its iterations. The bound on the norms of overconvergent $p$-adic multiple polylogarithms that we obtain is a prerequisite for our subsequent papers on $p$-adic cyclotomic multiple zeta values.
\newline This is Part I-1 of \emph{$p$-adic cyclotomic multiple zeta values and $p$-adic pro-unipotent harmonic actions}.
\end{abstract}

\tableofcontents

\setcounter{section}{-1}

\numberwithin{equation}{section}

\section{Introduction}

\subsection{Cyclotomic multiple zeta values and multiple polylogarithms}

Cyclotomic multiple zeta values are the following complex numbers, where $(n_{i})_{d}=(n_{1},\cdots,n_{d})$ is a tuple of positive integers, and $(\xi_{i})_{d}=(\xi_{1},\cdots,\xi_{d})$ is a tuple of $N$-th roots of unity, where $N$ is a positive integer and $(\xi_{d},n_{d}) \not= (1,1)$ :
\begin{equation}
\label{eq:multizetas} \zeta \big( (n_{i})_{d} ;
(\xi_{i})_{d} \big) =
\sum_{0<m_{1}<\ldots<m_{d}} \frac{\big( \frac{\xi_{2}}{\xi_{1}}\big)^{m_{1}} \ldots \big( \frac{1}{\xi_{d}}\big)^{m_{d}}}{m_{1}^{n_{1}} \ldots m_{d}^{n_{d}}}
\end{equation}
Denoting by $n=n_{d}+\cdots+n_{1}$ and $(\epsilon_{n},\ldots,\epsilon_{1}) =(\underbrace{0,\ldots,0}_{n_{d}-1},\xi_{d},\ldots,\underbrace{0,\ldots,0}_{n_{1}-1},\xi_{1})$, we have
\begin{equation} \label{eq:multizetas integral}
\zeta \big( (n_{i})_{d} ;
(\xi_{i})_{d} \big) = (-1)^{d}\int_{t_{n}=0}^{1} \frac{dt_{n}}{t_{n}-\epsilon_{n}} \int_{t_{n-1}=0}^{t_{n}} \ldots \int_{t_{1}=0}^{t_{2}} \frac{dt_{1}}{t_{1} - \epsilon_{1}}
\end{equation}

Equation (\ref{eq:multizetas integral}) shows that cyclotomic multiple zeta values are periods of the comparison between the Betti and De Rham realizations of the pro-unipotent fundamental groupoid ($\pi_{1}^{\un}$) of $\mathbb{P}^{1} \setminus \{0,\mu_{N},\infty\}$ (\cite{Deligne Goncharov}, \S5.16). 
\newline\indent Multiple polylogarithms (\cite{Goncharov}) are the following power series, which are convergent for $z \in \mathbb{C}$ such that $|z|<1$ :
\begin{equation} \label{eq: power series mpl} 
\Li\big((n_{i})_{d};(\xi_{i})_{d}\big)(z) = (-1)^{d} \sum_{0<m_{1}<\ldots<m_{d}} \frac{\big( \frac{\xi_{2}}{\xi_{1}}\big)^{m_{1}} \ldots \big( \frac{z}{\xi_{d}}\big)^{m_{d}}}{m_{1}^{n_{1}} \ldots m_{d}^{n_{d}}} 
\end{equation}
They are solutions to the KZ equation on $\mathbb{P}^{1} \setminus \{0,\mu_{N},\infty\}$ :
\begin{equation} \label{eq:KZ equation} d\Li\big((n_{i})_{d};(\xi_{i})_{d}\big)(z) = \left\{ \begin{array}{ll} \displaystyle\frac{dz}{z}\Li\big((n_{i})_{d-1},n_{d}-1;(\xi_{i})_{d}\big) &\text{ if }n_{d}\geqslant  2
\\ \displaystyle\frac{dz}{z-\xi_{d}}\Li\big((n_{i})_{d-1};(\xi_{i})_{d-1}\big) &\text{ if }n_{d}= 1 \end{array} \right.
\end{equation}
\indent Multiple polylogarithms can be extended to multivalued holomorphic functions on $(\mathbb{P}^{1} \setminus \{0,\mu_{N},\infty\})(\mathbb{C})$ \cite{Goncharov}, 
defined as iterated path integrals in the sense of Chen \cite{Chen} (with a regularization if needed) : for $\gamma : [0,1] \rightarrow \mathbb{P}^{1}(\mathbb{C})$ a differentiable topological path such that $\gamma\big((0,1)\big) \in (\mathbb{P}^{1} \setminus \{0,\mu_{N},\infty\})(\mathbb{C})$, let
\begin{equation} \label{eq:integral mpl} \Li \big((n_{i})_{d};(\xi_{i})_{d}\big)(\gamma) = \int_{t_{n}=0}^{1} \gamma^{\ast} (\frac{dz}{z-\epsilon_{n}})(t_{n})\int_{t_{n-1}=0}^{t_{n}} \ldots  \int_{t_{1}=0}^{t_{2}} \gamma^{\ast}(\frac{dz}{z-\epsilon_{1}})(t_{1}) 
\end{equation}

\subsection{$p$-adic cyclotomic multiple zeta values, $p$-adic multiple polylogarithms and their overconvergent variants}

Let $p$ be a prime number which does not divide $N$.
Let $K$ be a totally unramified finite extension of $\mathbb{Q}_{p}$ containing the $N$-th roots of unity in $\overline{\mathbb{Q}_{p}}$. The KZ equation (\ref{eq:KZ equation}) has a Frobenius structure over $K$, and this defines the crystalline pro-unipotent fundamental groupoid of $\mathbb{P}^{1} \setminus \{0,\mu_{N},\infty\}$ (\cite{Deligne} \S11, \cite{CLS}, \cite{Shiho 1}, \cite{Shiho 2}).
\newline\indent The theory of Coleman integration (\cite{Coleman}, \cite{Besser}, \cite{Vologodsky}) enables, by means of this Frobenius structure, to define $p$-adic analogues of the integrals (\ref{eq:multizetas integral}) and (\ref{eq:integral mpl}) : this defines $p$-adic cyclotomic multiple zeta values $\zeta_{p}^{\KZ}\big(
(n_{i})_{d} ; (\xi_{i})_{d} \big) \in K$ and $p$-adic multiple polylogarithms 
$\Li_{p}^{\KZ}\big( (n_{i})_{d};(\xi_{i})_{d} \big)$ (\cite{Furusho 1},\cite{Furusho 2} for $N=1$, \cite{Yamashita} for any $N$) (from now on we drop the assumption $(n_{d},\xi_{d})\not=(1,1)$ used above). The power series expansion of $\Li_{p}^{\KZ}\big( (n_{i})_{d};(\xi_{i})_{d} \big)$ at $0$ is convergent for $z \in K$ such that $|z|_{p}<1$ and is identical to (\ref{eq: power series mpl}).
\newline\indent Our goal in this paper is to give the foundations of an explicit theory of $p$-adic cyclotomic multiple zeta values. For computations and explicit formulas, it is more convenient to consider another type of $p$-adic analogues of the integrals (\ref{eq:multizetas integral}) and (\ref{eq:integral mpl}). Let $\alpha$ be a positive integer. We consider the "Frobenius of the KZ equation iterated $\alpha$ times" (Definition \ref{def iterated Frobenius}), and we define the associated $p$-adic cyclotomic multiple zeta values (Definition \ref{def of p-adic multiple zeta values}), as numbers 
\begin{equation} \zeta_{p,\alpha}\big((n_{i})_{d};(\xi_{i})_{d}) \in K 
\end{equation}
The simple relation between $\zeta_{p,\alpha}$ and $\zeta_{p}^{\KZ}$ is explained in \cite{I-3} (see also \cite{Yamashita}, Proposition 3.10 and \cite{Furusho 2}, Theorem 2.8) for particular cases. We also define (Definition \ref{def of overconvergent}) the overconvergent $p$-adic multiple polylogarithms
\begin{equation} \Li_{p,\alpha}^{\dagger}\big((n_{i})_{d};(\xi_{i})_{d}) \end{equation}
which are overconvergent rigid analytic functions on a certain convenient affinoid subspace of $\mathbb{P}^{1,\an}/K$. The definitions of $\zeta_{p,\alpha}$ and $\Li_{p,\alpha}^{\dagger}$ generalize definitions in \cite{Deligne Goncharov}, \cite{Unver MZV} for $N=1$, and in \cite{Yamashita}, \cite{Unver cyclotomic} for any $N$, which correspond to particular values of $\alpha$.

\subsection{Summary of the paper}

The setting and definitions are established in \S1. An expression of $\Li_{p,\alpha}^{\dagger}$ in terms of $\Li_{p}^{\KZ}$ and $\zeta_{p,\alpha}$, equation (\ref{eq:functional}), is established and studied in \S2. We introduce and study a notion of regularized $p$-adic iterated integrals (\S3), and, using \S2 and \S3, we characterize $\Li_{p,\alpha}^{\dagger}$ and $\zeta_{p,\alpha}$ in terms of regularized iterated integrals (\S4). This characterization gives our main result, which is a prerequisite for our explicit theory of $p$-adic cyclotomic multiple zeta values, established in \cite{I-2}, \cite{I-3} and subsequent papers.
\newline\indent In the statement below, the norm $||\text{ }||$ used is the norm on the $K$-algebra of global rigid analytic functions on the affinoid $U^{\an}=\mathbb{P}^{1,\an} \setminus \displaystyle\cup_{\xi \in \mu_{N}(K)} B(\xi,1)$ over $K$ defined in terms of the power series expansion at $0$ as $||\sum\limits_{m \in \mathbb{N}} a_{m}z^{m}|| = \sup\limits_{m \in \mathbb{N}} |a_{m}|_{p}$. Here $B(\xi,1)= \{z \in K\text{ }|\text{ }|z-\xi|_{p}<1\}$.
\newline A word of depth $d$ and weight $n$ is a sequence $w=\big((n_{i})_{d};(\xi_{i})_{d})$ as above, with $n_{1}+\ldots+n_{d}=n$.
\newline 
\newline \textbf{Theorem.} \emph{For any $d\in \mathbb{N}^{\ast}$, we have $\displaystyle\max_{\substack{w\text{ word}\\ \weight(w)=n \\ \depth(w)=d}}||\Li_{p,\alpha}^{\dagger}(w)|| \underset{n \rightarrow \infty}{\rightarrow} 0$.}
\newline 
\newline An explicit upper bound for any $||\Li_{p,\alpha}^{\dagger}(w)||$ can be obtained by our proof ; however, what we need for \cite{I-2} and subsequent papers is only the above information. Our method also gives bounds on the norms $|\zeta_{p,\alpha}(w)|_{p}$ (\S4.2), but this is not needed for our subsequent papers. In the $N=1$ case, a more precise bound on $|\zeta_{p,\alpha}(w)|_{p}$ can be deduced from the works of Akagi-Hirose-Yasuda \cite{Akagi Hirose Yasuda} and Chatsiztamatiou \cite{Chatzis}, who use completely different methods and work with lifts of Frobenius-invariant paths.
\newline\indent Appendix A is a description of the $K$-algebra of global rigid analytic function of the space $U^{\an}$ mentioned above, which arises as a reformulation of a classical theorem of Mahler \cite{Mahler} in terms of the rigid analytic space $U^{\an}$. It is a prerequisite for \S3 and \S4. We will investigate its further meaning in terms of more general rigid analytic spaces in a subsequent paper.
\newline\indent Appendix B is a byproduct of the proof : considering all $p$'s and $\alpha$'s at the same time, and we define a certain explicit space of adelic sums of series, which plays the role of a $p$-adic analogue of the type of series appearing in (\ref{eq:multizetas}) and enables to express the overconvergent $p$-adic multiple polylogarithms, resp. the $p$-adic cyclotomic multiple zeta values ; it will reappear in all our explicit theory of $p$-adic cyclotomic multiple zeta values.
\newline\indent The strategy of the paper is inspired by a strategy for computing $p$-adic cyclotomic multiple zeta values, suggested by Deligne in depth one in his foundational paper \cite{Deligne}, \S19.6, and then used by \"{U}nver in \cite{Unver MZV}, \cite{Unver cyclotomic}, in depth one and two. \"{U}nver's paper \cite{U3}, which deals with the $N=1$ case, was written at the same time with the present paper, and \cite{U4} appeared after the present paper ; the main result of
\cite{Unver MZV} \cite{Unver cyclotomic} \cite{U3} \cite{U4} is an algorithm to compute $p$-adic cyclotomic multiple zeta values (Theorem 1.1 of \cite{U4}),  which is implied by the results of the present paper. The explicit formulas for $p$-adic cyclotomic multiple zeta values which we can deduce from our proof, as well as those in \cite{I-2} and \cite{I-3}, are different from \"{U}nver's.
\newline 
\newline\indent \emph{Acknowledgments.} I thank Francis Brown for having asked me to compute $p$-adic multiple zeta values. I thank Sinan \"{U}nver for his encouragements and answers to my questions during discussions on this paper in 2013 and in 2014. I thank Pierre Cartier for encouragements in 2014 and 2015, and Jean-Pierre Serre for a remark communicated to me by Pierre Cartier in 2014, which led me to write the Appendix A. I have been helped in computer experiments independently by Aurel Page and by Annick Valibouze and also had helpful mathematical discussions with them in 2013. I thank an anonymous referee for important remarks and suggestions. This work has been written from 2012 to 2016 and revised afterwards, at Universit\'{e} Paris Diderot, Universit\'{e} de Strasbourg and Universit\'{e} de Gen\`{e}ve with, respectively, support of ERC grant 257638, Labex IRMIA and NCCR SwissMAP.

\numberwithin{equation}{subsection}

\section{Definitions}

We review some the pro-unipotent fundamental groupoid of $\mathbb{P}^{1} \setminus \{0,\mu_{N},\infty\}$ (\S1.1), and we generalize the definition of $p$-adic cyclotomic multiple zeta values (\S1.2) and of overconvergent $p$-adic multiple polylogarithms (\S1.3) by replacing the Frobenius by the "iterated Frobenius", which depends on the number of iterations $\alpha$ in $\mathbb{N}^{\ast} \cup - \mathbb{N}^{\ast}$.

\subsection{Preliminaries on the pro-unipotent fundamental groupoid of $\mathbb{P}^{1} \setminus \{0,\mu_{N},\infty\}$}

\subsubsection{Non-commutative formal power series and shuffle Hopf algebras}

For any $\mathbb{Q}$-algebra $C$, and any set of formal variables $\frak{a}= \{a_{1},\cdots,a_{r}\}$, let $C\langle a_{1},\ldots,a_{r}\rangle$, resp. $C\langle\langle a_{1},\ldots,a_{r}\rangle\rangle$ be the non-commutative $C$-algebra of polynomials, resp. of formal power series over the non-commuting variables $a_{1},\ldots,a_{r}$, with coefficients in $C$. Let $\text{Wd}(a)$ be the set of words over the alphabet $a=\{a_{1},\ldots,a_{r}\}$, including the empty word $\emptyset$. It is a basis of the free $C$-module $C\langle a_{1},\ldots,a_{r}\rangle$.

\begin{Notation}
For $f=f(a_{1},\ldots,a_{r}) \in C\langle \langle a_{1},\ldots,a_{r} \rangle\rangle$ and $w \in \text{Wd}(a)$, let $f[w]\in C$ be the coefficient of $w$ in $f$ ; i.e. we have $f = \sum\limits_{w \in \Wd(a)} f[w]w$. The notation $f[w]$ extends to linear combinations of words by linearity.
\end{Notation}

The shuffle Hopf algebra $\mathcal{O}^{\sh,\frak{a}}$ is the $\mathbb{Q}$-vector space $\mathbb{Q}\langle a_{1},\ldots,a_{r} \rangle$, graded by the length of words over $\frak{a}$, endowed with the shuffle product $\sh$ defined by $(a_{i_{n+n'}}\ldots a_{i_{n+1}})\text{ }\sh\text{ }(a_{i_{n}} \ldots a_{i_{1}}) =
\sum\limits_{\sigma}
a_{i_{\sigma^{-1}(n+n')}} \ldots a_{i_{\sigma^{-1}(1)}}$ where the sum is over permutations 
$\sigma$ of $\{1,\ldots,n+n'\}$ such that $\sigma(n)>\ldots>\sigma(1)$\ and  $\sigma(n+n')>\ldots>\sigma(n+1)$
; the deconcatenation coproduct $\Delta : a_{i_{n}}\ldots a_{i_{1}} \mapsto \sum_{n'=0}^{n} a_{i_{n}}\ldots a_{i_{n'+1}} \otimes a_{i_{n'}} \ldots a_{i_{1}}$ ;  the counit $\epsilon$ equal to the augmentation map ; the antipode $S : a_{i_{n}}\ldots a_{i_{1}} \mapsto (-1)^{n} a_{i_{1}}\ldots a_{i_{n}}$. (We order the words from the right to the left in order to follow the conventions in the literature on  $p$-adic multiple zeta values.) We have
\begin{equation} \label{eq:shuffle equation}
\begin{array}{cc}
\Spec(\mathcal{O}^{\sh,a})(C) & = \{ f \in C \langle\langle a_{1},\ldots,a_{r} \rangle\rangle \text{ }|\text{ }\forall w,w' \in \Wd(a), f[w\text{ }\sh\text{ }w']=f[w]f[w'],\text{ and }f[\emptyset] = 1 \}  \end{array}
\end{equation}

\subsubsection{The groupoid  $\pi_{1}^{\un,\DR}(\mathbb{P}^{1} \setminus \{0,\mu_{N},\infty\})$ and the canonical base-point $\omega_{\DR}$}

The De Rham pro-unipotent fundamental groupoid ($\pi_{1}^{\un,\DR}$) of smooth algebraic varieties over a field of characteristic zero is defined in \cite{Deligne}, \S10.27, \S10.30,(ii). Let $X$ be $\mathbb{P}^{1} \setminus \{0,\mu_{N},\infty\}$ over a field $K$ of characteristic zero which contains a primitive $N$-th root of unity.
The De Rham pro-unipotent fundamental groupoid of $X$ is a groupoid in pro-affine schemes over $X$, with the following base-points : the points of $X$, the points of punctured tangent spaces $T_{x} - \{0\}$, $x \in \{0,\infty\}\cup \mu_{N}(K)$ called tangential base-points (\cite{Deligne}, \S15), and the canonical base-point $\omega_{\DR}$ (\cite{Deligne}, (12.4.1)).
\newline \indent Namely, if $x,y$ are two base-points, we have a pro-affine scheme over $\mathbb{Z}$ denoted $\pi_{1}^{\un,\DR}(X,y,x)$, whose points are called the pro-unipotent De Rham paths from $x$ to $y$ ; we denote by $\pi_{1}^{\un,\DR}(X,x) = \pi_{1}^{\un,\DR}(X,x,x)$, which is a pro-unipotent group scheme ; if $x,y,z$ are three base-points, we have a groupoid multiplication which is a morphism of schemes $\pi_{1}^{\un,\DR}(X,z,y) \times \pi_{1}^{\un,\DR}(X,y,x) \rightarrow \pi_{1}^{\un,\DR}(X,z,x)$ ; the groupoid multiplication makes each $\pi_{1}^{\un,\DR}(X,y,x)$ into a bi-torsor under $(\pi_{1}^{\un,\DR}(X,y),\pi_{1}^{\un,\DR}(X,x))$. (Following the convention of \cite{Deligne} and \cite{Deligne Goncharov}, we denote the groupoid multiplication from the right to the left.)
\newline\indent By \cite{Deligne}, \S12, we have canonical isomorphisms of schemes $\pi_{1}^{\un,\DR}(X,y,x) \simeq \pi_{1}^{\un,\DR}(X,\omega_{\DR})$, which are compatible with the groupoid structure, and, canonical paths ${}_y 1_{x} \in \pi_{1}^{\un,\DR}(X,y,x)(K)$, which are compatible with the groupoid structure.
\newline\indent This reduces the description of the groupoid $\pi_{1}^{\un,\DR}(X)$ to the one of the group scheme $\pi_{1}^{\un,\DR}(X,\omega_{\DR})$.
\newline\indent By \cite{Deligne}, \S12.9, the affine group scheme $\Pi=\pi_{1}^{\un,\DR}(X,\omega_{\DR})$ is canonically isomorphic to the exponential of the completed free Lie algebra over the generators $e_{x}$, $x \in \{0\}\cup \mu_{N}(K)$. This is also $\Spec(\mathcal{O}^{\sh,e_{0\cup\mu_{N}}})$ where $e_{0 \cup \mu_{N}}$ is the alphabet $\{e_{x}\text{ }|\text{ }x \in \{0\}\cup \mu_{N}(K) \}$. This object is described explicitly by \S1.1.1.

\subsubsection{The KZ connection and its canonical formal solution}

By \cite{Deligne}, (12.5.5), (12.12.1), one has the following canonical connection on $\pi_{1}^{\un,\DR}(X,\omega_{\DR}) \times X$ called the Knizhnik-Zamolodchikov connection and denoted by $\nabla_{\KZ}$, which appeared in equation (\ref{eq:KZ equation}):
\begin{equation} f \mapsto df - \bigg( e_{0}\frac{dz}{z} +  \sum_{\xi \in \mu_{N}(K)} e_{\xi} \frac{dz}{z-\xi} \bigg) f 
\end{equation}
\indent It has a canonical formal solution $\Li \in K[[z]][\log(z)]\langle\langle e_{0\cup \mu_{N}} \rangle\rangle$, whose coefficient $\Li[e_{0}^{n_{d}-1}e_{\xi_{d}}\ldots e_{0}^{n_{1}-1}e_{\xi_{1}}]$ is the power series (\ref{eq: power series mpl}), and whose coefficient $\Li[e_{0}]$ is $\log(z)$.
\newline\indent The cyclotomic multiple harmonic sums are the numbers which arise as coefficients in the power series (\ref{eq: power series mpl}) :
\begin{equation} \label{eq:multiple harmonic sums}
\frak{h}_{m} \big((n_{i})_{d};(\xi_{i})_{d+1}\big) = m^{n_{1}+\ldots+n_{d}} \sum_{0<m_{1}<\ldots<m_{d}<m}
\frac{\big( \frac{\xi_{2}}{\xi_{1}} \big)^{m_{1}} \ldots \big(\frac{\xi_{d+1}}{\xi_{d}}\big)^{m_{d}}
\big(\frac{1}{\xi_{d+1}} \big)^{m}}{m_{1}^{n_{1}}\ldots m_{d}^{n_{d}}} 
\end{equation}
where $m$, $d$ and the $n_{i}$'s are positive integers and the $\xi_{i}$'s are $N$-th roots of unity.

\subsubsection{Iterations of the crystalline Frobenius of $\pi_{1}^{\un,\DR}(X_{K})$}

The notion of crystalline $\pi_{1}^{\un}$ is first defined in \cite{Deligne}, \S11, \S13.6, as an enrichment by a crystalline Frobenius structure of the notion of $\pi_{1}^{\un,\DR}$. Variants of this notion, defined more directly and in Tannakian terms, are given in \cite{CLS} and in \cite{Shiho 1} \cite{Shiho 2}. We follow the point of view of \cite{Deligne} which is the most elementary one.
\newline\indent Let $p$ be a prime number which does not divide $N$. Let $k=\mathbb{F}_{q}$, a finite field of characteristic $p$, which contains a primitive $N$-th root of unity. Let $R=W(k)$ its ring of Witt vectors ; it is the ring generated by $\mathbb{Z}_{p}$ and the roots of unity of order prime to $p$ whose reduction is in $\mathbb{F}_{q}$. Let $K$ be the field of fractions of $R$, equal to the field generated by $\mathbb{Q}_{p}$ and the same roots of unity. Let $X_{k}$, $X_{R}$ and $X_{K}$ be $\mathbb{P}^{1} \setminus \{0,\mu_{N},\infty\}$ over, respectively, $k$, $R$ and $K$. Let $\sigma : R \rightarrow R$ be the Frobenius automorphism of $R$, which lifts the $p$-th power Frobenius of $k$. We denote also by $\sigma : K\rightarrow K$ its extension to $K$. It is an automorphism of $K$ which is an isometry for the $p$-adic metric.
\newline\indent Let $\alpha$ be a positive integer, and $\sigma^{\alpha}=\underbrace{\sigma \circ \cdots \circ \sigma}_{\alpha}$. For any $\xi \in \mu_{N}(K)$, we have $\sigma^{\alpha}(\xi)=\xi^{p^{\alpha}}$. Let $X_{R}^{(p^{\alpha})}$ be the pull-back of $X_{R}$ by $\sigma^{\alpha}$ and $X_{K}^{(p^{\alpha})} = X_{R}^{(p^{\alpha})} \times_{\Spec R} \Spec K$.
By functoriality of the De Rham pro-unipotent fundamental groupoid, 
the description of $\pi_{1}^{\un,\DR}(X_{K})$ given in \S1.1.1 provides a similar description of  $\pi_{1}^{\un,\DR}(X_{K}^{(p^{\alpha})})$. We denote by $\omega^{(p^{\alpha})}_{\DR}$ the pull-back of the canonical base-point $\omega_{\DR}$ by $\sigma^{\alpha}$. We denote by $\mathcal{O}^{\sh,e_{0 \cup \mu_{N}^{(p^{\alpha})}}}$ the Hopf algebra of  $\pi_{1}^{\un,\DR}(X_{K}^{(p^{\alpha})},\omega_{\DR}^{(p^{\alpha})})$ : this is the shuffle Hopf algebra over the alphabet $e_{0 \cup \mu_{N}^{(p^{\alpha})}}= \{e_{0}\} \cup \{e_{\xi^{(p^{\alpha})}}\text{ }|\text{ }\xi \in \mu_{N}(K)\}$ (we write $\xi^{(p^{\alpha})}$ and not $\xi^{p^{\alpha}}$ in order to keep track of the fact that we are working on $X_{K}^{(p^{\alpha})}$ and not on $X_{K}$ ; this will be useful : see Definition \ref{the big alphabet}). We denote by  $\nabla_{\KZ}^{(p^{\alpha})}$ the pull-back of the KZ equation by $\sigma^{\alpha}$, namely, the connection on $\pi_{1}^{\un,\DR}(X_{K}^{(p^{\alpha})},\omega_{\DR}^{(p^{\alpha})}) \times X_{K}^{(p^{\alpha})}$ defined by
\begin{equation} f \mapsto df - \bigg( e_{0} \frac{dz'}{z'} + \sum_{\xi \in \mu_{N}(K)} e_{\xi^{(p^{\alpha})}} \frac{dz'}{z'-\xi^{p^{\alpha}}} \bigg) f
\end{equation}
\newline\indent By \cite{Deligne}, \S13.6, the crystalline Frobenius of $\pi_{1}^{\un,\DR}(X_{K})$ is a canonical $\sigma$-linear isomorphism of groupoids 
$$ \phi : \pi_{1}^{\un,\DR}(X_{K}^{(p)}) \simlra \pi_{1}^{\un,\DR}(X_{K}) $$
equal to the inverse of the Frobenius $F_{X/K\ast}$ defined in \cite{Deligne} \S11.11 (11.11.2), (11.11.3). We propose to study, more generally, the following :

\begin{Definition} \label{def iterated Frobenius}For $\alpha$ a positive integer, let
$\phi_{\alpha} = \phi \circ \sigma^{\ast}(\phi) \circ \cdots \circ (\sigma^{\alpha-1})^{\ast}(\phi)$ and $\phi_{-\alpha} =\phi_{\alpha}^{-1}$.
\end{Definition}
$\phi_{\alpha}$ is a canonical $\sigma^{\alpha}$-linear isomorphism $\pi_{1}^{\un,\DR}(X_{K}^{(p^{\alpha})}) \simlra \pi_{1}^{\un,\DR}(X_{K})$.
\newline\indent For $\alpha$ divisible by the integer $\frac{\log(q)}{\log(p)}$, i.e. $\alpha$ such that $p^{\alpha}$ is a power of $q$, we have $\sigma^{\alpha}=\id$ ; in particular, the source and target of $\phi_{\alpha}$ are equal, and $\phi_{\alpha}$ and is an iteration of $\phi_{\frac{\log(q)}{\log(p)}}$. More generally, for any positive integer $\alpha$, we will abusively call $\phi_{\pm\alpha}$ an "iteration of $\phi_{\pm 1}$".

\subsection{$p$-adic cyclotomic multiple zeta values, $p$-adic multiple polylogarithms and their overconvergent variants : definitions}

In the definitions, we use the canonical isomorphisms $\pi_{1}^{\un,\DR}(X,y,x) \simeq \pi_{1}^{\un,\DR}(X_{K},\omega_{\DR})=\Spec(\mathcal{O}^{\sh,e_{0\cup \mu_{N}}})$ mentioned in \S1.1.2.

\subsubsection{$p$-adic cyclotomic multiple zeta values}

We now generalize the notion of $p$-adic cyclotomic multiple zeta values, by replacing the Frobenius by the iterated Frobenius. Below, $\vec{v}_{x}$ means the tangent vector $\vec{v}$ at $x$ ; it is a tangential base-point of $\pi_{1}^{\un,\DR}(X_{K})$ in the sense reviewed in \S1.1.2.

\begin{Notation} (i) (\cite{Deligne Goncharov}, \S5)
For all $x \in \{0\}\cup \mu_{N}(K)$, let
$\Pi_{x,0} = \pi_{1}^{\un,\DR}(X_{K},\vec{1}_{x},\vec{1}_{0})(K)$ \newline (ii) For all $x \in \{0\}\cup \mu_{N}(K)$, let
$\Pi^{(p^{\alpha})}_{x^{p^{\alpha}},0} = \pi_{1}^{\un,\DR}(X_{K}^{(p^{\alpha})},\vec{1}_{x^{p^{\alpha}}},\vec{1}_{0})(K)$
\newline (iii) (\cite{Deligne Goncharov}, \S5.16) Let 
$\tau :\mathbb{G}_{m}(\mathbb{Q}) \times K \langle\langle e_{0\cup\mu_{N}} \rangle\rangle \rightarrow K \langle\langle e_{0\cup\mu_{N}}\rangle\rangle$ be the group action defined by 
\newline $(\lambda,f\big((e_{x})_{x \in \{0\}\cup \mu_{N}(K)}) \mapsto f\big((\lambda e_{x})_{x \in \{0\}\cup \mu_{N}(K)})$.
\end{Notation}

The map $\tau$ multiplies the coefficient of a word $w$ in a power series $f \in K \langle\langle e_{0\cup\mu_{N}}\rangle\rangle$ by $\lambda^{\weight(w)}$, where the weight of a word is its number of letters.
\newline\indent The $p$-adic cyclotomic multiple zeta values can be defined as coefficients of the unique element of $\Pi_{1,0}(K)$ which is invariant by $\phi_{\frac{\log(q)}{\log(p)}}$ (\cite{Furusho 1} Definition 2.17, \cite{Furusho 2} Theorem 2.5, for $N=1$ ; \cite{Yamashita} Definition 2.4 and \S3.3 for any $N$). That point of view follows the theory of Coleman integration \cite{Coleman}, \cite{Besser}, \cite{Vologodsky} ; we will study that notion and relate it to the one below in \cite{I-3}.
\newline\indent Alternatively, $p$-adic cyclotomic multiple zeta values can be defined as coefficients of the image by the Frobenius of the canonical path in $\Pi_{1,0}(K)$ in the sense introduced in \S1.1.2. This has been defined for $N \in \{1,2\}$, $\alpha=1$ (\cite{Deligne Goncharov} \S5.28) ; $N=1$, $\alpha=-1$ (\cite{Unver MZV} ,\S1) ; any $N$ and $\alpha = \frac{\log(q)}{\log(p)}$, (\cite{Yamashita}, Definition 3.1) ; any $N$ and $\alpha=-1$ (\cite{Unver cyclotomic}, \S2.2.3). We generalize these definitions. Below, the notation ${}_y 1_{x}$ refers to the canonical path from $x$ to $y$.

\begin{Definition} \label{def of p-adic multiple zeta values} For any positive integer $\alpha$, let
$$ \Phi_{p,\alpha} = \big(\tau(p^{\alpha})\circ \phi_{\alpha}\big) ( _{\vec{1}_{1}} 1 _{\vec{1}_{0}}) \in \Pi_{1,0}(K) $$ 
$$ \Phi_{p,-\alpha} = \phi_{-\alpha} ( _{\vec{1}_{1}} 1 _{\vec{1}_{0}}) \in \Pi^{(p^{\alpha})}_{1,0}(K) $$ Let $p$-adic cyclotomic multiple zeta values be the following numbers, where the $n_{i}$'s are positive integers and the $\xi_{i}$'s are $N$-th roots of unity :
$$ \zeta_{p,\alpha} \big( (n_{i})_{d}; (\xi_{i})_{d} \big) = (-1)^{d} \Phi_{p,\alpha}[e_{0}^{n_{d}-1}e_{\xi_{d}} \ldots e_{0}^{n_{1}-1}e_{\xi_{1}}] \in K $$ 
$$ \zeta_{p,-\alpha}
\big((n_{i})_{d}; (\xi_{i}^{(p^{\alpha})})_{d} \big) = (-1)^{d} \Phi_{p,-\alpha}[e_{0}^{n_{d}-1}e_{\xi_{d}^{(p^{\alpha})}} \ldots e_{0}^{n_{1}-1}e_{\xi_{1}^{(p^{\alpha})}}] \in K $$
\end{Definition}

Here, the presence of the factor $\tau(p^{\alpha})$ is our convention which is meant to be adapted to the computations in our subsequent papers.
\newline\indent 
We generalize to our $p$-adic cyclotomic multiple zeta values the conjectural description of their algebraic relations stated for $\alpha=1$, $N=1$, \cite{Deligne Goncharov}, and which is implicit in \cite{Yamashita} for any $N$ and $\alpha = \frac{\log(q)}{\log(p)}$.

\begin{Conjecture} \label{p-adic complex conjecture} For any positive integer $\alpha$  the correspondence $\zeta_{p,\alpha}(w) \mapsto \zeta(w)$ defines an isomorphism of algebras over the $N$-th cyclotomic field, from the algebra generated by the $p$-adic cyclotomic multiple zeta values $\zeta_{p,\alpha}(w)$, to the algebra generated by cyclotomic multiple zeta values (\ref{eq:multizetas}) moded out by the ideal generated by $\zeta(2)$.
\end{Conjecture}

\subsubsection{$p$-adic multiple polylogarithms}

The $p$-adic multiple polylogarithms are defined as coefficients of Frobenius-invariant paths, namely, i.e. they are Coleman integrals, as are the numbers $\zeta_{p}^{\KZ}$ evoked in \S1.2. Let us fix a determination $\log_{p}$ of the $p$-adic logarithm. Let $\Li_{p,X_{K}}^{\KZ}$ resp. $\Li_{p,X_{K}^{(p^{\alpha})}}^{\KZ}$ be the unique non-commutative generating series of Coleman functions on $X_{K}$, resp. $X_{K}^{(p^{\alpha})}$, relatively to the chosen determination of $\log_{p}$, which is a horizontal section of $\nabla_{\KZ}$ resp. $\nabla_{\KZ}^{(p^{\alpha})}$ and satisfies the asymptotics $\Li_{p,X_{K}}^{\KZ}(z) \underset{z \rightarrow 0}{\sim} e^{e_{0} \log_{p}(z)}$, resp. 
$\Li_{p,X_{K}^{(p^{\alpha})}}^{\KZ}(z) \underset{z \rightarrow 0}{\sim} e^{e_{0} \log_{p}(z)}$. The coefficients of these formal power series are called $p$-adic multiple polylogarithms (for $N=1$ and depth 1, \cite{Coleman} ; for $N=1$, \cite{Furusho 1} Theorem 3.3 ; for any $N$, \cite{Yamashita}, Proposition 2.6).
\newline\indent For $z \in K$ such that $|z|_{p}<1$, one has the following power series expansion (compare with (\ref{eq: power series mpl}) and see also equation (\ref{eq:multiple harmonic sums})) : 
$$\Li_{p,X_{K}}^{\KZ}[e_{0}^{n_{d}-1}e_{\xi_{d}} \ldots e_{0}^{n_{1}-1}e_{\xi_{1}}](z) = (-1)^{d}\sum_{0<m_{1}<\ldots<m_{d}}
\frac{\big(\frac{\xi_{2}}{\xi_{1}}\big)^{n_{1}} \ldots \big(\frac{\xi_{d}}{\xi_{d-1}}\big)^{n_{d-1}}\big(\frac{z}{\xi_{d}}\big)^{n_{d}}}{m_{1}^{n_{1}}\ldots m_{d}^{n_{d}}} $$
$$ \Li_{p,X_{K}^{(p^{\alpha})}}^{\KZ}[e_{0}^{n_{d}-1}e_{\xi_{d}^{(p^{\alpha})}} \ldots e_{0}^{n_{1}-1}e_{\xi_{1}^{(p^{\alpha})}}](z) = (-1)^{d}\sum_{0<m_{1}<\ldots<m_{d}}
\frac{\big(\frac{\xi_{2}^{p^{\alpha}}}{\xi_{1}^{p^{\alpha}}}\big)^{n_{1}} \ldots \big(\frac{\xi_{d}^{p^{\alpha}}}{\xi_{d-1}^{p^{\alpha}}}\big)^{n_{d-1}}\big(\frac{z}{\xi_{d}^{p^{\alpha}}}\big)^{n_{d}}}{m_{1}^{n_{1}}\ldots m_{d}^{n_{d}}} $$
The expression of $\Li_{p,X_{K}}^{\KZ}$ and $\Li_{p,X_{K}^{(p^{\alpha})}}^{\KZ}$ in terms of Frobenius-invariant paths can be found in \cite{Furusho 2}, Theorem 2.3 for $N=1$. We note that, in this paper, $\Li_{p,X_{K}}^{\KZ}$ and $\Li_{p,X_{K}^{(p^{\alpha})}}^{\KZ}$ are the only objects which depend on a choice of determination of the $p$-adic logarithm.

\subsubsection{Overconvergent variants of $p$-adic multiple polylogarithms associated with the iterated Frobenius}

The overconvergent variants of $p$-adic multiple polylogarithms are defined by means of the images of the canonical De Rham paths by the Frobenius, like the numbers $\zeta_{p,\alpha}$ (Definition \ref{def of p-adic multiple zeta values}), but at variable base-points. This requires to choose an affinoid subspace of $\mathbb{P}^{1,\an}/K$. 

\begin{Notation} \label{the open affine}
Let $U^{\an}$ be the affinoid rigid analytic space $(\mathbb{P}^{1,\an} - \cup_{\xi \in \mu_{N}(K)} B(\xi,1) )/ K$, where $B(\xi,1)$ is the disk of points whose reduction modulo $p$ is $\xi$.
\newline Let $A(U^{\an})$ resp. $A^{\dagger}(U^{\an})$ be the $K$-algebra of global rigid analytic functions resp. overconvergent global rigid analytic functions over $U^{\an}$.  They are Banach $K$-algebras with the norm defined in terms of the power series expansion at $0$ by $||\sum c_{m}z^{m} ||=\sup_{m\in \mathbb{N}} |c_{m}|_{p}$.
\newline Let $F$ be the $\sigma$-linear lift of Frobenius on $A(U^{\an})$ defined by $f(z) \mapsto f(z^{p})$.
\end{Notation}

The space $U^{\an}$ is considered first in \cite{Deligne} \S19.6 for $N=1$, and in \cite{Yamashita} for any $N$, and used later in Furusho's and \"{U}nver's papers on $p$MZV$\mu_{N}$'s. 
\newline\indent In the next definition, the left multiplication by a canonical path below is a convention will be practical for the computations in \S2.1. This definition generalizes particular cases in \cite{Deligne}, \S19.6, \cite{Unver MZV}, \cite{Unver cyclotomic} and \cite{Yamashita}.

\begin{Definition} \label{def of overconvergent}
For $\alpha$ a positive integer, for $z \in U^{\an}(K)$,
$$ \Li_{p,\alpha}^{\dagger}(z) = {}_{\vec{1}_{0}} 1_{z^{p^{\alpha}}}. \phi^{\alpha}({}_z 1_{\vec{1}_{0}}) \in \Pi_{0,0}(K) $$
$$ \Li_{p,-\alpha}^{\dagger}(z) = {}_{\vec{1}_{0}} 1_{z^{p^{\alpha}}}. \phi^{-\alpha}({}_z 1_{\vec{1}_{0}}) \in \Pi^{(p^{\alpha})}_{0,0}(K) $$  
The overconvergent $p$-adic multiple polylogarithm are the following functions, where the $n_{i}$'s are positive integers and the $\xi_{i}$'s are $N$-th roots of unity :
$$ \Li^{\dagger}_{p,\alpha} \big( (n_{i})_{d}; (\xi_{i})_{d} \big) = (-1)^{d} \Li^{\dagger}_{p,\alpha}[e_{0}^{n_{d}-1}e_{\xi_{d}} \ldots e_{0}^{n_{1}-1}e_{\xi_{1}}] \in A^{\dagger}(U^{\an}) $$ 
$$ \Li^{\dagger}_{p,-\alpha}
\big((n_{i})_{d}; (\xi_{i}^{(p^{\alpha})})_{d} \big) = (-1)^{d} \Li^{\dagger}_{p,-\alpha}[e_{0}^{n_{d}-1}e_{\xi_{d}^{(p^{\alpha})}} \ldots e_{0}^{n_{1}-1}e_{\xi_{1}^{(p^{\alpha})}}] \in A^{\dagger}(U^{\an}) $$
\end{Definition}

\begin{Remark} In the case of $\mathbb{P}^{1} \setminus \{0,1,\infty\}$, i.e. $N=1$, the space $U^{\an}$ is equal to $\mathbb{P}^{1,\an} \setminus B(1,1)$, and there exists another natural choice of affinoid subspace, namely $\mathbb{P}^{1,an} \setminus B(\infty,1)$, equal to the analytic unit disk $\mathbb{Z}_{p}^{\an}$. Those two spaces differ by the homography $z \mapsto \frac{z}{z-1}$, which sends $(0,1,\infty)\mapsto (0,\infty,1)$.
\newline The reason for choosing $U^{\an}$ instead of $\mathbb{Z}_{p}^{\an}$ is that working with $\mathbb{Z}_{p}^{\an}$ would require to replace $f(z) \mapsto f(z^{p})$ in Notation \ref{the open affine} by $f(z) \mapsto f(\frac{z^{p}}{z^{p} - (z-1)^{p}})$, i.e. the conjugation of $f(z) \mapsto f(z^{p})$ by $f(z) \mapsto f(\frac{z}{z-1})$, which would make more complicated the computations in the next sections.
\end{Remark}

In the rest of this text, we focus of the objects defined through positive number of iterations of the Frobenius, i.e. $\zeta_{p,\alpha}$ and $\Li_{p,\alpha}^{\dagger}$ with $\alpha$ a a positive integer. One has similar results and proofs for $\zeta_{p,-\alpha}$ and $\Li_{p,-\alpha}^{\dagger}$, which are left to the reader. We will refer to $\zeta_{p,-\alpha}$ and $\Li_{p,-\alpha}^{\dagger}$ in \cite{I-3}, where they will appear naturally.

\section{Overconvergent $p$-adic multiple polylogarithms and iterated integrals}

The formal properties of the Frobenius imply an expression of the overconvergent $p$-adic multiple polylogarithms in terms of $p$-adic cyclotomic multiple zeta values and $p$-adic multiple polylogarithms (Proposition \ref{first formulation}) ; we encode it as an expression of  overconvergent $p$-adic multiple polylogarithms in terms of certain $p$-adic iterated integrals, which is inductive with respect to the weight (Proposition \ref{definition iterated integrals}).

\begin{Notation} \label{notation omega} $\omega_{z_{0}}(z) = \frac{dz}{z-z_{0}}$ for all $z_{0} \in \{0\} \cup \mu_{N}(K)$.
\end{Notation}

\subsection{The differential equation satisfied by overconvergent $p$-adic multiple polylogarithms\label{equation hor}}

By \cite{Deligne}, \S11.11, the Frobenius $\phi$ is characterized by the fact that it commutes with the connections $\nabla_{\KZ}$ and $\nabla^{(p^{\alpha})}_{\KZ}$. In this paragraph, we make this commutation property explicit. According to \cite{Deligne}, \S11.9, this amounts to write a collection of differential equations satisfied by the Frobenius, one for each element of an open affinoid covering of $\mathbb{P}^{1,\an}/K$, and gluing them. For our purposes, it will actually be enough to consider the affinoid subspace $U^{\an}$ of \S1.3, and only the values of the Frobenius at canonical paths, namely, the overconvergent $p$-adic multiple polylogarithms.
\newline\indent We will view the KZ connection as a connection on $\Pi_{0,0}\times X$. Thus, let us write it the Frobenius on $\Pi_{0,0}(K)$, in terms of $p$-adic cyclotomic multiple zeta values. 

\begin{Notation} For any $\xi \in \mu_{N}(K)$, let $f \mapsto f^{(\xi)}$ the map $\Pi_{1,0}(K) \rightarrow \Pi_{\xi,0}(K)$ induced by the functoriality of $\pi_{1}^{\un,\DR}(X_{K})$ with respect to the automorphism $z \mapsto \xi z$ of $X_{K}$.
\end{Notation}

\begin{Lemma} \label{Frobenius at 0,0}
The map $\tau(p^{\alpha})\phi^{\alpha}$ sends $e_{0}\in \Lie(\Pi_{0,0}^{(p^{\alpha})})$ to $e_{0}$ and, for all $\xi \in \mu_{N}(K)$, it sends 
$e_{\xi^{(p^{\alpha})}} \in \Lie(\Pi_{0,0}^{(p^{\alpha})})$ to $\Ad_{\Phi^{(\xi)}_{p,\alpha}}(e_{\xi})=\{\Phi^{(\xi)}_{p,\alpha}\}^{-1}e_{\xi}\Phi^{(\xi)}_{p,\alpha}$.
\end{Lemma}

(We denote $\{\Phi^{(\xi)}_{p,\alpha}\}^{-1}e_{\xi}\Phi^{(\xi)}_{p,\alpha}$ by $\Ad_{\Phi^{(\xi)}_{p,\alpha}}(e_{\xi})$ and not by $\Ad_{\{\Phi^{(\xi)}_{p,\alpha}\}^{-1}}(e_{\xi})$ in order to be coherent with the fact that we read the groupoid multiplication from the right to the left.)

\begin{proof} Let $x \in \{0\} \cup \mu_{N}(K)$ and $T_{x}$ be the tangent space to $\mathbb{P}^{1}$ at $x$. We have $T_{x} \setminus \{0\}\simeq \mathbb{G}_{m}$, and thus by \cite{Deligne}, \S12, $\pi_{1}^{\un,\DR}(T_{x} \setminus \{0\})$ admits a canonical base-point $\omega_{\DR}$ satisfying all the properties of \S1.1 such that we have $\pi_{1}^{\un,\DR}(T_{x} \setminus \{0\},\omega_{\DR})=\Spec(\mathcal{O}^{\sh,\{e_{0}\}})$. It is endowed with the canonical connection on $\pi_{1}^{\un,\DR}(T_{x} \setminus \{0\},\omega_{\DR}) \times T_{x} \setminus \{0\}$ defined by $f \mapsto \frac{dz}{z} e_{x}f$, whose crystalline Frobenius structure is the trivial one given by $e_{x} \mapsto \frac{1}{p}e_{x}$.
\newline\indent By the compatibility between the Frobenius and the tangential base-points (\cite{Deligne}, \S15), for each $\xi \in \mu_{N}(K)$, we have
\begin{equation} \label{eq: Frob at e_xi Lie xi xi} \phi^{\alpha} : e_{\xi^{(p^{\alpha})}} \in \Lie(\Pi^{(p^{\alpha})}_{\xi,\xi}) \mapsto \frac{1}{p^{\alpha}} e_{\xi} \in \Lie(\Pi_{\xi,\xi})
\end{equation}
\indent By the compatibility between canonical paths and the groupoid structure of $\pi_{1}^{\un,\DR}(X_{K})$, we have the following equality, where, in the left-hand side $e_{\xi^{(p^{\alpha})}} \in \Lie(\Pi_{0,0})$ and, in the right-hand side, $e_{\xi^{(p^{\alpha})}} \in \Lie(\Pi_{\xi^{(p^{\alpha})},\xi^{(p^{\alpha})}})$ :
\begin{equation} \label{eq: composition des chemins e xi}e_{\xi^{(p^{\alpha})}} = ({}_0 1_{\xi^{(p^{\alpha})}}) e_{\xi^{(p^{\alpha})}} ({}_{\xi^{(p^{\alpha})}} 1_{0})
\end{equation}
\indent We apply $\tau(p^{\alpha}) \circ \phi^{\alpha}$ to equation (\ref{eq: composition des chemins e xi}) : by the compatibility between the Frobenius with the groupoid structure of $\pi_{1}^{\un,\DR}(X_{K}^{(p^{\alpha})})$, this gives $(\tau(p^{\alpha}) \circ \phi^{\alpha})(e_{\xi}^{(p^{\alpha})}) = (\tau(p^{\alpha}) \circ \phi^{\alpha})({}_0 1_{\xi^{(p^{\alpha})}})\text{ }.\text{ }(\tau(p^{\alpha}) \circ \phi^{\alpha})(e_{\xi^{(p^{\alpha})}})\text{ }.\text{ } (\tau(p^{\alpha}) \circ \phi^{\alpha})({}_{\xi^{(p^{\alpha})}} 1_{0})$ ; by Definition \ref{def of p-adic multiple zeta values} and equation (\ref{eq: Frob at e_xi Lie xi xi}), we obtain the result.
\end{proof}

The above Lemma generalizes a known fact ($N\in \{1,2\}$ and $\alpha=1$, \cite{Deligne Goncharov} \S5.28 ; $N=1$ and $\alpha=-1$, \cite{Unver MZV}, \S4.3 ; any $N$ and $\alpha = \frac{\log(q)}{\log(p)}$, \cite{Yamashita} ; any $N$ and $\alpha=-1$, \cite{Unver cyclotomic}, equation (2.2.6)).
\newline\indent The differential equation which characterizes the Frobenius can be formulated as a differential equation satisfied by $\Li_{p,\alpha}^{\dagger}$, or as a functional relation involving $p$-adic multiple polylogarithms and their overconvergent variants. We use Notation \ref{notation omega}. 

\begin{Proposition} \label{first formulation} We have :
\begin{equation} \label{eq:diff}
d\Li_{p,\alpha}^{\dagger} = \bigg( p^{\alpha}\omega_{0}(z)e_{0} + \sum_{\xi \in \mu_{N}(K)} p^{\alpha} \omega_{\xi}(z) e_{\xi} \bigg) \Li_{p,\alpha}^{\dagger} - \Li_{p,\alpha}^{\dagger} \bigg(  \omega_{0}(z^{p^{\alpha}}) e_{0} + \sum_{\xi \in \mu_{N}(K)} \omega_{\xi^{p^{\alpha}}}(z^{p^{\alpha}}) \{\Phi^{(\xi)}_{p,\alpha}\}^{-1}e_{\xi}\Phi^{(\xi)}_{p,\alpha} \bigg) 
\end{equation} 

\label{second formulation} Equivalently,
\begin{equation} \label{eq:functional} \Li_{p,\alpha}^{\dagger}(z) = 
\Li_{p,X_{K}}^{\KZ}(z)\big((p^{\alpha}e_{x})_{x \in \{0\} \cup \mu_{N}(K)}) \text{ } \text{ }
\Li_{p,X_{K}^{(p^{\alpha})}}^{\KZ}(z^{p^{\alpha}})\big(e_{0},(\{\Phi^{(\xi)}_{p,\alpha}\}^{-1}e_{\xi}\Phi^{(\xi)}_{p,\alpha})_{\xi \in \mu_{N}(K)}) \big)^{-1}
\end{equation}
\end{Proposition}

\begin{proof} \label{prop horizontality diagram}
By \cite{Deligne}, \S7.30.2$_{S}$, we view the KZ connection as the connection on the trivial bundle $\Pi_{0,0}\times X_{K}$ with values in $\Lie \Pi_{0,0} \otimes \Omega^{1}(X_{K})$ given by
\begin{equation} \label{eq:KZ Lie 0,0} f \mapsto f^{-1}\bigg(df - \big( e_{0}  \omega_{0} + \sum\limits_{\xi \in \mu_{N}(K)} e_{\xi} \omega_{\xi}\big) f \bigg)
\end{equation}
\indent By \cite{Deligne} \S11.11, the Frobenius $\phi$ commutes with the connections, which gives the commutativity of the following diagram (where $F$ is as in Notation \ref{the open affine})
\begin{equation} \label{eq:commutativity Frob KZ}
\begin{array}{ccccc}
(F^{\alpha})^{\ast}\big(\Pi_{0,0}^{(p^{\alpha})}(K) \times A(U^{\an})\big)
& \overset{\phi^{\alpha}}{\xrightarrow{\hspace*{1.5cm}}}
& \Pi_{0,0}(K) \times A(U^{\an})
\\
\Big\downarrow{(F^{\alpha})^{\ast}\nabla_{\KZ}^{(p^{\alpha})}} &&
\Big\downarrow{\footnotesize{\nabla_{\KZ}}} 
\\
(F^{\alpha})^{\ast}\big(\Lie(\Pi_{0,0}^{(p^{\alpha})}(K)) \otimes \Omega^{1}(X_{K}^{\an}) \big)
& \overset{\phi^{\alpha}  \otimes \id}{\xrightarrow{\hspace*{1.5cm}}}
& 
\Lie \Pi_{0,0}(K) \otimes \Omega^{1}(X_{K}^{\an})
\end{array}
\end{equation}
\indent We apply the arrows of (\ref{eq:commutativity Frob KZ}) to the canonical De Rham path ${}_{z}1_{\vec{1}_{0}}$ :
\newline \indent $(a)$ The canonical De Rham path ${}_{z}1_{\vec{1}_{0}}$ is mapped by $(F^{\alpha})^{\ast}\nabla_{\KZ}^{(p^{\alpha})}$ to $-\omega_{0}(z^{p^{\alpha}})e_{0} - \sum_{\xi \in \mu_{N}(K)} \omega_{\xi}(z^{p^{\alpha}}) e_{\xi^{(p^{\alpha})}}$ (by equation (\ref{eq:KZ Lie 0,0})), which is mapped by  $\phi^{\alpha} \otimes \id$ to $\displaystyle\frac{1}{p^{\alpha}}\omega_{0}(z)e_{0}  - \sum_{\xi \in \mu_{N}(K)}  \omega_{\xi_{j}}(z)\frac{1}{p^{\alpha}} \Ad_{\tau(p^{-\alpha})\Phi^{(\xi)}_{p,\alpha}}(e_{\xi})$ (by Lemma \ref{Frobenius at 0,0}).
\newline\indent $(b)$ 
The canonical De Rham path ${}_{z}1_{\vec{1}_{0}}$ is mapped by $\phi^{\alpha}$ to  $\tau(p^{-\alpha})\Li_{p,\alpha}^{\dagger}$ (by Definition \ref{def of overconvergent}), which is mapped by $\nabla_{\KZ}$ to 
$\displaystyle(\tau(p^{-\alpha})\Li_{p,\alpha}^{\dagger})^{-1} \bigg( d\tau(p^{-\alpha})\Li_{p,\alpha}^{\dagger} - \big( \omega_{0}(z)e_{0} + \sum_{\xi \in \mu_{N}(K)} \omega_{\xi}(z) e_{\xi_{j}} \big) \tau(p^{-\alpha})\Li_{p,\alpha}^{\dagger} \bigg)$ (by equation (\ref{eq:KZ Lie 0,0})).
\newline\indent The results of (a) and (b) above are the same by the commutativity of (\ref{eq:commutativity Frob KZ}) ; applying $\tau(p^{\alpha})$ to the equality between them, we obtain equation (\ref{eq:diff}).
\newline\indent Now, by the definition of $p$-adic multiple polylogarithms (\S1.3.1) we have $\nabla_{\KZ}\Li_{p,X_{K}}^{\KZ}=0$ and $\nabla_{\KZ}^{(p^{\alpha})} \Li_{p,X_{K}^{(p^{\alpha})}}^{\KZ}=0$, whence the right-hand side of (\ref{eq:functional}) is, like $\Li_{p,\alpha}^{\dagger}$, a solution to the following differential equation on $L$ :
\begin{equation} \label{eq:diff L} dL = \bigg( p^{\alpha}\omega_{0}(z)e_{0} + \sum_{\xi \in \mu_{N}(K)} p^{\alpha} \omega_{\xi}(z) e_{\xi} \bigg)L -L\bigg(  \omega_{0}(z^{p^{\alpha}}) e_{0} + \sum_{\xi \in \mu_{N}(K)} \omega_{\xi^{p^{\alpha}}}(z^{p^{\alpha}}) \Ad_{\Phi^{(\xi)}_{p,\alpha}} (e_{\xi}) \bigg) 
\end{equation}
\newline\indent By Definition \ref{def of overconvergent}, we have $\Li_{p,\alpha}^{\dagger}(0)=1$, and by the definition of $\Li_{p,X_{K}}^{\KZ}$ and $\Li_{p,X_{K}^{(p^{\alpha})}}^{\KZ}$ (\S1.3.1) the right-hand side of (\ref{eq:functional}) is equivalent when $z \rightarrow 0$ to $e^{p^{\alpha} \log_{p}(z)e_{0}} e^{- \log_{p}(z^{p^{\alpha}})e_{0}}=1$, where we have used $\log(z^{p^{\alpha}}) = p^{\alpha}\log_{p}(z)$. The equation (\ref{eq:diff L}) being pro-unipotent, it has a unique solution $L\in K[[z]]\langle\langle e_{0\cup \mu_{N}}\rangle\rangle$ such that $L(0)=1$ (for any word $w$, equation (\ref{eq:diff L}) plus $L(0)=1$ determine $L[w]$ by induction on the weight of $w$).
\end{proof}

Variants of equation (\ref{eq:diff}) for $N=1$, $\alpha=-1$ are written, \cite{Deligne}, \S19.6, equation (19.6.2), \cite{Unver MZV} \S5.2 Proposition 1 ; any $N$, $\alpha=-1$, \cite{Unver cyclotomic} equation (2.2.9) ; the particular case for any $N$, $\alpha=\frac{\log(q)}{\log(p)}$, \cite{Yamashita} Remark 3.2.
Equation (\ref{eq:functional}) is written in the case $N=1$, $\alpha=1$, in \cite{Furusho 2} Theorem 2.14.

\subsection{Decomposition of overconvergent $p$-adic multiple polylogarithms in terms of iterated integrals\label{first consequences}}

The goal of this paragraph is to rewrite (\ref{eq:diff}), as a "decomposition" of $\Li_{p,\alpha}^{\dagger}$ in terms of certain iterated integrals. It arises initially as given by an induction on the weight. We will turn this induction into an induction on the depth in \S4.

\begin{Definition} \label{def alphabet} \label{the big alphabet} Let $e_{0\cup\mu_{N} \cup \mu_{N}^{(p^{\alpha})}}$ be the alphabet $\{e_{0},e_{\xi_{1}},\ldots,e_{\xi^{N}},e_{\xi^{(p^{\alpha})}},\ldots,e_{{\xi^{N}}^{(p^{\alpha})}}\}$.
\newline The weight, resp. depth of a word over the alphabet $e_{0\cup \mu_{N} \cup \mu_{N}^{(p^{\alpha})}}$ is its number of letters, resp. its number of letters distinct from $e_{0}$.
\end{Definition}

Identifying the alphabet $e_{0\cup\mu_{N} \cup \mu_{N}^{(p^{\alpha})}}$ to $\{p^{\alpha}\omega_{0}(z)\} \cup \{ p^{\alpha}\omega_{\xi}(z) \text{ }|\text{ }\xi \in \mu_{N}(K) \} \cup  \{ \omega_{\xi^{p^{\alpha}}}(z^{p^{\alpha}})\text{ }|\text{ }\xi \in \mu_{N}(K)\}$ in the natural way, we regard the right-hand side of equation (\ref{eq:functional}) as a map $\mathcal{O}^{\sh,e_{0\cup\mu_{N} \cup \mu_{N}^{(p^{\alpha})}}} \rightarrow K[[z]]$, which we are going to factorise by certain iterated integrals (Definition \ref{definition iterated integrals}) and a certain map of "decomposition" (Definition \ref{definition of dec}). In order to define the map of decomposition, let us write explicitly the coefficients of equation (\ref{eq:diff}) : for any $w=e_{0}^{n_{d}-1}e_{\xi_{d}} \ldots e_{0}^{n_{1}-1}e_{\xi_{1}}e_{0}^{n_{0}-1}$, where the $n_{i}$'s are positive integers and the $\xi_{i}$'s are $N$-th roots of unity, we have

\begin{multline} \label{eq:diff coefficient}
d\Li^{\dagger}_{p,\alpha}[e_{0}^{n_{d}-1}e_{\xi_{d}} \ldots e_{0}^{n_{1}-1}e_{\xi_{1}}e_{0}^{n_{0}-1}] =
\sum_{\substack{w_{1},w_{2}\text{ words } \\ w_{1}w_{2}=w\\ 
 		w_{2}\not=\emptyset}} \sum_{\xi \in \mu_{N}(K)}
\bigg(\Ad_{\Phi_{p,\alpha}^{(\xi)}}(e_{\xi})[w_{2}]\bigg) \omega_{\xi^{p^{\alpha}}}(z^{p^{\alpha}})\Li^{\dagger}_{p,\alpha}[w_{1}] \text{ }+
\\ 
\left\{ 
\begin{array}{ll} \displaystyle
p^{\alpha}\omega_{0}(z)\Li^{\dagger}_{p,\alpha}[e_{0}^{n_{d}-2}e_{\xi_{d}} \ldots e_{0}^{n_{1}-1}e_{\xi_{1}}e_{0}^{n_{0}-1}] - p^{\alpha}\omega_{0}(z) \Li^{\dagger}_{p,\alpha}[e_{0}^{n_{d}-1}e_{\xi_{d}} \ldots e_{0}^{n_{1}-1}e_{\xi_{1}}e_{0}^{n_{0}-2}]
&  \text{ if }n_{d}\geqslant 2,\text{ } n_{0}\geqslant 2
\\ \displaystyle
p^{\alpha}\omega_{0}(z)\Li^{\dagger}_{p,\alpha}[e_{0}^{n_{d}-2}e_{\xi_{d}} \ldots e_{0}^{n_{1}-1}e_{\xi_{1}}]
& \text{ if }n_{d}\geqslant 2,\text{ } n_{0}=1
\\ \displaystyle
p^{\alpha}\omega_{\xi_{d}}(z)\Li^{\dagger}_{p,\alpha}[e_{0}^{n_{d-1}-1}e_{\xi_{d-1}}\ldots e_{\xi_{1}} e_{0}^{n_{0}-1}] - p^{\alpha}\omega_{0}(z)
\Li_{p,\alpha}^{\dagger}[e_{\xi_{d}}e_{0}^{n_{d-1}-1}e_{\xi_{d-1}} \ldots e_{0}^{n_{1}-1}e_{\xi_{1}}e_{0}^{n_{0}-2}]
& \text{ if }n_{d}= 1,\text{ }n_{0}\geqslant 2
\\ \displaystyle p^{\alpha}\omega_{\xi_{d}}(z)\Li_{p,\alpha}^{\dagger}[e_{0}^{n_{d-1}-1}e_{\xi_{d-1}}\ldots e_{0}^{n_{1}-1}e_{\xi_{1}}]
 & \text{ if }n_{d}=1,\text{ }n_{0}=1
 \\
\end{array}	\right.
\end{multline}
 
where the restriction to $w_{2}\not=\emptyset$ in the first line of (\ref{eq:diff coefficient}) above is due to the fact that the term indexed by $(w_{1},w_{2})=(w,\emptyset)$ is $0$ : the coefficient of the empty word in  $\Ad_{\Phi_{p,\alpha}^{(\xi)}}(e_{\xi})$ is zero because $\Ad_{\Phi_{p,\alpha}^{(\xi)}}(e_{\xi})$ is obtained by multiplying two non-commutative formal power series $e_{\xi}$, which is of weight 1.

\begin{Definition} \label{definition of dec}
	Let $\dec_{p,\alpha} : \mathcal{O}^{\sh,e_{0\cup \mu_{N}}} \rightarrow\mathcal{O}^{\sh,e_{0\cup \mu_{N}\cup \mu_{N}^{(p^{\alpha})}}} \otimes_{\mathbb{Q}} K$ be defined by induction on the weight, as follows : $\dec_{p,\alpha}(\emptyset)=1$ ; for all $n \in \mathbb{N}^{\ast}$, 
$\dec_{p,\alpha}(e_{0}^{n})=e_{0}\dec_{p,\alpha}(e_{0}^{n-1})-e_{0}\dec_{p,\alpha}(e_{0}^{n-1})=0$, and for any $w= e_{0}^{n_{d}-1}e_{\xi_{d}} \ldots e_{0}^{n_{1}-1}e_{\xi_{1}}e_{0}^{n_{0}-1}$ word over $e_{0\cup\mu_{N}}$, where the $n_{i}$'s are positive integers and the $\xi_{i}$'s are $N$-th roots of unity $(d\geqslant 1)$,

\begin{multline} \label{eq:definition dec}
\dec_{p,\alpha}(e_{0}^{n_{d}-1}e_{\xi_{d}} \ldots e_{0}^{n_{1}-1}e_{\xi_{1}}e_{0}^{n_{0}-1}) = \sum_{\substack{w_{1},w_{2}\text{ words } \\ w_{1}w_{2}=w\\ 
w_{2}\not=\emptyset \\ \depth(w_{1})\geqslant 1\text{ or }w_{1}=\emptyset}} \sum_{\xi \in \mu_{N}(K)} \bigg( \Ad_{\Phi_{p,\alpha}^{(\xi)}}(e_{\xi})[w_{2}] \bigg) e_{\xi^{(p^{\alpha})}}\dec_{p,\alpha}(w_{1}) +
\\ 
\left\{ 
\begin{array}{ll} 
e_{0}\dec_{p,\alpha}(e_{0}^{n_{d}-2}e_{\xi_{d}} \ldots e_{0}^{n_{1}-1}e_{\xi_{1}}e_{0}^{n_{0}-1}) - e_{0} \dec_{p,\alpha}(e_{0}^{n_{d}-1}e_{\xi_{d}} \ldots e_{0}^{n_{1}-1}e_{\xi_{1}}e_{0}^{n_{0}-2}) 
& \text{ if }n_{d}\geqslant 2, n_{0}\geqslant 2
\\
e_{0}\dec_{p,\alpha}(e_{0}^{n_{d}-2}e_{\xi_{d}} \ldots e_{0}^{n_{1}-1}e_{\xi_{1}}) 
& \text{ if }n_{d}\geqslant 2, n_{0}=1
\\
e_{\xi_{d}}\dec_{p,\alpha}(e_{0}^{n_{d-1}-1}e_{\xi_{d}}\ldots e_{\xi_{1}}e_{0}^{n_{0}-1}) - e_{0} \dec_{p,\alpha}(e_{\xi_{d}}e_{0}^{n_{d-1}-1}e_{\xi_{d-1}} \ldots e_{0}^{n_{1}-1}e_{\xi_{1}}e_{0}^{n_{0}-2})
& \text{ if }n_{d}= 1, n_{0}\geqslant 2
\\e_{\xi_{d}}\dec_{p,\alpha}(e_{0}^{n_{d-1}-1}e_{\xi_{d-1}}\ldots e_{\xi_{1}})
& \text{ if }n_{d}=1, n_{0}=1
\\
 \end{array}	 \right.
\end{multline}
\end{Definition}

Given that $\Li_{p,\alpha}^{\dagger}$ has no pole at $0$ by Definition \ref{def of overconvergent}, we can consider only iterated integrals which have no pole at $0$ :

\begin{Definition} \label{definition iterated integrals}Let $\mathcal{O}^{\sh,e_{0\cup \mu_{N}\cup \mu_{N}^{(p^{\alpha})}}}_{\conv}\subset\mathcal{O}^{\sh,e_{0\cup \mu_{N}\cup \mu_{N}^{(p^{\alpha})}}}$ be the vector subspace generated by words whose rightmost letter is not $e_{0}$. Let $$\displaystyle\int_{p,\alpha} : \mathcal{O}^{\sh,e_{0\cup \mu_{N}\cup \mu_{N}^{(p^{\alpha})}}}_{\conv} \rightarrow K[[z]] $$
be the linear map defined by induction by $\int \emptyset = 1$ and the following variant of the KZ equation (\ref{eq:KZ equation}) :
\begin{equation} \label{eq:equation KZ with congruences} d \int_{p,\alpha} e_{0}^{n_{d}-1}e_{{\xi_{d}}^{u_{d}}} \ldots e_{0}^{n_{1}-1}e_{{\xi_{1}}^{u_{1}}} 
= \left\{ \begin{array}{ll} \displaystyle p^{\alpha}\frac{dz}{z} \int_{p,\alpha} e_{0}^{n_{d}-2}e_{{\xi_{d}}^{u_{d}}} \ldots e_{0}^{n_{1}-1}e_{{\xi_{1}}^{u_{1}}} &\text{ if } n_{d}\geqslant 2 
\\ \displaystyle p^{\alpha}\frac{dz}{z-\xi_{d}} \int_{p,\alpha} e_{0}^{n_{d-1}-1}e_{{\xi_{d-1}}^{u_{d-1}}} \ldots e_{0}^{n_{1}-1}e_{{\xi_{1}}^{u_{1}}} &\text{ if } n_{d}=1 \text{ and }u_{d}=1  
\\\displaystyle \frac{d(z^{p^{\alpha}})}{z^{p^{\alpha}}-\xi_{d}^{p^{\alpha}}} \int_{p,\alpha} e_{0}^{n_{d-1}-1}e_{{\xi_{d-1}}^{u_{d-1}}} \ldots e_{0}^{n_{1}-1}e_{{\xi_{1}}^{u_{1}}} &\text{ if } n_{d}=1 \text{ and }u_{d}=(p^{\alpha}) 
\end{array}\right. 
\end{equation}
and the condition that $\displaystyle\int_{p,\alpha} e_{0}^{n_{d}-1}e_{{\xi_{d}}^{u_{d}}} \ldots e_{0}^{n_{1}-1}e_{{\xi_{1}}^{u_{1}}}$ has constant coefficient equal to $0$.
\end{Definition}

Explicitly, these iterated integrals are given as the power series expansions, which are convergent on $\{z \in K\text{ }|\text{ }|z|_{p}<1\}$ ; this is a variant of (\ref{eq: power series mpl}) :

\begin{equation} \label{eq:power series it int}  \int_{p,\alpha} e_{0}^{n_{d}-1}e_{{\xi_{d}}^{u_{d}}} \ldots e_{0}^{n_{1}-1}e_{{\xi_{1}}^{u_{1}}} = (p^{\alpha})^{n_{1}+\cdots+n_{d}} \sum_{\substack{0=m_{0}<m_{1}<\ldots<m_{d} \\ 
		\forall i \text{ s.t.} u_{i} = (p^{\alpha}),\text{ }m_{i-1} \equiv m_{i}[p^{\alpha}]}}
\frac{(\frac{{\xi_{2}}}{{\xi_{1}}})^{n_{1}} \ldots (\frac{z}{{\xi_{d}}})^{n_{d}}}{m_{1}^{n_{1}}\ldots m_{d}^{n_{d}}} \in K[[z]]
\end{equation}

We can now reformulate equation (\ref{eq:diff}) ; below we view $\Li_{p,\alpha}^{\dagger}$ as a map $\mathcal{O}^{\sh,e_{0\cup \mu_{N}}} \mapsto K[[z]]$, which is enabled by the fact that each $\Li_{p,\alpha}^{\dagger}[w]$, being an element of $A^{\dagger}(U^{\an})$, is characterized by its power series expansion at $0$.

\begin{Proposition} \label{prop decomposition}The map $\dec_{p,\alpha}$ take values in $\mathcal{O}^{\sh,e_{0\cup \mu_{N}\cup \mu_{N}^{(p^{\alpha})}}}_{\conv}$ and we have
\begin{equation} \label{eq:resolution} 
\Li_{p,\alpha}^{\dagger} = \int_{p,\alpha} \text{ }\circ\text{ }\dec_{p,\alpha}
\end{equation}
\end{Proposition}

\begin{proof} The differential equation (\ref{eq:diff coefficient}) together with $\Li_{p,\alpha}^{\dagger}(z=0)=1$ determine $\Li_{p,\alpha}^{\dagger}[w]$ for all words $w$  by induction on the weight of $w$. By equations (\ref{eq:definition dec}) and (\ref{eq:equation KZ with congruences}), $\displaystyle\int_{p,\alpha} \text{ }\circ\text{ }\dec_{p,\alpha}$ satisfies the same differential equation than (\ref{eq:diff coefficient}). Moreover, by Definition \ref{eq:definition dec} and Definition \ref{eq:power series it int}, for any word $w$, $\displaystyle\bigg(\int_{p,\alpha}\circ \dec_{p,\alpha}\bigg)(w)(z=0)$ is equal to $0$ if $w$ is non-empty and 1 if $w$ is the empty word. This concludes the proof.
\end{proof}

\begin{Remark} One can also write another formula for $\dec_{p,\alpha}$, by using operations in the shuffle Hopf algebra $\mathcal{O}^{\sh,e_{0 \cup \mu_{N}\cup \mu_{N}^{(p^{\alpha})}}}$ arising as duals of the product, inversion and composition of non-commutative formal power series which appear in equation (\ref{eq:functional}).
\end{Remark}

\section{Regularized $p$-adic iterated integrals and bounds on their norms\label{iterated integrals}}

Let the affinoid space $U^{\an}=\mathbb{P}^{1,\an} \setminus \cup_{\xi^{N}=1} B(\xi,1)$ over $K$, and $A(U^{\an})$ its $K$-algebra of global rigid analytic functions (Notation \ref{the open affine}). We define in \S3.2 regularized $p$-adic iterated integrals (Definition \ref{def regularization integrals}), which lie in a specific subalgebra of $\{f \in A(U^{\an}) \text{ }|\text{ }f(0)=0 \}$ defined in \S3.1 (Definition \ref{def of the subalgebra}); we show that they can be computed by induction on the depth (\S3.3, Proposition \ref{prop comp reg it int}) and we prove bounds on their norms (\S3.4, Proposition \ref{bound on valuation}).

\subsection{A subalgebra of the algebra of rigid analytic functions on a subspace of $\mathbb{P}^{1,\an}/K$\label{subspace}}

By Proposition \ref{prop Mahler general}, the map "coefficients at $0$", $\Cf_{0} : \bigg( f = \sum\limits_{m\in \mathbb{N}^{\ast}} c_{m}z^{m}\bigg) \mapsto \bigg( m \in \mathbb{N}^{\ast} \mapsto c_{m} \in K\bigg) $ induces an isometric isomorphism of Banach spaces over $K$, $\{f \in A(U^{\an}) \text{ }|\text{ }f(0)=0 \} \rightarrow \mathcal{C}(\mathbb{Z}_{p}^{(N)},K)$
where the target is the space of continuous functions on $\mathbb{Z}_{p}^{(N)}=\varprojlim \mathbb{Z}/Np^{l}\mathbb{Z} \simeq \mathbb{Z}/N\mathbb{Z} \times \mathbb{Z}_{p}$, with values in $K$, endowed with the supremum norm $||\text{ }||_{\infty}$, and where we consider the canonical inclusion $\mathbb{N}^{\ast}=\amalg_{1\leqslant m_{0}\leqslant N} (m_{0}+N\mathbb{N}) \subset \mathbb{Z}_{p}^{(N)}$.
\newline\indent We define a subspace of $\mathcal{C}(\mathbb{Z}^{(N)}_{p},K)$ and prove that it is the image by $\Cf_{0}$ of a subalgebra of $\{f \in A(U^{\an})\text{ }|\text{ }f(0)=0\}$ (Proposition \ref{multiplication}).

\begin{Definition} \label{def of the subalgebra}(i)  Let $\mathcal{S}_{\alpha}\subset K^{\mathbb{N}}$ be the subspace of the sequences 
$b=(b_{l})_{l \in \mathbb{N}}$ such that there exists a polynomial $P$ with coefficients in $\mathbb{R}$ such that we have, for all $l\in \mathbb{N}$, 
$|b_{l}|_{p} \leqslant P(l) p^{(\alpha-1)l}$.
\newline (ii) Let $\LA_{\mathcal{S}_{\alpha}}(\mathbb{Z}_{p}^{(N)},K)$ be the space of functions $c : \mathbb{Z}_{p}^{(N)} \rightarrow K$, such that, for all $r_{0} \in \{0,1,\ldots,p^{\alpha}-1\}$ and $\xi \in \mu_{N}(K)$, there exist sequences $(c^{(l,\xi)}(r_{0}))_{l\in\mathbb{N}} \in \mathcal{S}_{\alpha}$, such that, for $m \in \mathbb{N}^{\ast}$ such that $|m - r|_{p} \leqslant p^{-\alpha}$ we have the absolutely convergent series expansion :
\begin{equation} \label{eq:condition series expansion} c(m) = \sum_{l \in \mathbb{N}} \sum_{\xi \in \mu_{N}(K)} c^{(l,\xi)}(r_{0}) \xi^{-m}(m-r_{0})^{l}
\end{equation}
\end{Definition}

Any element of $\LA_{\mathcal{S}_{\alpha}}(\mathbb{Z}_{p}^{(N)},K)$ is locally analytic, by the fact that $\mathbb{N}^{\ast}$ is dense in $\mathbb{Z}_{p}^{(N)}$ and because for any $\xi \in \mu_{N}(K)$, the map $m \mapsto \xi^{m}$ is constant on each class of congruence modulo $N$. We have a natural isomorphism of vector spaces $\LA_{\mathcal{S}_{\alpha},\xi}(\mathbb{Z}_{p},K)\simeq \mathcal{S}_{\alpha}^{Np^{\alpha}}$ defined by 
$c \mapsto \bigg( (c^{(l,\xi)}(r_{0}))_{\l \in \mathbb{N}} \bigg)_{\substack{\xi\in\mu_{N}(K)\\ r_{0}) \in\{0,\ldots,p^{\alpha}-1\}}}$. In the $\alpha=1$ case, the space $\LA_{\mathcal{S}_{\alpha}}(\mathbb{Z}_{p}^{(N)},K)$ is a subspace of the space of "$N$-power series functions" introduced in \cite{Unver cyclotomic}, Definition 5.0.2. (We add the condition that the sequences of coefficients in $\mathcal{S}_{\alpha}$, and another difference is that in \cite{Unver cyclotomic} one considers $N$ unrelated power series expansions on each $m + N\mathbb{N}$, $m \in \{0,\cdots,N-1\}$, unrelated to each other. The space of "$N$-power series functions" does not enable us to prove our results.)

\begin{Example} \label{example elements of LA}For any $\xi_{0} \in \mu_{N}(K)$, $l \in \mathbb{N}$ and $\xi \in \mu_{N}(K)$,
\newline (i) $\Cf_{0}\bigg( \displaystyle\frac{p^{\alpha}z}{z-\xi_{0} }=-\sum_{m\in\mathbb{N}^{\ast}} \big(\frac{p^{\alpha}z}{\xi_{0} }\big)^{m}\bigg) \in \LA_{\mathcal{S}_{\alpha}}(\mathbb{Z}_{p}^{(N)},K)$ and  $\displaystyle(\Cf_{0}\frac{z}{z-\xi_{0} })^{(l,\xi)}(r_{0})=\left\{\begin{array}{ll}-p^{\alpha} &\text{ if }(l,\xi)=(0,\xi_{0})
\\ 0 & \text{ otherwise }\end{array}\right.$
\newline (ii) $\displaystyle\Cf_{0}\bigg(\displaystyle\frac{z^{p^{\alpha}}}{z^{p^{\alpha}}-\xi_{0}^{p^{\alpha}}}= -\sum_{m\in\mathbb{N}^{\ast}} \big(\frac{z}{\xi_{0}}\big)^{p^{\alpha}m}\bigg)  \in  \LA_{\mathcal{S}_{\alpha}}(\mathbb{Z}_{p}^{(N)},K)$ and $(\Cf_{0}\frac{z^{p^{\alpha}}}{z^{p^{\alpha}}-\xi_{0}^{p^{\alpha}}} )^{(l,\xi)}(r_{0})=\left\{\begin{array}{ll}-p^{\alpha} &\text{ if }(l,\xi,r)=(0,\xi_{0},0)
	\\ 0 & \text{ otherwise }\end{array}\right.$.
\end{Example}

\begin{Lemma} \label{def of coeff B}Let $l,l' \in \mathbb{N}$, $\xi \in \mu_{N}(K)$.
\newline (i) 
The map sending $u \in \mathbb{N}^{\ast}$ to, respectively, 
$\displaystyle\sum_{u_{1}=1}^{u-1} u_{1}^{l}(u-u_{1})^{l'} \xi^{u_{1}}$,
$\displaystyle\sum_{u_{1}=1}^{u} u_{1}^{l}(u-u_{1})^{l'} \xi^{u_{1}}$,
$\displaystyle\sum_{u_{1}=0}^{u-1} u_{1}^{l}(u-1-u_{1})^{l'} \xi^{u_{1}}$,
$\displaystyle\sum_{u_{1}=0}^{u} u_{1}^{l}(u-u_{1})^{l'} \xi^{u_{1}}$ is a $K$-linear combination of the functions $u \mapsto u^{l_{0}} \xi^{u}$, $l_{0} \in \{0,\cdots,l+l'+1\}$.
\newline(ii) If we denote, respectively, by  $ \mathcal{B}^{(l,l',\xi)}_{l_{0}}$, $ {\mathcal{B}^{\ast}}^{(l,l',\xi)}_{l_{0}}$
$\widetilde{\mathcal{B}^{\ast}}^{(l,l',\xi)}_{l_{0}}$
$\widetilde{\widetilde{\mathcal{B}^{\ast}}}^{(l,l',\xi)}_{l_{0}}$ the coefficients of $u^{l_{0}} \xi^{u}$ in the above expression, then we have
$$ \min\bigg( 
v_{p}\big(\mathcal{B}_{l_{0}}^{(l,l',\xi)}\big), v_{p}\big({\mathcal{B}^{\ast}}_{l_{0}}^{(l,l',\xi)}\big),
v_{p}\big(\widetilde{\mathcal{B}^{\ast}}_{l_{0}}^{(l,l',\xi)}\big), v_{p}\big(\widetilde{\widetilde{\mathcal{B}^{\ast}}}_{l_{0}}^{(l,l',\xi)}\big) \bigg) \geqslant  - 1 - \frac{\log(1+l+l')}{\log(p)} $$
\end{Lemma}

\begin{proof} Let us prove the result for $\displaystyle\sum_{u_{1}=1}^{u-1} u_{1}^{l}(u-u_{1})^{l'} \xi^{u_{1}}$, the other cases being similar. By writing $(u-u_{1})^{l'}=\sum_{l''=0}^{l'} {l' \choose l''} u^{l''}u_{1}^{l'-l''}$, we are reduced to the $l'=0$ case. We distinguish two cases : 
\newline - If $\xi=1$, we can write $\sum_{u_{1}=1}^{u-1} u_{1}^{l} = \sum_{l_{0}=1}^{l+1} \frac{1}{l+1} {l+1 \choose l_{0}} B_{l+1-l_{0}} u^{l_{0}}-\delta_{l,0}$ with $B$ denoting Bernoulli numbers. By Von Staudt-Clausen's theorem, we have $v_{p}(B_{l+1-l_{0}})\geqslant  -1$ and, given that $p^{v_{p}(l+1)} \leqslant l+1$ thus $v_{p}(\frac{1}{l+1})=-v_{p}(l+1) \geqslant  -\frac{\log(l+1)}{\log(p)}$, we have $v_{p}\big(  \frac{1}{l+1} {l+1 \choose l_{0}} B_{l+1-l_{0}}\big)\geq -1 - \frac{\log(l+1)}{\log(p)}$.
\newline - If $\xi_{1} \not=1$, we consider the equation $\sum_{u_{1}=1}^{u-1} T^{u_{1}} = \frac{1-T^{u_{1}}}{1-T}-1$, where $T$ is a formal variable, we apply to it $(T\frac{d}{dT})^{l}$, and we substitute $\xi$ to $T$ : this gives, by induction on $l$, the existence of a polynomial $P_{l}$, with coefficients in $\mathbb{Z}[\xi,\frac{1}{\xi},\frac{1}{\xi-1}] \subset \mathbb{Z}_{p}$, of degree $\leqslant l$, such that 
$\sum_{u_{1}=1}^{u-1} u_{1}^{l} \xi^{u_{1}} =P(l)(u_{1}) \xi^{u_{1}}$.
\end{proof}

We now prove that $\Cf_{0}^{-1}(\LA_{\mathcal{S}_{\alpha}}(\mathbb{Z}_{p}^{(N)},K))$ (where $\Cf_{0}$ is the "map coefficients of the power series expansion at $0$", see Proposition \ref{prop Mahler general}) is a subalgebra of $A(U^{\an})$, with an explicit formula for the multiplication.

\begin{Proposition} \label{multiplication} Let  $f_{1}, f_{2} \in A(U^{\an})$ with $f_{1}(0)=f_{2}(0)=0$, such that $c_{1}=\Cf_{0}(f_{1})$ and $c_{2}=\Cf_{0}(f_{2})$ are elements of $\LA_{\mathcal{S}_{\alpha}}(\mathbb{Z}_{p}^{(N)},K)$. Then the map $c=\Cf_{0}(f_{1}f_{2})$ is the element of $\LA_{\mathcal{S}_{\alpha}}(\mathbb{Z}_{p}^{(N)},K)$ determined as follows : for all $L \in \mathbb{N}$ and $\xi \in \mu_{N}(K)$,
\begin{multline} \label{eq:product at 0}
c^{(L,\xi)}(0) =
\\
\sum_{\substack{l,l' \in \mathbb{N} \\ l+l'+1\geqslant  L}} (p^{\alpha})^{l+l'-L} \sum_{\xi' \in \mu_{N}(K)}
\bigg( \mathcal{B}_{L}^{l,l',(\frac{\xi'}{\xi})^{p^{\alpha}}} \text{ }c_{1}^{(l,\xi)}(0)\text{ }c_{2}^{(l',\xi')}(0) 
+ 
\sum_{r=1}^{p^{\alpha}-1} \widetilde{\widetilde{\mathcal{B}^{\ast}}}_{L}^{l,l',(\frac{\xi'}{\xi})^{p^{\alpha}}} (\xi')^{r} \text{ } c_{1}^{(l,\xi)}(r)\text{ } c_{2}^{(l',\xi')}(p^{\alpha}-r) \bigg)
\end{multline}
 and, for all $r_{0} \in \{1,\ldots,p^{\alpha}-1\}$.
\begin{multline} \label{eq:product at r0}
c^{(L,\xi)}(r_{0})= 
\sum_{\substack{l,l' \in \mathbb{N} \\ l+l'+1\geqslant  L}} (p^{\alpha})^{l+l'-L} \sum_{\xi'\in \mu_{N}(K)} \bigg(  
\big(\frac{\xi'}{\xi}\big)^{r_{0}}{\mathcal{B}^{\ast}}_{L}^{l,l',(\frac{\xi'}{\xi})^{p^{\alpha}}}\text{ } c_{1}^{(l,\xi)}(0) \text{ }c_{2}^{l',\xi'}(r_{0}) +
\\ 
\sum_{r=1}^{r_{0}-1}  \widetilde{\mathcal{B}^{\ast}}_{L}^{l,l',(\frac{\xi}{\xi'})^{p^{\alpha}}} \bigg(\frac{\xi}{\xi'}\bigg)^{(r_{0}-r)}  \text{ } c_{1}^{(l,\xi)}(r) \text{ }c_{2}^{(l',\xi')}(r_{0}-r)
+ \sum_{r=r_{0}}^{p^{\alpha}-1}   \widetilde{\widetilde{\mathcal{B}^{\ast}}}_{L}^{l,l',\big(\frac{\xi'}{\xi}\big)^{p^{\alpha}}} (\xi')^{p^{\alpha}} \bigg(\frac{\xi}{\xi'}\bigg)^{(r_{0}-r)}\text{ }
c_{1}^{(l,\xi)}(r) \text{ }
c_{2}^{(l',\xi')}(p^{\alpha}+r_{0}-r) \bigg)
\end{multline}
\end{Proposition}

\begin{proof} Let us prove first the part of result concerning the power series expansion of $c$ at $r_{0}=0$. For any $u \in \mathbb{N}^{\ast}$, we have 

\begin{equation*} c(p^{\alpha}u) = \sum_{m=1}^{p^{\alpha}u-1} c_{1}(m)c_{2}(p^{\alpha}u-m)
 =
\sum_{u' = 1}^{u-1}
c_{1}\big(p^{\alpha}u'\big) c_{2}\big(p^{\alpha}(u-u')\big) + 
\sum_{u' = 0}^{u-1} \sum_{r=1}^{p^{\alpha}-1} c_{1}\big(p^{\alpha}u'+r\big)
c_{2}\big(p^{\alpha}(u-1-u')+p^{\alpha} - r\big) 
\end{equation*}

Using the hypothesis that $c_{1}$ and $c_{2}$ are in $\LA_{\mathcal{S}_{\alpha}}(\mathbb{Z}_{p}^{(N)},K)$, we obtain, with the notations of Definition \ref{def of the subalgebra},

\begin{multline} \label{eq: un} c(p^{\alpha}u)  = \sum_{u' = 1}^{u-1}
\sum_{l,l'\in \mathbb{N}}\sum_{\xi,\xi'\in \mu_{N}(K)} c_{1}^{(l,\xi)}(0)\text{ }\xi^{-p^{\alpha}u'}\text{ }(p^{\alpha}u')^{l}\text{ }\text{ }c_{2}^{(l',\xi')}(0)\text{ }{\xi'}^{-p^{\alpha}(u-u')}\text{ }(p^{\alpha}(u-u'))^{l'}
\\ + \sum_{u' = 0}^{u-1} \sum_{r=1}^{p^{\alpha}-1} 
\sum_{l,l'\in \mathbb{N}}\sum_{\xi,\xi'\in \mu_{N}(K)}
c_{1}^{(l,\xi)}(r)\text{ }\xi^{-p^{\alpha}u'}\text{ }(p^{\alpha}u')^{l}\text{ }
\text{ }c_{2}^{(l',\xi')}(p^{\alpha}-r)\text{ }{\xi'}^{-p^{\alpha}(u-1-u')}\text{ }(p^{\alpha}(u-u'))^{l'}
\end{multline}

Expressing the sums $\displaystyle\sum_{u'=1}^{u-1}\bigg(\frac{\xi'}{\xi}\bigg)^{p^{\alpha}u'} {u'}^{l}(u-u')^{l'}$ and $\displaystyle\sum_{u'=0}^{u-1}\bigg(\frac{\xi'}{\xi}\bigg)^{p^{\alpha}u'} {u'}^{l}(u-1-u')^{l'}$, which appear in (\ref{eq: un}), by means of Lemma \ref{def of coeff B}, and inverting a double absolutely convergent series, we obtain
$$ c(p^{\alpha}u) = \sum_{L\in\mathbb{N}}\sum_{\xi \in \mu_{N}(K)} c^{(L,\xi)}(0) \xi^{-p^{\alpha}u} (p^{\alpha}u)^{L}$$
where $c^{(L,\xi)}(0)$ is as in equation (\ref{eq:product at 0}).
\newline\indent Let us prove that the sequences $(c^{(L,\xi)}(0))_{L\in\mathbb{N}}$ of equation (\ref{eq:product at 0}) are in $\mathcal{S}_{\alpha}$. For convenience, let us characterize $\mathcal{S}_{\alpha}\subset K^{\mathbb{N}}$ as the subspace of the sequences 
$b=(b_{l})_{l \in \mathbb{N}}$ such that there exist $\kappa_{b},\text{ }\kappa_{b}',\text{ }\kappa_{b}'' \in \mathbb{R}_{+}^{\ast}$ such that, for all $l \in \mathbb{N}$, we have :
$v_{p}(b_{l})\geqslant -\kappa_{b} - \frac{\kappa_{b}'}{\log(p)}\log(l+\kappa_{b}'')-(\alpha-1)l$. By the hypothesis that $c_{1}$ and $c_{2}$ are in $\LA_{\mathcal{S}_{\alpha}}(\mathbb{Z}_{p}^{(N)}),K)$, we can find, for $i\in \{1,2\}$ $\kappa_{c_{i}},\text{ }\kappa'_{c_{i}},\text{ }\kappa''_{c_{i}}\in \mathbb{R}_{+}^{\ast}$ such that the inequality $v_{p}(c_{i}^{(l,\xi)}(r))\geqslant -\kappa_{c_{i}} - \frac{\kappa_{c_{i}}'}{\log(p)}\log(l+\kappa_{c_{i}}'')-(\alpha-1)l$ is satisfied for all $\xi \in \mu_{N}(K)$ and all $r \in \{0,\ldots,p^{\alpha}-1\}$. Then, by equation (\ref{eq:product at 0}),
\begin{multline*} v_{p} ( c^{(L,\xi)}(0) ) \geqslant  \inf_{\substack{l,l' \in \mathbb{N} \\ l+l'+1 \geqslant  L}} \bigg( \alpha(l+l'-L) - (\alpha-1)(l+l')
\\ - (1+ \frac{\log(l+l'+1)}{\log(p)} + \kappa_{c_{1}} + \frac{\kappa'_{c_{1}}}{\log(p)}\log(l+\kappa''_{c_{1}}) + \kappa_{c_{2}} + \frac{\kappa'_{c_{2}}}{\log(p)}\log(l'+\kappa_{c_{2}})) \bigg)
\end{multline*}
$$\text{ }\text{ }\text{ }\text{ }\text{ }
\text{ }\text{ }\text{ }\text{ }\text{ }
\text{ }\text{ }\text{ }\text{ }\text{ }  
\text{ }\text{ }\text{ }
\geqslant  \inf_{\delta\in [-1,+\infty[} \bigg( \delta  - (1+ \kappa_{c_{1}}+ \kappa_{c_{2}}) - \frac{1+\kappa'_{c_{1}}+\kappa'_{c_{2}}}{\log(p)}\log(L+\delta+\max(1,\kappa''_{c_{1}},\kappa''_{c_{2}})) - (\alpha-1)L \bigg) $$
 where we have set $\delta = l+l'-L$. The map $\delta \mapsto \delta - \frac{1+\kappa'_{c_{1}}+\kappa'_{c_{2}}}{\log(p)}\log(l+\delta+\max(1,\kappa''_{c_{1}},\kappa''_{c_{2}}))$ is increasing on 
$\big[\frac{1+ \kappa'_{c_{1}} + \kappa'_{c_{2}}}{\log(p)} - l - \max(1,\kappa''_{c_{1}},\kappa''_{c_{2}}), +\infty\Big)$, and for $L\geqslant  \frac{1+ \kappa'_{c_{1}} + \kappa'_{c_{2}}}{\log(p)} -  \max(1,\kappa''_{c_{1}},\kappa''_{c_{2}})+1$, we have :
$\big[\frac{1+ \kappa'_{c_{1}} + \kappa'_{c_{2}}}{\log(p)} - l - \max(1,\kappa''_{c_{1}},\kappa''_{c_{2}}), +\infty\Big) \subset \big[-1,+\infty\big)$ ; thus 
\begin{equation} \label{eq:the inequality proof} v_{p} ( c^{(L,\xi)}(0)) \geqslant  - \big(2 + \kappa_{c_{1}} + \kappa_{c_{2}} \big) - \frac{1+\kappa'_{c_{1}}+\kappa'_{c_{2}}}{\log(p)}\log(l-1+\max(1,\kappa''_{c_{1}},\kappa''_{c_{2}})) - (\alpha-1)L 
\end{equation}
Since the inequality (\ref{eq:the inequality proof}) is true for $L \in \mathbb{N}$ outside the set $\big[0,\frac{1+ \kappa'_{c_{1}} + \kappa'_{c_{2}}}{\log(p)} -  \max(1,\kappa''_{c_{1}},\kappa''_{c_{2}})+1\big]\cap \mathbb{N}$ which is finite and independent of $L$, we deduce that $(c^{(L,\xi)}(0))\in \mathcal{S}_{\alpha}$.
\newline\indent The proof of the analogous results concerning the power series expansion of $c$ at $r_{0} \in \{1,\ldots,p^{\alpha}-1\}$ is entirely similar, starting with the following equation : for any $u \in \mathbb{N}$ and $r_{0} \in \{1,\cdots,p^{\alpha}-1\}$,
\begin{multline*} c(r_{0}+p^{\alpha}u)
= \sum_{u' = 1}^{u}
c_{1}\big(p^{\alpha}u'\big) c_{2}\big(p^{\alpha}(u-u')+r_{0}\big)
\\ + 
\sum_{u' = 0}^{u} \sum_{r=1}^{r_{0}-1} c_{1}\big(p^{\alpha}u'+r\big)
c_{2}\big(p^{\alpha}(u-u')+ r_{0} - r\big) 
+ 
\sum_{u' = 0}^{u-1} \sum_{r=r_{0}}^{p^{\alpha}-1} c_{1}\big(p^{\alpha}u'+r\big)
c_{2}\big(p^{\alpha}(u-1-u')+r_{0}+p^{\alpha} - r\big) 
\end{multline*}
\end{proof}

\subsection{Regularized $p$-adic iterated integrals ; definition and computation by induction on the depth}

Let the operator "unique primitive which vanishes at $0$" on power series with coefficients in $K$ :
\begin{equation} \label{eq:unique primitive at 0}
\int_{0}^{z} : f= \sum_{m \in \mathbb{N}} c_{m}z^{m} \in K[[z]]\mapsto \sum_{m \in \mathbb{N}} c_{m}\frac{z^{m+1}}{m+1} \in K[[z]]
\end{equation}
The variant of the KZ equation (\ref{eq:equation KZ with congruences}) satisfied by the iterated integrals of Definition \ref{definition iterated integrals} can be reformulated as 
$$ \int_{p,\alpha} e_{0}^{n_{d}-1}e_{{\xi_{d}}^{u_{d}}} \ldots e_{0}^{n_{1}-1}e_{{\xi_{1}}^{u_{1}}} = 
\left\{ \begin{array}{ll} \displaystyle \int_{0}^{z} p^{\alpha}\frac{dz}{z} \bigg( \int_{p,\alpha} e_{0}^{n_{d}-2}e_{{\xi_{d}}^{u_{d}}} \ldots e_{0}^{n_{1}-1}e_{{\xi_{1}}^{u_{1}}}\bigg) & \text{ if } n_{d} \geqslant 2 \\
\displaystyle\int_{0}^{z} p^{\alpha}\frac{dz}{z-\xi} \bigg( \int_{p,\alpha} e_{0}^{n_{d-1}-1}e_{{\xi_{d-1}}^{u_{d-1}}} \ldots e_{0}^{n_{1}-1}e_{{\xi_{1}}^{u_{1}}} \bigg) &\text{ if }n_{d}=1\text{ and }u_{d}=1  \\
\displaystyle\int_{0}^{z} \frac{d(z^{p^{\alpha}})}{z^{p^{\alpha}}-\xi^{p^{\alpha}}} \bigg( \int_{p,\alpha} e_{0}^{n_{d-1}-1}e_{{\xi_{d-1}}^{u_{d-1}}} \ldots e_{0}^{n_{1}-1}e_{{\xi_{1}}^{u_{1}}} \bigg) & \text{ if }n_{d}=1\text{ and }u_{d}=(p^{\alpha})
\end{array}  \right. $$
\indent We are going to define a regularized variant of the above system of equations whose solution lies in the subalgebra introduced in \S3.1. We start by defining a variant of the space $\LA_{\mathcal{S}_{\alpha}}(\mathbb{Z}_{p}^{(N)},K)$ with poles.

\begin{Definition} Let $\LA_{\mathcal{S}_{\alpha}}^{\text{pol}}(\mathbb{Z}_{p}^{(N)},K)$ be the set of functions $\mathbb{Z}_{p} \setminus \{0\} \rightarrow K$ satisfying series expansion as in (\ref{eq:condition series expansion}), such that, if $r_{0} \in \{0,\cdots,p^{\alpha}-1\}$, we have a power series expansion, $|m-r_{0}|_{p} \leqslant p^{-\alpha}$, of the following form 
$$ c(m) = \left\{ \begin{array}{ll} \sum\limits_{l \geqslant 0} \sum\limits_{\xi \in \mu_{N}(K)} c^{(l,\xi)}(r_{0}) \xi^{-m}(m-r_{0})^{l} & r_{0}\in \{1,\cdots,p^{\alpha}-1\}
\\  \sum\limits_{l \geqslant -1} \sum\limits_{\xi \in \mu_{N}(K)} c^{(l,\xi)}(0) \xi^{-m}(m-r_{0})^{l} & r_{0}=0
\end{array}\right.$$
where $(c^{(l,\xi)}(r_{0}))_{l\in \mathbb{N}}$ is a sequence in $\mathcal{S}_{\alpha}$ for any $r_{0} \in \{0,\cdots,p^{\alpha}-1\}$ and any $\xi \in \mu_{N}(K)$.
\end{Definition}

\begin{Example} \label{example of polar element LA}For any $\xi_{0} \in \mu_{N}(K)$, $\Cf_{0}\bigg( \displaystyle\int_{0}^{z}\omega_{\xi_{0}^{p^{\alpha}}}(z^{p^{\alpha}}) = -p^{\alpha}\sum_{m\in \mathbb{N}^{\ast}} \frac{z^{p^{\alpha}m}}{p^{\alpha}m \xi_{0}^{p^{\alpha}m}} \bigg)\in  \LA_{\mathcal{S}_{\alpha}}^{\text{pol}}(\mathbb{Z}_{p}^{(N)},K)$ and 
$\Cf_{0}\bigg( \displaystyle\int_{0}^{z}\omega_{\xi^{p^{\alpha}}}(z^{p^{\alpha}})^{(l,\xi)}(r_{0})= \left\{\begin{array}{ll} -p^{\alpha} &\text{ if }(l,\xi,r_{0})=(-1,\xi_{0},0)
\\ 0 & \text{ otherwise }
\end{array}\right.$
\end{Example}

Now we decompose $\LA_{\mathcal{S}_{\alpha}}^{\text{pol}}(\mathbb{Z}_{p}^{(N)},K)$ into a regular part and a "pure pole" part.

\begin{Lemma} \label{lemma direct sum decomposition} We have $ \LA_{\mathcal{S}_{\alpha}}^{\text{pol}}(\mathbb{Z}_{p}^{(N)},K) = \displaystyle \LA_{\mathcal{S}_{\alpha}}(\mathbb{Z}_{p}^{(N)},K) \bigoplus \bigg(\oplus_{\xi \in \mu_{N}(K)} K.\int_{0}^{z}\omega_{\xi^{p^{\alpha}}}(z^{p^{\alpha}})\bigg)$
\end{Lemma}

\begin{proof} Let us denote by $\xi_{0}$ a primitive $N$-th root of unity in $K$. Let $V_{\xi}=(\xi^{ij})_{\substack{1 \leqslant i \leqslant N \\ 1 \leqslant j \leqslant N}} \in M_{N}(\mathbb{Q}(\xi))$ the Vandermonde matrix associated with the sequence $(\xi^{1},\ldots,\xi^{N})$. An element $c \in \LA_{\mathcal{S}_{\alpha}}(\mathbb{Z}_{p}^{(N)},K) \cap \Vect\bigg(\int_{p,\alpha}\omega_{\xi^{p^{\alpha}}}(z^{p^{\alpha}}),\text{ }|\text{ }\xi \in \mu_{N}(K)\bigg)$ is in particular continuous when $|m|_{p}$ tends to $0$ and $m$ remains in a given congruence class modulo $N$. Thus,  $\lim_{\substack{|m|_{p} \rightarrow 0 \\ m \equiv m_{0} \mod N}}  m c(m)=0$ for all $m_{0} \in \{1,\ldots,M\}$. Moreover, we have 
$\begin{pmatrix}
\lim_{\substack{|m|_{p} \rightarrow 0 \\ m \equiv 1 \mod N}}  m c(m)
\\ \vdots 
\\  \lim_{\substack{|m|_{p} \rightarrow 0 \\ m\equiv N \mod N}}  m c(m)
\end{pmatrix}
= V_{\xi}
\begin{pmatrix}
c^{(-1,\xi^{1})}(0)
\\ \vdots 
\\ c^{(-1,\xi^{N})}(0)
\end{pmatrix}$ and $V_{\xi}$ is invertible. This proves that $\LA_{\mathcal{S}_{\alpha}}(\mathbb{Z}_{p}^{(N)},K) \cap \Vect\bigg(\int_{p,\alpha}\omega_{\xi^{p^{\alpha}}}(z^{p^{\alpha}}),\text{ }|\text{ }\xi \in \mu_{N}(K)\bigg)=\{0\}$. The rest of the statement follows easily, using Example \ref{example of polar element LA}.
\end{proof}

\begin{Lemma} \label{lemma int dz/z polar} Let $f \in A(U^{\an})$ such that $f(0)=0$. If $f \in \Cf_{0}^{-1}\LA_{\mathcal{S}_{\alpha}}(\mathbb{Z}_{p}^{(N)},K)$, then  $\int^{z}_{0} f \omega_{0} \in \Cf_{0}^{-1}\LA_{\mathcal{S}_{\alpha}}^{\text{pol}}(\mathbb{Z}_{p}^{(N)},K)$ (where $\int_{0}^{z}$ is defined in equation (\ref{eq:unique primitive at 0})). Moreover, 
$\Cf_{0}\bigg( \int^{z}_{0} f \omega_{0}\bigg)^{(-1,\xi)}(0)=\Cf_{0}(f)^{(0,\xi)}(0)$.
\end{Lemma}

\begin{proof} Let $c=\Cf_{0}(f)$. We have $\displaystyle f(z)=\sum_{m \in \mathbb{N}^{\ast}} c(m)z^{m}$, and $\displaystyle\int^{z}_{0} f \omega_{0} = \sum_{m \in \mathbb{N}^{\ast}} \frac{c(m)}{m}z^{m}$. By the assumption, there exists, for any $r_{0} \in \{0,1,\ldots,p^{\alpha}-1\}$ and $\xi \in \mu_{N}(K)$, there exist sequences $(c^{(l,\xi)}(r_{0}))_{l\in\mathbb{N}} \in \mathcal{S}_{\alpha}$, such that, for $|m-r_{0}|_{p}\leqslant p^{-\alpha}$,
$\displaystyle c(m) = \sum_{l \in \mathbb{N}} \sum_{\xi \in \mu_{N}(K)} c^{(l,\xi)}(r_{0}) \xi^{-m}(m-r_{0})^{l}$. By $(c^{(l,\xi)}(r_{0}))_{l\in\mathbb{N}} \in \mathcal{S}_{\alpha}$, there exists a polynomial $P_{\xi,r_{0}}$ such that, for all $l \in \mathbb{N}$, $|c^{(l,\xi)}(r_{0})|_{p} \leqslant P_{\xi,r_{0}}(l)p^{(\alpha-1)l}$. Then :
\newline (a) For $|m|_{p}\leqslant p^{-\alpha}$, we have
\begin{equation} \label{eq:computation cm sur m r0=0}\displaystyle\frac{c(m)}{m}=  \sum_{l \geqslant -1} \sum_{\xi \in \mu_{N}(K)} c^{(l+1,\xi)}(r_{0}) \xi^{-m}m^{l}
\end{equation}
and the sequences $(c^{(l+1,\xi)}(r_{0}))_{l \in \mathbb{N}}$ are clearly in $\mathcal{S}_{\alpha}$.
\newline (b) For $|m-r_{0}|_{p}\leqslant p^{-\alpha}$ with $r_{0} \in \{1,\cdots,p^{\alpha}-1\}$, we have
\begin{equation} \label{eq:computation cm sur m r0 not 0}\displaystyle\frac{c(m)}{m}= \bigg(\sum_{l \in \mathbb{N}} \sum_{\xi \in \mu_{N}(K)} c^{(l,\xi)}(r_{0}) \xi^{-m}(m-r_{0})^{l}\bigg)\bigg( \sum_{l'\in\mathbb{N}} \frac{\big(m-r_{0}\big)^{l'}}{r_{0}^{l'+1}}\bigg) = \sum_{l''\in \mathbb{N}} \sum_{\xi \in \mu_{N}(K)}
\sum_{l=0}^{l''} \frac{c^{(l,\xi)}(r_{0})}{r_{0}^{l''-l}}  \xi^{-m}(m-r_{0})^{l}
\end{equation}
Moreover, we have $|\frac{1}{r_{0}}|_{p}\leqslant p^{\alpha-1}$. Thus, with $P_{\xi,r_{0}}$ defined above,
$$\displaystyle\bigg|\sum_{l=0}^{l''} \frac{c^{(l,\xi)}(r_{0})}{r_{0}^{l''-l}}\bigg|_{p}\leqslant \sum_{l=0}^{l''}P_{\xi,r_{0}}(l) p^{(\alpha-1)(l+l''-l)}=\sum_{l=0}^{l''}P_{\xi,r_{0}}(l) p^{(\alpha-1)l''}$$ and there exists a polynomial $Q_{\xi,r_{0}}$ such that for all $l'' \in \mathbb{N}$, $\sum_{l=0}^{l''}P_{\xi,r_{0}}(l)=Q_{\xi,r_{0}}(l'')$.
\end{proof}

\begin{Definition} \label{def regularization integrals} Let $\reg^{\LA} : \LA^{\text{pol}}(\mathbb{Z}_{p}^{(N)},K) \twoheadrightarrow \LA_{\mathcal{S}_{\alpha}}(\mathbb{Z}_{p}^{(N)},K)$ be the projection associated with the direct sum decomposition of Lemma {lemma direct sum decomposition}. Let $\mathcal{O}^{\sh,e_{0\cup \mu_{N}\cup \mu_{N}^{(p^{\alpha})}}}_{\conv}\subset\mathcal{O}^{\sh,e_{0\cup \mu_{N}\cup \mu_{N}^{(p^{\alpha})}}}$ be the vector subspace generated by words whose rightmost letter is not $e_{0}$. We define an operator
$$ \displaystyle\Reg\int_{p,\alpha} : \mathcal{O}^{\sh,e_{0\cup \mu_{N} \cup \mu_{N}^{(p^{\alpha})}}} \rightarrow \Cf_{0}^{-1}\LA_{\mathcal{S}_{\alpha}}(\mathbb{Z}_{p}^{(N)},K) $$ 
as the linear map defined by induction on the weight as follows : for any word $w=e_{0}^{n_{d}-1}e_{{\xi_{d}}^{u_{d}}} \ldots e_{0}^{n_{1}-1}e_{{\xi_{1}}^{u_{1}}}$,
\begin{multline} \label{eq:equation KZ with congruences regularized} \Reg\int_{p,\alpha} e_{0}^{n_{d}-1}e_{{\xi_{d}}^{u_{d}}} \ldots e_{0}^{n_{1}-1}e_{{\xi_{1}}^{u_{1}}} 
= 
\\ \left\{ \begin{array}{ll} \displaystyle \Cf_{0}^{-1}\reg_{\LA}\Cf_{0} \int_{0}^{z} p^{\alpha}\frac{dz}{z} \bigg(  \Reg\int_{p,\alpha} e_{0}^{n_{d}-2}e_{{\xi_{d}}^{u_{d}}} \ldots e_{0}^{n_{1}-1}e_{{\xi_{1}}^{u_{1}}} \bigg) &\text{ if } n_{d}\geqslant 2 
\\ \displaystyle \Cf_{0}^{-1}\reg_{\LA}\Cf_{0}\int_{0}^{z} p^{\alpha} \frac{dz}{z-\xi} \bigg(  \Reg\int_{p,\alpha} e_{0}^{n_{d-1}-1}e_{{\xi_{d-1}}^{u_{d-1}}} \ldots e_{0}^{n_{1}-1}e_{{\xi_{1}}^{u_{1}}}\bigg) &\text{ if } n_{d}=1 \text{ and }u_{d}=1  
\\\displaystyle \Cf_{0}^{-1}\reg_{\LA}\Cf_{0}  \int_{0}^{z} \frac{d(z^{p^{\alpha}})}{z^{p^{\alpha}}-\xi^{(p^{\alpha})}} \bigg(  \Reg\int_{p,\alpha} e_{0}^{n_{d-1}-1}e_{{\xi_{d-1}}^{u_{d-1}}} \ldots e_{0}^{n_{1}-1}e_{{\xi_{1}}^{u_{1}}} \bigg) &\text{ if } n_{d}=1 \text{ and }u_{d}=(p^{\alpha}) 
\end{array}\right.
\end{multline}
\end{Definition}

This definition is consistent by the previous Lemma \ref{lemma int dz/z polar} : namely, this Lemma it proves by induction on the weight of a word $w$ that $\Reg\int_{p,\alpha}(w)$ is well-defined as an element of $\Cf_{0}^{-1}\LA_{\mathcal{S}_{\alpha}}(\mathbb{Z}_{p}^{(N)},K)$.
\newline\indent We now compute the regularized $p$-adic iterated integrals of Definition \ref{def regularization integrals}. This is a regularized analogue of equation (\ref{eq:power series it int}). The formula is inductive on the depth.
\newline\indent In the notation $ \mathcal{B}^{(l,l',\xi)}_{l_{0}}$, $ {\mathcal{B}^{\ast}}^{(l,l',\xi)}_{l_{0}}$
$\widetilde{\mathcal{B}^{\ast}}^{(l,l',\xi)}_{l_{0}}$
$\widetilde{\widetilde{\mathcal{B}^{\ast}}}^{(l,l',\xi)}_{l_{0}}$ defined in Lemma \ref{def of coeff B}, we omit $l'$ when it is equal to $0$. We note that the definition implies $\widetilde{\widetilde{\mathcal{B}^{\ast}}}_{l_{0}}^{(l,\xi)} = \widetilde{\mathcal{B}^{\ast}}_{l_{0}}^{(l,\xi)} - \delta_{l,l_{0}}$.

\begin{Proposition} \label{prop comp reg it int} 
Let $w_{d} = e_{0}^{n_{d}-1}e_{\xi_{d}^{u_{d}}}\cdots e_{0}^{n_{1}-1}e_{\xi_{1}^{u_{1}}} \in \mathcal{O}_{\conv}^{\sh,e_{0 \cup \mu_{N} \cup \mu_{N}^{(p^{\alpha})}}}$ and let  $c_{d}=\Cf_{0}\Reg\int_{p,\alpha} w_{d}$. 
\newline Let $w_{d-1}=e_{0}^{n_{d-1}-1}e_{\xi_{d-1}^{u_{d-1}}}\cdots e_{0}^{n_{1}-1}e_{\xi_{1}^{u_{1}}}$ and $c_{d-1}=\Cf_{0}\Reg\int_{p,\alpha} w_{d-1}$. For any $l \in \mathbb{N}$, $\xi \in \mu_{N}(K)$ and $r_{0} \in \{1,\cdots,p^{\alpha}-1\}$, we have 
\newline (i) if $u_{d}=1$, then
\begin{multline} \label{eq:depth 0 1} c_{d}^{(l,\xi)}(0) =
-(p^{\alpha})^{n_{d}} \sum_{L\geqslant  -1} (p^{\alpha})^{L}
\bigg[ \mathcal{B}_{l+n_{d}}^{l+n_{d}+L,\big( \frac{\xi_{d}}{\xi}\big)p^{\alpha}} \text{ }c_{d-1}^{(l+n_{d}+L,\xi)}(0)
+ \sum_{r=1}^{p^{\alpha}-1}  \widetilde{\widetilde{\mathcal{B}^{\ast}}}_{l+n_{d}}^{l+n_{d}+L,\big(\frac{\xi_{d}}{\xi}\big)^{p^{\alpha}}}\xi^{r} \text{ } c_{d-1}^{(l+n_{d}+L,\xi)}(r) \bigg]
\end{multline}
\begin{multline} \label{eq:depth r 1}
c_{d}^{(l,\xi)}(r_{0}) = 
- (p^{\alpha})^{n_{d}} \sum_{l_{1}=0}^{l} \frac{{-n_{d} \choose l-l_{1}}}{r^{l-l_{1}+n_{d}}} \Bigg[ \sum_{L\geqslant  -1} (p^{\alpha})^{L} \bigg[  {\mathcal{B}^{\ast}}_{l_{1}}^{l_{1}+L,(\frac{\xi_{d}}{\xi})^{p^{\alpha}}} \bigg(\frac{\xi_{d}}{\xi}\bigg)^{r_{0}} \text{ } c_{d-1}^{(l_{1}+L,\xi)}(0)
\\ +  \widetilde{\widetilde{\mathcal{B}^{\ast}}}_{l_{1}}^{l_{1}+L,(\frac{\xi_{d}}{\xi})^{p^{\alpha}}} \xi_{d}^{p^{\alpha}} 
\sum_{r=1}^{p^{\alpha}-1} \bigg(\frac{\xi_{d}}{\xi}\bigg)^{(r_{0}-r)} \text{ }c_{d-1}^{(l_{1}+L,\xi)}(r_{0}) \bigg]
 - \bigg(\frac{\xi_{d}}{\xi}\bigg)^{p^{\alpha}}
\sum_{r=1}^{r_{0}}\bigg(\frac{\xi_{d}}{\xi}\bigg)^{(r_{0}-r)} \text{ }c_{d-1}^{(l_{1},\xi)}(r) \Bigg]
\end{multline}
(ii) if $u_{d}=(p^{\alpha})$, then
\begin{equation} \label{eq:depth 0 p a}c_{d}^{(l,\xi)}(0)
= - (p^{\alpha})^{n_{d}} \sum_{L\geqslant  -1} (p^{\alpha})^{L} \bigg[  \widetilde{\mathcal{B}^{\ast}}^{l+n_{d};l+n_{d}+L}_{(\frac{\xi}{\xi_{d}})^{p^{\alpha}}} \text{ }c_{d-1}^{(l+n_{d}+L,\xi)}(0)] \bigg]
\end{equation}
\begin{equation} \label{eq:depth r p a}
c_{d}^{(l,\xi)}(r_{0}) =
- (p^{\alpha})^{n_{d}} \sum_{l_{1}=0}^{l} \frac{{-n_{d} \choose l-l_{1}}}{r_{0}^{l-l_{1}+n_{d}}} \sum_{L\geqslant  -1} (p^{\alpha})^{L} \bigg[  \widetilde{\mathcal{B}^{\ast}}_{l_{1}}^{(l_{1}+L,( \frac{\xi}{\xi_{d}})^{p^{\alpha}})} \text{ }c_{d-1}^{(l_{1}+L,\xi)}(r_{0})] \bigg]
\end{equation}
\end{Proposition}

\begin{proof} (a)
We apply Proposition \ref{multiplication} in the particular case $\displaystyle(f_{1},f_{2})=(f,\frac{z}{z-\xi_{d}})$, using Example \ref{example elements of LA} (i). This gives that $c=\Cf_{0}(f \frac{z}{z - \xi_{d}})$ is in $\LA_{\mathcal{S}_{\alpha}}(\mathbb{Z}_{p}^{(N)},K)$ and satisfies, for all $l \in \mathbb{N}$, $\xi \in \mu_{N}(K)$, $r_{0} \in \{1,\ldots,p^{\alpha}-1\}$,
\begin{equation} \label{eq:mult by z sur z - xi 1}
c^{(l,\xi)}(0) = - \sum_{L\geqslant  -1} (p^{\alpha})^{L}
\bigg[ \mathcal{B}_{l}^{l+L,\big(\frac{\xi_{d}}{\xi}\big)^{p^{\alpha}}} \text{ }c^{(l,\xi)}(0)
+ \sum_{r=1}^{p^{\alpha}-1} \widetilde{\widetilde{\mathcal{B}^{\ast}}}_{l}^{l+L,\big(\frac{\xi_{d}}{\xi}\big)^{p^{\alpha}}}\xi^{r} \text{ } c^{(l+L,\xi)}(r) \bigg]
\end{equation}
\begin{multline} \label{eq:mult by z sur z - xi 2}
c^{(l,\xi)}(r_{0}) = - \sum_{L\geqslant  -1} (p^{\alpha})^{L} \bigg[ {}_{\ast} \mathcal{B}_{l}^{l+L,\big(\frac{\xi_{d}}{\xi}\big)^{p^{\alpha}}}\bigg(\frac{\xi}{\xi_{d}}\bigg)^{r_{0}}\text{ } c^{(l+L,\xi)}(0)
+ \xi_{d}^{p^{\alpha}} \widetilde{\widetilde{\mathcal{B}^{\ast}}}_{l}^{l+L,\big(\frac{\xi_{d}}{\xi}\big)^{p^{\alpha}}}  
\sum_{r=1}^{p^{\alpha}-1}  \big(\frac{\xi}{\xi_{d}}\big)^{(r_{0}-r)} \text{ }c^{(l+L,\xi)}(r)  \bigg]
\\ - \big(\frac{\xi}{\xi_{d}}\big)^{p^{\alpha}} 
\sum_{r=1}^{r_{0}}\big(\frac{\xi}{\xi_{d}}\big)^{(r_{0}-r)} \text{ }c^{(l+L,\xi)}(r)
\end{multline}
\newline (b) We apply Proposition \ref{multiplication} in the particular case $\displaystyle(f_{1},f_{2})=(f,\frac{z^{p^{\alpha}}}{z^{p^{\alpha}}-\xi_{d}^{p^{\alpha}}})$, using Example \ref{example elements of LA} (ii). This gives that $c=\Cf_{0} (f \frac{z^{p^{\alpha}}}{z^{p^{\alpha}} - \xi_{d}^{p^{\alpha}}})$ is in $\LA_{\mathcal{S}_{\alpha}}(\mathbb{Z}_{p}^{(N)},K)$ and satisfies, for all $l \in \mathbb{N}$, $\xi \in \mu_{N}(K)$, and $r_{0} \in \{1,\ldots,p^{\alpha}-1\}$,
\begin{equation} \label{eq:mult by z sur z - xi p 1}
c^{(l,\xi)}(0) = - \sum_{L\geqslant  -1} (p^{\alpha})^{L} 
\mathcal{B}_{l}^{l+L,\big(\frac{\xi_{d}}{\xi}\big)^{p^{\alpha}}} \text{ }c^{(l+L,\xi)}(0)
\end{equation}
\begin{equation} \label{eq:mult by z sur z - xi p 2}
c^{(l,\xi)}(r_{0}) = - \sum_{L \geqslant  -1} (p^{\alpha})^{L}  \widetilde{\mathcal{B}^{\ast}}_{l}^{l+L,\big(\frac{\xi_{d}}{\xi}\big)^{p^{\alpha}}} \text{ }c^{(l+L,\xi)}(r_{0})
\end{equation}
(c) Then we apply the following statement to the results of (a) and (b) : let $w \in \mathcal{O}_{\conv}^{\sh,e_{0 \cup \mu_{N} \cup \mu_{N}^{(p^{\alpha})}}}$, $n\in \mathbb{N}$, $c=\Cf_{0} \Reg\int_{p,\alpha}w$ and $\tilde{c}_{n}=\Cf_{0} \Reg\int_{p,\alpha} e_{0}^{n}w$. Then $\tilde{c}$ is the element of  $\LA_{\mathcal{S}_{\alpha}}(\mathbb{Z}_{p}^{(N)},K)$ determined by, for all $l,\xi$ and for all $r_{0} \in \{1,\ldots,p^{\alpha}-1\}$,
\begin{equation} \label{eq:int reg omega0 iteration n times r0=0}
\tilde{c}^{(l,\xi)}(0) = (p^{\alpha})^{n} c^{(l+n,\xi)}(0) 
\end{equation}
\begin{equation} \label{eq:int reg omega0 iteration n times r0}\tilde{c}^{(l,\xi)}(r_{0}) = 
(p^{\alpha})^{n} \sum_{l_{1}=0}^{l} \frac{{-n \choose l-l_{1}+1}c^{(l_{1},\xi)}(r_{0})}{r_{0}^{l-l_{1}+n}} 
\end{equation}
Indeed, (\ref{eq:int reg omega0 iteration n times r0=0}) is deduced from (\ref{eq:computation cm sur m r0=0}) and the definition of the regularization (Definition \ref{def regularization integrals}), by induction on $n$ ; (\ref{eq:int reg omega0 iteration n times r0}) is obtained like (\ref{eq:computation cm sur m r0 not 0}) with $\displaystyle\frac{c(m)}{m^{n}}$ instead of $\displaystyle\frac{c(m)}{m}$, given that the regularization affects only the coefficients at $r_{0}=0$.
\end{proof}

\subsection{Bounds on the norms of regularized $p$-adic iterated integrals}

We deduce from \S3.1 and \S3.2 bounds on the norms of regularized $p$-adic iterated integrals. The norm is the one of Notation \ref{the open affine}.

\begin{Proposition} \label{bound on valuation}For any word $w=e_{0}^{n_{d}-1}e_{\xi_{d}^{u_{d}}}\cdots e_{0}^{n_{1}-1}e_{\xi_{1}^{u_{1}}}$ in $\mathcal{O}^{e_{0\cup\mu_{N}\cup \mu_{N}^{(p^{\alpha})}}}_{\conv}$, letting $n=\sum_{i=1}^{d} n_{i}$ be its weight, we have 
$$ \bigg|\bigg| \Reg\int_{p,\alpha} w
\bigg|\bigg| \leqslant p^{- \big(n - 2d - \frac{d}{\log(p)}\log ( \frac{d}{\log(p)}) -\frac{d}{\log(p)} \log(n+3d) \big)} $$ 
\end{Proposition}

\begin{proof} We proceed in three steps.
\newline\indent (a) We prove by induction of the depth that for each $r_{0} \in \{0,\ldots,p^{\alpha}-1\}$, $l\in\mathbb{N}$, $\xi \in \mu_{N}(K)$ :
\begin{multline} \label{eq:first bound on valuation}v_{p}\bigg(\Cf_{0} \Reg\int_{p,\alpha} w\bigg)^{(l,\xi)}(r)
\geqslant 
\\
\inf_{\substack{L_{d}\in \mathbb{Z}\\ L_{d}\geqslant -d}} \bigg[ \alpha\bigg(\sum_{i=1}^{d}n_{i} +  L_{d}\bigg) - d\bigg(1+  \frac{\log(l+\sum_{i=1}^{d}n_{i}+L_{d}+d)}{\log(p)} \bigg) - (\alpha-1)(l+L_{d}+\sum_{i=1}^{d}n_{i}) \bigg]
\end{multline}
This follows by induction on the depth by applying the result to the word $w_{d-1}=e_{0}^{n_{d-1}-1}e_{\xi_{d-1}^{u_{d-1}}}\cdots e_{0}^{n_{1}-1}e_{\xi_{1}^{u_{1}}}$ by Propostion \ref{prop comp reg it int} and the bounds on the valuations of the coefficients $\mathcal{B}$ in Lemma \ref{def of coeff B}.
\newline\indent The second line in equation (\ref{eq:first bound on valuation}) is also equal to 
\begin{equation} \label{eq:first bound prime on valuation} \displaystyle\inf_{\substack{L_{d}\in \mathbb{Z}\\ L_{d}\geqslant -d}} \bigg[ \sum_{i=1}^{d}n_{i} +  L_{d} - d\bigg(1+  \frac{\log(l+\sum_{i=1}^{d}n_{i}+L_{d}+d)}{\log(p)} \bigg) - (\alpha-1)l \bigg]
\end{equation}
\newline\indent (b) We bound the $\inf$ in equation (\ref{eq:first bound prime on valuation}). Let the function $f : L \in (-l-\sum_{i=1}^{d}n_{i}-d,+\infty) \mapsto L - \frac{d}{\log(p)} \log(L+l+\sum_{i=1}^{d}n_{i}+d)\in \mathbb{R}$ ; it is increasing over $[t_{0},+\infty)$, where $t_{0} = \frac{d}{\log(p)} - (l+\sum_{i=1}^{d}n_{i}+d)$, and it reaches its minimum at $t_{0}$, which is $f(t_{0})=t_{0} - \frac{d}{\log(p)}\log \bigg( \frac{d}{\log(p)}\bigg)$. We distinguish two cases :
\newline - If $t_{0} \geqslant  -d$, then $f(t_{0})\geqslant - d - \frac{d}{\log(p)}\log \big( \frac{d}{\log(p)}\big)$, and thus $\inf_{[-d,\infty)} f \geqslant - d - \frac{d}{\log(p)}\log \big( \frac{d}{\log(p)}\big)$.
\newline - If $t_{0} \leqslant -d$, then $f$ is increasing in $[-d,\infty)$ and 
$\inf_{[-d,\infty)} f \geqslant f(-d) = -d - \frac{d}{\log(p)} \log(l+\sum_{i=1}^{d}n_{i}+1)$. We deduce
$$\begin{array}{ll} \inf_{[-d,\infty)} f &\geqslant \min\bigg(- d - \frac{d}{\log(p)}\log \big( \frac{d}{\log(p)}\big), -d - \frac{d}{\log(p)} \log(l+\sum_{i=1}^{d}n_{i}+1) \bigg) 
\\ & \geqslant -d - \frac{d}{\log(p)} \bigg(  \log(l+\sum_{i=1}^{d}n_{i}+L_{d}+d) + \log(\frac{d}{ \log(p)}) \bigg)
\end{array} $$
We deduce the following bound on the $\inf$ in equation (\ref{eq:first bound prime on valuation})
\begin{multline} \label{eq:second bound on valuation}
\displaystyle\inf_{\substack{L_{d}\in \mathbb{Z}\\ L_{d}\geqslant -d}} \bigg[ \sum_{i=1}^{d}n_{i} +  L_{d} - d\bigg(1+  \frac{\log(l+\sum_{i=1}^{d}n_{i}+L_{d}+d)}{\log(p)} \bigg) - (\alpha-1)l \bigg] 
\\ \geqslant \sum_{i=1}^{d} n_{i} - 2d -\frac{d}{\log(p)} \bigg( \log\big(l+\sum_{i=1}^{d}n_{i}+1\big) + \log\big(\frac{d}{ \log(p)}\big) \bigg) - (\alpha-1)l
\end{multline}
(c) For all $r_{0} \in \{0,\ldots,p^{\alpha}-1\}$, $m \in \mathbb{N}^{\ast}$ such that $|m-r_{0}|_{p}\leqslant p^{-\alpha}$ (i.e. $v_{p}((m-r_{0})^{l})\geqslant \alpha l$), we have $\displaystyle (\Cf_{0}\Reg \int_{p,\alpha}w)(m) = \sum_{l\geqslant \mathbb{N}} \sum_{\xi \in \mu_{N}(K)} c^{(l,\xi)}(m-r_{0})^{l}$. By (\ref{eq:second bound on valuation}), we deduce
\begin{equation} \label{eq:third bound on valuation} v_{p}(c(m)) \geqslant \inf_{l\in \mathbb{N}} \bigg[ l + \sum_{i=1}^{d} n_{i} - 2d -\frac{d}{\log(p)} \bigg(  \log(l+\sum_{i=1}^{d}n_{i}+1) + \log\big(\frac{d}{ \log(p)}\big) \bigg) \bigg]
\end{equation}
The function $g : l \mapsto l - \frac{d}{\log(p)} \log(l+\sum_{i=1}^{d}n_{i}+1)$ is increasing on $[t_{1},\infty)$ with $t_{1}=\frac{d}{\log(p)}-\sum_{i=1}^{d} n_{i}-1$, and reaches its minimum at $t_{1}$, which is equal to $g(t_{1})=t_{1} - \frac{d}{\log(p)}\log\big(\frac{d}{\log(p)}\big)$.
Since we have $n_{i}\geqslant 1$ for all $i$, we have $t_{1} \leqslant \frac{d}{\log(p)} - d - 1\leqslant 4d-d-1=3d-1 $. We again distinguish two cases.
\newline - If $t_{1}\geqslant 0$, then since $t_{1}\leqslant 3d-1$ we have $\inf_{l \in \mathbb{N}}g(l) \geqslant \displaystyle\min_{l\in [0,3d-1] \cap \mathbb{N}} g(l) \geqslant 0 - \frac{d}{\log(p)} \log(3d + \sum_{i=1}^{d} n_{i})$.
\newline - If $t_{1} \leqslant 0$ then $\inf_{l \in \mathbb{N}}g(l) = g(0) = -\frac{d}{\log(p)}\log(\sum_{i=1}^{d} n_{i}+1)$.
\newline In both cases, we have $\inf_{l \in \mathbb{N}} g(l) \geqslant - \frac{d}{\log(p)} \log(\sum_{i=1}^{d} n_{i}+3d)$. Thus equation (\ref{eq:third bound on valuation}) becomes 
\begin{equation} v_{p}\bigg(\big(\Cf_{0}\Reg \int_{p,\alpha}w\big)(m)\bigg) \geqslant \sum_{i=1}^{d} n_{i} - 2d - \frac{d}{\log(p)}\log \big( \frac{d}{\log(p)}\big) -\frac{d}{\log(p)} \log(\sum_{i=1}^{d} n_{i}+3d)
\end{equation}
This implies the result, since the map $\Cf_{0}$ is an isometry (Proposition \ref{prop Mahler general}).
\end{proof}

\begin{Remark} The bound in Proposition \ref{bound on valuation} becomes better when $p$ becomes large.
\end{Remark}

\section{Regularity of overconvergent $p$-adic multiple polylogarithms and end of the proof\label{conclusion}}

We prove (Proposition \ref{prop characterization in terms of regularized}) that, in the decomposition of overconvergent $p$-adic multiple polylogarithms of Proposition \ref{prop decomposition}, the iterated integrals of Definition \ref{definition iterated integrals} can be replaced by their regularized variants of Definition \ref{def regularization integrals} and that the decomposition map of Definition \ref{definition of dec} can be replaced by a "regularized" variant (Definition \ref{definition of reg dec}). This is a consequence of the analytic nature of the overconvergent $p$-adic multiple polylogarithms.
This gives a characterization of $\Li_{p,\alpha}^{\dagger}$ and  $\Ad_{\Phi_{p,\alpha}^{(\xi)}}(e_{\xi})$ $(\xi \in \mu_{N}(K))$ in terms of regularized iterated integrals. This characterization can be written by an induction on the depth (Proposition \ref{induction on depth for reg dec}). Combining this characterization and the bounds on the norms of regularized iterated integrals (Proposition \ref{bound on valuation}) this enables to finish the proof of the main theorem (Proposition \ref{end of proof of theorem}).

\subsection{Regularized decomposition of overconvergent $p$-adic multiple polylogarithms in terms of regularized iterated integrals}

\begin{Lemma} We have \label{intersection of spaces} $\LA_{\mathcal{S}_{\alpha}}^{\text{polar}}(\mathbb{Z}_{p}^{(N)},K) \cap \mathcal{C}(\mathbb{Z}_{p}^{(N)},K) = \LA_{\mathcal{S}_{\alpha}}(\mathbb{Z}_{p}^{(N)},K)$.
\end{Lemma}

\begin{proof} This is a byproduct of the proof of Lemma \ref{lemma direct sum decomposition} and the formula in Example \ref{example of polar element LA}.
\end{proof}

The differential equation of $\Li_{p,\alpha}^{\dagger}$ (\ref{eq:diff coefficient}) amounts to (we are again using (equation \ref{eq:unique primitive at 0}))
\begin{multline} \label{eq:diff coefficient prime}
\Li^{\dagger}_{p,\alpha}[e_{0}^{n_{d}-1}e_{\xi_{d}} \ldots e_{0}^{n_{1}-1}e_{\xi_{1}}e_{0}^{n_{0}-1}] =
\sum_{\substack{w_{1},w_{2}\text{ words } \\ w_{1}w_{2}=w\\ 
		w_{2}\not=\emptyset}} \sum_{\xi \in \mu_{N}(K)}
\bigg(\Ad_{\Phi_{p,\alpha}^{(\xi)}}(e_{\xi})[w_{2}]\bigg)\int_{0}^{z} \omega_{\xi^{p^{\alpha}}}(z^{p^{\alpha}})\Li^{\dagger}_{p,\alpha}[w_{1}] \text{ }+
\\ 
\left\{
\begin{array}{ll} \displaystyle
\int_{0}^{z}p^{\alpha}\omega_{0}(z)\Li^{\dagger}_{p,\alpha}[e_{0}^{n_{d}-2}e_{\xi_{d}} \ldots e_{0}^{n_{1}-1}e_{\xi_{1}}e_{0}^{n_{0}-1}] - p^{\alpha}\omega_{0}(z) \int_{0}^{z}\Li^{\dagger}_{p,\alpha}[e_{0}^{n_{d}-1}e_{\xi_{d}} \ldots e_{0}^{n_{1}-1}e_{\xi_{1}}e_{0}^{n_{0}-2}]
&  \text{ if }n_{d}\geqslant 2,\text{ } n_{0}\geqslant 2
\\ \displaystyle
\int_{0}^{z}p^{\alpha}\omega_{0}(z)\Li^{\dagger}_{p,\alpha}[e_{0}^{n_{d}-2}e_{\xi_{d}} \ldots e_{0}^{n_{1}-1}e_{\xi_{1}}]
& \text{ if }n_{d}\geqslant 2,\text{ } n_{0}=1
\\ \displaystyle
\int_{0}^{z}p^{\alpha}\omega_{\xi_{d}}(z)\Li^{\dagger}_{p,\alpha}[e_{0}^{n_{d}-1}e_{\xi_{d}}\ldots e_{\xi_{1}} e_{0}^{n_{0}-1}] - \int_{0}^{z}p^{\alpha}\omega_{0}(z)
\Li_{p,\alpha}^{\dagger}[e_{\xi_{d-1}} \ldots e_{0}^{n_{1}-1}e_{\xi_{1}}e_{0}^{n_{0}-2}]
& \text{ if }n_{d}=1,\text{ }n_{0}\geqslant 2
\\ \displaystyle \int_{0}^{z}p^{\alpha}\omega_{\xi_{d}}(z)\Li_{p,\alpha}^{\dagger}[e_{0}^{n_{d-1}-1}e_{\xi_{d-1}}\ldots e_{0}^{n_{1}-1}e_{\xi_{1}}]
& \text{ if }n_{d}=1,\text{ }n_{0}=1
\\
\end{array}	\right.
\end{multline}

We note that the first line of (\ref{eq:diff coefficient prime}) contains the term $\displaystyle\sum_{\xi \in \mu_{N}(K)} (\Ad_{\Phi_{p,\alpha}^{(\xi)}}[w]) \omega_{\xi^{p^{\alpha}}}(z^{p^{\alpha}})$, which corresponds to $(w_{1},w_{2})=(\emptyset,w)$. By Example \ref{example of polar element LA}, the integral $\displaystyle\int_{0}^{z}\displaystyle\sum_{\xi \in \mu_{N}(K)} (\Ad_{\Phi_{p,\alpha}^{(\xi)}}[w]) \omega_{\xi^{p^{\alpha}}}(z^{p^{\alpha}}) \in \LA_{\mathcal{S}_{\alpha}}^{\text{pol}}(\mathbb{Z}_{p}^{(N)},K)$ and is a "pure pole" in the sense of the decomposition of Lemma \ref{lemma direct sum decomposition}. We remove it in the next definition.

\begin{Definition} \label{definition of reg dec}Let $\Reg\dec_{p,\alpha} : \mathcal{O}^{\sh,e_{0\cup \mu_{N}}} \rightarrow\mathcal{O}^{\sh,e_{0\cup \mu_{N}\cup \mu_{N}^{(p^{\alpha})}}} \otimes_{\mathbb{Q}} K$ be defined by induction on the weight, as follows : $\Reg\dec_{p,\alpha}(\emptyset)=1$ ; for all $n \in \mathbb{N}^{\ast}$, 
	$\Reg\dec_{p,\alpha}(e_{0}^{n})=e_{0}\Reg\dec_{p,\alpha}(e_{0}^{n-1})-e_{0}\Reg\dec_{p,\alpha}(e_{0}^{n-1})(=0)$, and for any $w= e_{0}^{n_{d}-1}e_{\xi_{d}} \ldots e_{0}^{n_{1}-1}e_{\xi_{1}}e_{0}^{n_{0}-1}$ word over $e_{0\cup\mu_{N}}$, where the $n_{i}$'s are positive integers and the $\xi_{i}$'s are $N$-th roots of unity, 
\begin{multline} \label{eq:definition dec}
\Reg\dec_{p,\alpha}(e_{0}^{n_{d}-1}e_{\xi_{d}} \ldots e_{0}^{n_{1}-1}e_{\xi_{1}}e_{0}^{n_{0}-1}) = \sum_{\substack{w_{1},w_{2}\text{ words } \\ w_{1}w_{2}=w\\ w_{2}\not=\emptyset \\ w_{1}\not=\emptyset}} \sum_{\xi \in \mu_{N}(K)} \bigg( \Ad_{\Phi_{p,\alpha}^{(\xi)}}(e_{\xi})[w_{2}] \bigg) e_{\xi^{(p^{\alpha})}}\Reg\dec_{p,\alpha}(w_{1}) +
\\ 
\left\{ 
\begin{array}{ll} 
e_{0}\Reg\dec_{p,\alpha}(e_{0}^{n_{d}-2}e_{\xi_{d}} \ldots e_{0}^{n_{1}-1}e_{\xi_{1}}e_{0}^{n_{0}-1}) - e_{0}\Reg\dec_{p,\alpha}(e_{0}^{n_{d}-1}e_{\xi_{d}} \ldots e_{0}^{n_{1}-1}e_{\xi_{1}}e_{0}^{n_{0}-2})) 
& \text{ if }n_{d}\geqslant 2, n_{0}\geqslant 2
\\
e_{0}\Reg\dec_{p,\alpha}(e_{0}^{n_{d}-2}e_{\xi_{d}} \ldots e_{0}^{n_{1}-1}e_{\xi_{1}}) 
& \text{ if }n_{d}\geqslant 2, n_{0}=1
\\
e_{\xi_{d}}\Reg\dec_{p,\alpha}(e_{0}^{n_{d-1}-1}e_{\xi_{d-1}}\ldots e_{0}^{n_{0}-1}) - e_{0} \Reg\dec_{p,\alpha}(e_{\xi_{d}} \ldots e_{0}^{n_{1}-1}e_{\xi_{1}}e_{0}^{n_{0}-2})
& \text{ if }n_{d}= 1, n_{0}\geqslant 2
\\e_{\xi_{d}}\Reg\dec_{p,\alpha}(e_{0}^{n_{d-1}-1}e_{\xi_{d-1}}\ldots e_{0}^{n_{1}-1}e_{\xi_{1}})
& \text{ if }n_{d}=1, n_{0}=1
\\
\end{array}	 \right.
\end{multline}
\end{Definition}

The next Proposition is a regularized variant of the statement $\displaystyle\Li_{p,\alpha}^{\dagger}=\int_{p,\alpha}\circ \dec_{p,\alpha}$ (Proposition \ref{prop decomposition}) and, at the same time, a characterization of $\Li_{p,\alpha}^{\dagger}$ and $\Ad_{\Phi_{p,\alpha}^{(\xi)}}(e_{\xi})$ for any $\xi \in \mu_{N}(K)$ in terms of regularized iterated integrals. As in Proposition \ref{prop decomposition} we regard $\Li_{p,\alpha}^{\dagger}$ as a function on $\mathcal{O}^{\sh,e_{0 \cup \mu_{N}\cup \mu_{N}^{(p^{\alpha})}}}$.

\begin{Proposition} \label{prop characterization in terms of regularized} We have :
\begin{equation} \label{eq:regularity of Li} \Li_{p,\alpha}^{\dagger} = \Reg\int_{p,\alpha}   \circ \Reg\dec_{p,\alpha}
\end{equation}
and, for any word $w=e_{0}^{n_{d}-1}e_{\xi_{d}}\ldots e_{0}^{n_{1}-1}e_{\xi_{1}}e_{0}^{n_{0}-1}$, with the $n_{i}$'s positive integers and the $\xi_{i}$'s roots of unity, for all $\xi \in \mu_{N}(K)$,
\begin{multline} \label{eq:vanishing of the polar part} \Ad_{\Phi_{p,\alpha}^{(\xi)})}(e_{\xi})[w] =
\\
\sum_{\substack{w_{1},w_{2} \text{ words } \\ w_{1}w_{2}=w\\ 
		w_{2}\not=\emptyset}} \sum_{\xi \in \mu_{N}(K)}
\bigg(\Ad_{\Phi_{p,\alpha}^{(\xi)}}(e_{\xi})[w_{2}]\bigg)\Cf_{0}\bigg(
\frac{z^{p^{\alpha}}}{z^{p^{\alpha}}-\xi^{p^{\alpha}}}\Li^{\dagger}_{p,\alpha}[w_{1}]\bigg)^{(0,\xi)}(0) \text{ }+
\\ 
\left\{ 
\begin{array}{ll} \displaystyle
p^{\alpha}(\Cf_{0}\Li^{\dagger}_{p,\alpha}[e_{0}^{n_{d}-2}e_{\xi_{d}} \ldots e_{0}^{n_{1}-1}e_{\xi_{1}}e_{0}^{n_{0}-1}])^{(0,\xi)}(0) - p^{\alpha}\Cf_{0}\big(\Li^{\dagger}_{p,\alpha}[e_{0}^{n_{d}-1}e_{\xi_{d}} \ldots e_{0}^{n_{1}-1}e_{\xi_{1}}e_{0}^{n_{0}-2}]\big)^{(0,\xi)}(0)
\\ \text{ if }n_{d}\geqslant 2,\text{ } n_{0}\geqslant 2
\\ \displaystyle
p^{\alpha}(\Cf_{0}\Li^{\dagger}_{p,\alpha}[e_{0}^{n_{d}-2}e_{\xi_{d}} \ldots e_{0}^{n_{1}-1}e_{\xi_{1}}])^{(0,\xi)}(0) \text{ if }n_{d}\geqslant 2,\text{ } n_{0}=1
\\ \displaystyle
p^{\alpha}(\Cf_{0} \frac{z}{z-\xi_{d}} \Li^{\dagger}_{p,\alpha}[e_{0}^{n_{d}-1}e_{\xi_{d}}\ldots e_{\xi_{1}} e_{0}^{n_{0}-1}])^{(0,\xi)}(0) - \Cf_{0}\big( p^{\alpha}
\Li_{p,\alpha}^{\dagger}[e_{0}^{n_{d}-1}e_{\xi_{d}} \ldots e_{0}^{n_{1}-1}e_{\xi_{1}}e_{0}^{n_{0}-2}]\big)^{(0,\xi)}(0)
\\ \text{ if }n_{d}= 1,\text{ }n_{0}\geqslant 2
\\ \displaystyle p^{\alpha}(\Cf_{0} \frac{z}{z-\xi_{d}}\Li_{p,\alpha}^{\dagger}[e_{0}^{n_{d-1}-1}e_{\xi_{d-1}}\ldots e_{\xi_{1}}])^{(0,\xi)}(0) \text{ if }n_{d}=1,\text{ }n_{0}=1
\end{array}	\right.
\end{multline}
\end{Proposition}

\begin{proof} In weight 0, we have $\Li_{p,\alpha}^{\dagger}(\emptyset)=1$ and the result is clear. Assume that the result holds in weight up to $0,\ldots,n-1$. Let $w$ be a word on $e_{0 \cup \mu_{N}}$ of weight $n$. 
For any $\xi \in \mu_{N}(K)$, $\displaystyle\Cf_{0}(\frac{p^{\alpha}z}{z-\xi})$ and $\displaystyle\Cf_{0}(\frac{z^{p^{\alpha}}}{z^{p^{\alpha}}-\xi^{p^{\alpha}}})$ are in $\LA_{\mathcal{S}_{\alpha}}(\mathbb{Z}_{p}^{(N)},K)$ (Example \ref{example elements of LA}). By Proposition \ref{multiplication} and the induction hypothesis, all the products of the form $\displaystyle\Li_{p,\alpha}^{\dagger}[w']\frac{p^{\alpha}z}{z-\xi}$ and $\displaystyle\Li_{p,\alpha}^{\dagger}[w']\frac{z^{p^{\alpha}}}{z^{p^{\alpha}}-\xi^{p^{\alpha}}}$ with $w'$ a word of weight strictly lower than the weight of $w$, are in $\Cf_{0}^{-1}\LA_{\mathcal{S}_{\alpha}}(\mathbb{Z}_{p}^{(N)},K)$. Moreover, for all $\xi \in \mu_{N}(K)$, we have $\displaystyle\omega_{\xi}(z)=\frac{z}{z-\xi}\omega_{0}$ and 
$\displaystyle\omega_{\xi^{p^{\alpha}}}(z^{p^{\alpha}})=\frac{z^{p^{\alpha}}}{z^{p^{\alpha}}-\xi^{p^{\alpha}}}\omega_{0}$. 
By Lemma \ref{lemma int dz/z polar}, we deduce that $\Cf_{0}(\Li_{p,\alpha}^{\dagger}[w])$ is in $\LA_{\mathcal{S}_{\alpha}}^{\text{pol}}(\mathbb{Z}_{p}^{(N)},K)$.
\newline\indent On the other hand, since $\Li_{p,\alpha}^{\dagger}[w] \subset A^{\dagger}(U^{\an}) \subset A(U^{\an})$, 
by Proposition \ref{prop Mahler general}, $\Cf_{0}(\Li_{p,\alpha}^{\dagger}[w]) \in \mathcal{C}(\mathbb{Z}_{p}^{(N)},K)$. Thus,  $\Cf_{0}(\Li_{p,\alpha}^{\dagger}[w]) \in \LA_{\mathcal{S}_{\alpha}}^{\text{pol}}(\mathbb{Z}_{p}^{(N)},K) \text{ }\cap\text{ }\mathcal{C}(\mathbb{Z}_{p}^{(N)},K)$. By Lemma \ref{intersection of spaces}, we deduce $\Cf_{0}(\Li_{p,\alpha}^{\dagger}[w]) \in \LA_{\mathcal{S}_{\alpha}}(\mathbb{Z}_{p}^{(N)},K)$.
\newline By Lemma \ref{lemma direct sum decomposition}, this implies that the polar part, in the sense of the direct sum decomposition of $\LA_{\mathcal{S}_{\alpha}}^{\text{pol}}(\mathbb{Z}_{p}^{(N)},K)$ in Lemma \ref{lemma direct sum decomposition}, vanishes. By Lemma \ref{lemma int dz/z polar}, the vanishing of this polar part is equation (\ref{eq:vanishing of the polar part}).
\end{proof}

\subsection{Computation of the regularized decomposition by induction on the depth}

We write a formula for $\Reg\dec_{p,\alpha}$ which is inductive with respect to the depth (Proposition \ref{induction on depth for reg dec}). 

\begin{Lemma} \label{lemma depth 0}For all $n \in \mathbb{N}^{\ast}$, we have $\Li^{\dagger}_{p,\alpha}[e_{0}^{n}] = 0$, $\Reg\dec_{p,\alpha}(e_{0}^{n})=\dec_{p,\alpha}(e_{0}^{n})=0$, and  $\Ad_{\Phi^{(\xi)}_{p,\alpha}}(e_{\xi})[e_{0}^{n}]=0$ for all $\xi \in \mu_{N}(K)=0$
\end{Lemma}

\begin{proof} The fact that $\Reg\dec_{p,\alpha}(e_{0}^{n})=\dec_{p,\alpha}(e_{0}^{n})=0$ for all $n \in \mathbb{N}^{\ast}$ has already been observed as an immediate byproduct of the definitions (Definition \ref{definition of dec} and Definition \ref{definition of reg dec}). It implies $\Li_{p,\alpha}^{\dagger}[e_{0}^{n}]=0$ for all $n \in \mathbb{N}^{\ast}$ by Proposition \ref{prop decomposition}. Finally, $\Ad_{\Phi^{(\xi)}_{p,\alpha}}(e_{\xi})=\{\Phi^{(\xi)}_{p,\alpha}\}^{-1}e_{\xi}\Phi^{(\xi)}_{p,\alpha}$ is the product of some formal power series and $e_{\xi}$ which is of depth $1$, and it thus contains only terms of depth $\geqslant 1$ ; thus $\Ad_{\Phi^{(\xi)}_{p,\alpha}}(e_{\xi})[e_{0}^{n}]=0$ for all $n \in \mathbb{N}^{\ast}$ (and for $n=0$, which we already used since equation (\ref{eq:diff coefficient})).
\end{proof}
 
\begin{Proposition} \label{induction on depth for reg dec}We have, for any word $w=e_{0}^{n_{d}-1}e_{\xi_{d}}\ldots e_{0}^{n_{1}-1}e_{\xi_{1}}e_{0}^{n_{0}-1}$, with the $n_{i}$'s positive integers and the $\xi_{i}$'s roots of unity,
\begin{multline} \label{eq:premiere equation prime}
\Reg\dec_{p,\alpha}(w) =
\\
\sum_{\substack{0 \leqslant l_{0} \leqslant n_{0}-1 \\ 0 \leqslant l_{d} \leqslant n_{d}-1}} 
\sum_{\xi \in \mu_{N}(K)}
\sum_{\substack{w_{1},w_{2} \text{words }
	\\w = e_{0}^{l_{d}}w_{1}w_{2}e_{0}^{l_{0}}
	\\ \depth(w_{1})\geqslant 1
	\\ \depth(w_{2}) \geqslant 1}}
 \Ad_{\Phi^{(\xi)}_{p,\alpha}}(e_{\xi})[w_{1}] \text{ }\text{ } (-1)^{l_{0}} {l_{0}+l_{d}-1 \choose l_{0}} e_{0}^{l_{d}+l_{0}} e_{\xi^{(p^{\alpha})}} \Reg\dec_{p,\alpha}(w_{2})
\\ - \sum_{0 \leqslant l_{0} \leqslant n_{0}-1}
(-1)^{l_{0}} {l_{0}+l_{d}-1 \choose l_{0}}
\text{ }\text{ } e_{0}^{n_{d}-1+l_{0}}e_{\xi_{d}}
\Reg\dec_{p,\alpha}(e_{0}^{n_{d-1}-1}\ldots e_{0}^{n_{1}-1} e_{\xi_{1}})
\end{multline}
\end{Proposition}

\begin{proof} In the first line of equation (\ref{definition of reg dec}), the terms such that $\depth(w_{2})=0$ of $\depth(w_{1})=0$ vanish by Lemma \ref{lemma depth 0}. Thus, the first line of equation (\ref{definition of reg dec}) can be replaced by 
$$  \sum_{\substack{w_{1},w_{2}\text{ words } \\ w_{1}w_{2}=w\\ \depth(w_{2})\geqslant 1\\ \depth(w_{1})\geqslant 1}} \sum_{\xi \in \mu_{N}(K)} \bigg( \Ad_{\Phi_{p,\alpha}^{(\xi)}}(e_{\xi})[w_{2}] \bigg) e_{\xi^{(p^{\alpha})}}\Reg\dec_{p,\alpha}(w_{1}) $$
The result is then obtained by induction on $(n_{d},n_{0})$ for the lexicographical order on $(\mathbb{N}^{\ast})^{2}$, noting that if if $(n_{0},n_{d})=(1,1)$ the formula for $\Reg\dec_{p,\alpha}$ is already inductive with respect to the depth and that $(n_{0},n_{d})$ decreases when applying $\Reg\dec_{p,\alpha}$ in the other cases.
\end{proof}

\subsection{End of the proof of the main theorem}

Given that  $\Ad_{\Phi^{(\xi)}_{p,\alpha}}(e_{\xi})=\{\Phi^{(\xi)}_{p,\alpha}\}^{-1}e_{\xi}\Phi^{(\xi)}_{p,\alpha}$, for any word $w$ of weight $n$ and depth $d$, $\Ad_{\Phi^{(\xi)}_{p,\alpha}}(e_{\xi})[w]$ is a $\mathbb{Z}$-linear combination of $p$-adic cyclotomic multiple zeta values of weight $n-1$ and depth $d-1$ : more precisely, $\Ad_{\Phi^{(\xi)}_{p,\alpha}}(e_{\xi})[w] = \sum_{\substack{w_{1},w_{2} \text{words }\\ w_{1}e_{\xi}w_{2}=w}} \Phi_{p,\alpha} [S(w_{1})\text{  }\sh\text{ }w_{2}]$ (by the formula for the antipode and product of shuffle Hopf algebras, \S1.1.1).
This is the reason for the shift in the weight and the depth below.

\begin{Corollary} \label{end of proof of theorem}
(i) Each $\Li_{p,\alpha}^{\dagger}[w]$ is a $\mathbb{Z}$-linear combination of products $\big(\prod_{i=1}^{r} \Ad_{\Phi_{p,\alpha}^{(\xi_{i})}}(e_{\xi_{i}})[w_{i}]\big) \Reg\int_{p,\alpha}\tilde{w}$, where the $\xi_{i}$'s are in $\mu_{N}(K)$, the $w_{i}$'s are words on $e_{0\cup \mu_{N}}$ of depth $\geqslant 2$, $r \in \mathbb{N}$, and $\tilde{w}$ is a word over $e_{0\cup \mu_{N}\cup \mu_{N}^{(p^{\alpha})}}$ such that $\sum_{i=1}^{r} \big(\depth(w_{i})-1\big) + \depth(\tilde{w}) \leqslant \depth(w)$ and $\sum_{i=1}^{r} \big(\weight(w_{i})-1\big) + \weight(\tilde{w}) = \weight(w)$.
\newline (ii) Each $\Ad_{\Phi_{p,\alpha}^{(\xi)}}(e_{\xi})[w]$ is a $\mathbb{Z}$-linear combination of products
$\big(\prod_{i=1}^{r} \Ad_{\Phi_{p,\alpha}^{(\xi_{i})}}(e_{\xi_{i}})[w_{i}]\big) \Cf_{0}(\frac{p^{\alpha}z}{z-\xi}\big(\Reg\int_{p,\alpha}(\tilde{w})\big)$ and 
$\big(\prod_{i'=1}^{r'} \Ad_{\Phi_{p,\alpha}^{({\xi'}_{i'})}}(e_{{\xi'}_{i'}})[{w'}_{i'}]\big) \Cf_{0}(\frac{z^{p^{\alpha}}}{z^{p^{\alpha}}-\xi^{p^{\alpha}}}\big(\Reg\int_{p,\alpha}(\tilde{w})\big)$
where the $\xi_{i}$'s and ${\xi'}_{i'}$'s are in $\mu_{N}(K)$, the $w_{i}$'s and ${w'}_{i'}$'s are words on $e_{0\cup \mu_{N}}$ of depth $\geqslant 2$, $r \in \mathbb{N}$, and $\tilde{w}$ is a word on $e_{0\cup \mu_{N}\cup \mu_{N}^{(p^{\alpha})}}$ such that $\sum_{i=1}^{r} \big(\depth(w_{i})-1\big) + \depth(\tilde{w}) \leqslant \depth(w)$ and $\sum_{i=1}^{r} \big(\weight(w_{i})-1\big) + \weight(\tilde{w}) = \weight(w)$.
\end{Corollary}

\begin{proof} (i) follows by induction on the depth from Proposition \ref{prop characterization in terms of regularized} and Proposition \ref{induction on depth for reg dec}, noting that equation (\ref{eq:premiere equation prime}) is homogeneous for the weight and depth, in a sense which agrees with the shift explained before the statement : in the first line of (\ref{eq:premiere equation prime}), we have, for all $(w_{1},w_{2})$,
$\weight(w) = (\weight(w_{1})-1) + \weight(w_{2}) + \weight(e_{0}^{l_{d}+l_{0}}e_{\xi^{(p^{\alpha})}})$ and $\depth(w) = (\depth(w_{1})-1) + \depth(w_{2}) + \depth(e_{0}^{l_{d}+l_{0}}e_{\xi^{(p^{\alpha})}})$
and, in the second line of (\ref{eq:premiere equation prime}),
$\weight(w) = \weight(e_{0}^{n_{d}-1+l_{0}}e_{\xi_{d}}) + \weight(e_{0}^{n_{d-1}-1}\cdots e_{0}^{n_{1}-1}e_{\xi_{1}})$ and 
$\depth(w) = \depth(e_{0}^{n_{d}-1+l_{0}}) + 
\depth(e_{0}^{n_{d-1}-1}\cdots e_{0}^{n_{1}-1}e_{\xi_{1}})$. The fact that the $w_{i}$'s can be chosen of depth $\geqslant 2$ follows from the fact that the only non-zero coefficient of $\Ad_{\Phi_{p,\alpha}^{(\xi)}}(e_{\xi})$ in depth $\leqslant 1$ is $\Ad_{\Phi_{p,\alpha}^{(\xi)}}(e_{\xi})[e_{\xi}]=1$, for any $\xi \in \mu_{N}(K)$.
\newline (ii) Follows from (i) by induction on the weight, using equation (\ref{eq:vanishing of the polar part}).
\end{proof}

\begin{Corollary} \label{last corollary}For any $d$, There exists $\kappa_{d},\kappa'_{d},\kappa''_{d} \in \mathbb{R}_{+}^{\ast}$ such that, for any word $w$ of weight $n$ and depth $d$ over $e_{0 \cup \mu_{N}}$, we have $v_{A(U^{\an})}\big(\Li_{p,\alpha}^{\dagger}[w]\big) \geqslant n- \kappa_{d}-\kappa'_{d} \log(n+\kappa''_{d})$ and 
$v_{p}\big(\Ad_{\Phi_{p,\alpha}^{(\xi)}}[w]\big) \geqslant n- \kappa_{d}-\kappa'_{d} \log(n+\kappa''_{d})$ where, for $f\in A(U^{\an})$, $v_{A(U^{\an})}(f)$ is defined by $||f||=p^{-v_{A(U^{\an})}(f)}$.
\end{Corollary}

\begin{proof} This follows by induction on $d$ from Corollary \ref{end of proof of theorem}, in which we note that $r \leqslant \depth(w)$, Proposition \ref{bound on valuation}, equations (\ref{eq:mult by z sur z - xi 1}), (\ref{eq:mult by z sur z - xi 2}), (\ref{eq:mult by z sur z - xi p 1}), (\ref{eq:mult by z sur z - xi p 2}) which describe the multiplication by $\Cf_{0}\big(\frac{z}{z-\xi}\big)$ and $\Cf_{0}\big(\frac{z^{p^{\alpha}}}{z^{p^{\alpha}}-\xi^{p^{\alpha}}}\big)$ in $\LA_{\mathcal{S}_{\alpha}}(\mathbb{Z}_{p}^{(N)},K)$, and the bounds on the valuations of the coefficients $\mathcal{B}$, which appear in those equations, in Lemma \ref{def of coeff B}.
\end{proof}

The Corollary \ref{last corollary} implies the theorem, knowing that, at $d$ fixed, 
$$ n - \kappa_{d}-\kappa'_{d} \log(n+\kappa''_{d}) \underset{n\rightarrow \infty}{\longrightarrow} \infty $$

\appendix

\section{Recasting a theorem of Mahler as a  characterization of an affinoid subspace of $\mathbb{P}^{1,\an}/K$, preliminary to \S3, \S4}

This section is an interpretation of a classical theorem of Mahler in terms of rigid analytic spaces.
The result (Proposition \ref{prop Mahler general}) is a characterization of the $K$-Banach space $A(U^{\an})$ (Notation \ref{the open affine}) in terms of the power series expansion at $0$ of its elements.
\newline\indent 
Let $L$ be a $\mathbb{Q}_{p}$-Banach space. The sequence of Mahler coefficients of a sequence $(c_{m})_{m \in \mathbb{N}} \in L^{\mathbb{N}}$, or more generally, of a function $c : \mathbb{Z}_{p} \rightarrow L$ for which we denote by $c_{m}=c(m)$, is the sequence $\displaystyle\bigg(\sum_{j=0}^{m}(-1)^{j}{m \choose j}c_{m-j}\bigg)_{m \in \mathbb{N}}\in L^{\mathbb{N}}$.  
The following statement is a classical theorem of Mahler \cite{Mahler} ; the formulation which we reproduce can be found, as well as a simple proof, in \cite{Colmez}, Theorem 2.8 : \emph{The map $L^{\mathbb{N}} \rightarrow L^{\mathbb{N}}$ which sends a sequence to its Mahler coefficients induces an isometry between $C(\mathbb{Z}_{p},L)$ equipped with the norm of uniform convergence and the space $l_{0}^{\infty} = \{ (c'_{m}) \in L^{\mathbb{N}} \text{ }|\text{ } \displaystyle c'_{m} \underset{m \rightarrow \infty}{\longrightarrow} 0 \}$ of sequences tending to $0$, equipped with the norm $||.||_{\infty}$.} In particular, a map $\mathbb{N}^{\ast} \rightarrow L$ can be interpolated by an element of $C(\mathbb{Z}_{p},L)$ if and only if its sequence of Mahler coefficients tends to $0$.

\begin{Notation} For any $L$ a $K$-Banach space and $M$ a metric space we denote by $\mathcal{C}(M,L)$ the $K$-vector space of continuous functions $M \rightarrow L$.
\end{Notation}

\begin{Notation} \label{notation coeff Taylor}For any $f \in A(U^{\an})$, and $m \in \mathbb{N}$, the coefficient of degree $m$ in the power series expansion of $f$ at $0$ is denoted by $f[z^{m}]$.
\end{Notation}

We explain a relation between Mahler's theorem, rigid analytic functions and $\mathbb{P}^{1} \setminus \{0,1,\infty\}$. Let the operator $(z \mapsto \frac{z}{z-1}) : L^{\mathbb{N}} \rightarrow L^{\mathbb{N}}$ defined as the conjugation of the automorphism $\sum_{m \in \mathbb{N}} c_{m}z^{m} \mapsto \sum_{m \in \mathbb{N}} c_{m}\big( \frac{z}{z-1} \big)^{m}$ of $L[[z]]$ by the natural isomorphism $L[[z]] \simeq L^{\mathbb{N}}$. A simple computation shows the following :

\begin{Lemma} \label{lemma Mahler} For all $(c_{m}) \in L^{\mathbb{N}}$, ${\big((z \mapsto \frac{z}{z-1}\big)c\big)}_{m} =
\left\{
\begin{array}{ll} c_{0} & \text{ if } m=0
\\  
(-1)^{m} \sum_{m'=0}^{m} (-1)^{m'} {m-1 \choose m'-1} c_{m'} & \text{ if }m>0
\end{array} \right.$.
\end{Lemma}
The Proposition below is a new formulation and a new proof of Proposition 1 in \cite{Unver MZV}, \S5. We use Notation \ref{notation coeff Taylor}. We recall that $A(U^{\an})$ is equipped with the norm defined by $||\sum a_{m} z^{m}|| = \sup_{m \in \mathbb{N}} |a_{m}|_{p}$.
\newline 
\newline Let 
$\mathbb{Z}_{p}^{(N)} = \underset{l \rightarrow \infty}{\varprojlim} 
\mathbb{Z}/ Np^{l}\mathbb{Z}\simeq \mathbb{Z}/N\mathbb{Z}\times \mathbb{Z}_{p}$, this isomorphism being an homeomorphism where $\mathbb{Z}/N\mathbb{Z}$ is equipped with the discrete topology. Thus, $\mathbb{Z}_{p}^{(N)}$ is the disjoint union of $N$ copies $\mathbb{Z}_{p}$, regarded as the closures of  $m_{0}+N\mathbb{N}$, $m_{0} \in \{1,\ldots,N\}$, respectively.
\newline\indent Following \cite{Koblitz}, \S2, let, for $z \in U^{\an}(K)$, for all $u \in \mathbb{N}^{\ast}$ : $\displaystyle \mu_{z,N}(n+p^{u}N\mathbb{Z}_{p}) =\frac{z^{n}}{1 - z^{p^{u}N}},\text{ } n \in \{1,\ldots,p^{u}N\}$ (when $z=\infty$, this is $-1$ if $n=p^{u}N$ and $0$ otherwise). This defines a measure on $\mathbb{Z}_{p}^{(N)}$. The integral of $c \in C(\mathbb{Z}_{p}^{(N)},L)$ with respect to $\mu_{z,N}$ is $\displaystyle \int_{\mathbb{Z}_{p}} c \mu_{z,N} = \lim_{u \rightarrow \infty} \sum_{r=1}^{p^{u}N} c(r) \mu_{z,N}(r+p^{u}N\mathbb{Z}_{p})$. The statement below gives a new formulation and proof to Proposition 3.0.3 of \cite{Unver cyclotomic}, \S3.

\begin{Proposition} \label{prop Mahler general} For any complete complete field extension $L$ of $K$, we have an isometric isomorphism of Banach spaces over $K$
$$ \{f \in A(U^{\an})\otimes_{K}L\text{ }|\text{ }f(0)=0\} \simlra \mathcal{C}(\mathbb{Z}_{p}^{(N)},L)$$
defined by 
\begin{equation} \label{eq:mahler transform cyclotomic} \Cf_{0} :  f \mapsto \text{the interpolation of }(m \in\mathbb{N}^{\ast} \mapsto f[z^{m}]\in L) 
\end{equation}
whose inverse is
\begin{equation} \label{eq: mahler transform inverse cyclotomic} \Cf_{0}^{-1} : c \mapsto \text{ the map }(z \mapsto \int_{\mathbb{Z}_{p}^{(N)}} c d\mu_{z,N} \big) 
\end{equation}
\end{Proposition}

\begin{proof} (a) By definition, we have $\{f \in A(\mathbb{Z}_{p}^{\an})\otimes_{\mathbb{Q}_{p}} L\text{ }|\text{ }f(0)=0\} = \{  \sum_{m \in \mathbb{N}} c_{m}z^{m} 
\text{ }|\text{ } c_{m} \in L^{\mathbb{N}}, c_{m}\underset{m \rightarrow +\infty}{\longrightarrow} 0 \}$. We apply to this equality the map $z \mapsto \frac{z}{z-1}$.
\newline\indent On the one hand, since $z \mapsto \frac{z}{z-1}$ is the unique homography that sends $(0,1,\infty) \mapsto (0,\infty,1)$, it induces an involutive isomorphism between $U^{\an} = (\mathbb{P}^{1,\an} - B(1,1))/\mathbb{Q}_{p}$ and the rigid analytic unit disk $\mathbb{Z}_{p}^{\an} = (\mathbb{P}^{1,\an} - B(\infty,1))/ \mathbb{Q}_{p}$, and induces an isomorphism between $\{f \in A(\mathbb{Z}_{p}^{\an})\otimes_{\mathbb{Q}_{p}} L\text{ }|\text{ }f(0)=0\}$ and 
$\{f \in A(U^{\an})\otimes_{\mathbb{Q}_{p}} L\text{ }|\text{ }f(0)=0\}$, and the formula in Lemma \ref{lemma Mahler} implies easily that this map is an isometry.
\newline\indent On the other hand, by Lemma \ref{lemma Mahler} and Mahler's theorem, the same map induces an isometric isomorphism between $\{ \sum_{m \in \mathbb{N}} c_{m}z^{m} 
\text{ }|\text{ } c_{m} \in L^{\mathbb{N}}, c_{m}\underset{m \rightarrow +\infty}{\longrightarrow} 0 \}$ and
$\{ \sum_{m \in \mathbb{N}} c_{m}z^{m} 
\text{ }|\text{ } c_{m} \in L^{\mathbb{N}}, m \mapsto c_{m} \in  \mathcal{C}(\mathbb{Z}_{p},L)\}$.
\newline\indent Thus, the map $z \mapsto \frac{z}{z-1}$ sends the equality with which we started to the statement that the map $\Cf_{0}$ is a linear isometry $\{f \in A(U^{\an}\otimes_{\mathbb{Q}_{p}} L)\text{ }|\text{ }f(0)=0\} \simlra  \mathcal{C}(\mathbb{Z}_{p},L)$.
\newline (b) We first prove that is $f \in A(U^{\an}) \otimes_{K}L$, has power series expansion $\sum c_{m} z^{m}$ at $0$, then the map $m \mapsto c_{m}$ interpolates to an element of $\mathcal{C}(\mathbb{Z}_{p}^{(N)},L)$.
\newline\indent The elements of $A(U^{\an})\otimes_{K}L$ are the uniform limits of sequences of rational fractions over $\mathbb{P}^{1}$ whose poles are in $\cup_{\xi \in \mu_{N}(K)} B(\xi,1)$ ; moreover, the uniform norm $||f||_{p}=\sup_{z} |f(z)|_{p}$ is equivalent to our chosen norm $||\sum a_{m}z^{m}||=\sup_{m}|a_{m}|_{p}$. Thus, by density, it is enough to prove that the statement is true for these rational fractions. Without loss of generality, we can assume $L$ algebraically closed.
\newline\indent We are reduced to prove the statement for a function $f(z)=\displaystyle\frac{1}{(z-z_{0})^{\nu}}$ with $z_{0} \in B(\xi,1)$ for a certain $\cup_{\xi \in \mu_{N}(K)}$ and $\nu \in \mathbb{N}^{\ast}$. (Indeed, the polynomial part of a partial fraction decomposition of a rational function in $A(U^{\an})\otimes_{K}L$ is necessarily constant, because such a rational fraction is bounded at $\infty$.)
We have $\xi^{-1}z_{0} \in B(\xi,1)$ and  
$\displaystyle\frac{1}{(z-z_{0})^{\nu}} = 
\frac{\xi^{-\delta}}{(\xi^{-1}z-\xi^{-1}z_{0})^{\nu}} = f(\xi^{-1}z)$ where $f \in A(\mathbb{P}^{1,an} \setminus B(1,1))\otimes_{\mathbb{Q}_{p}} L$.
\newline\indent By (a), letting $f=\sum_{m \in \mathbb{N}} c_{m}z^{m}$ be the power series expansion of $f$ at $0$, the map $m \in \mathbb{N}^{\ast} \mapsto c_{m}\in L$ defines an element of $\mathcal{C}(\mathbb{Z}_{p},L)$. We have $f(\xi^{-1}z)= = \sum_{m \in \mathbb{N}} c_{m}\xi^{-m}z^{m}$ and, since $m \mapsto \xi^{-m}$ is constant on each class of congruence modulo $N$, we deduce the map $m \in \mathbb{N}^{\ast} \mapsto c_{m}\xi^{-m} \in L$ defines an element of $\mathcal{C}(\mathbb{Z}_{p}^{(\mathbb{N})},L)$.
\newline\indent This proves that the map defined in (\ref{eq:mahler transform cyclotomic}) sends $A(U^{\an})\otimes_{K}L\text{ }|\text{ }f(0)=0\}$ to $\mathcal{C}(\mathbb{Z}_{p}^{(N)},L)$. The fact that it is an isometry its injectivity follow from its definition.
\newline\indent (c) Let $c \in \mathcal{C}(\mathbb{Z}_{p}^{(N)},L)$. For all $z \in U^{\an}(L)$. We have 
$\int_{\mathbb{Z}_{p}} c d\mu_{z,N}
= \lim_{u\rightarrow \infty} \sum_{m=1}^{Np^{u}} \frac{c(m) z^{m}}{1- z^{Np^{u}}} = 
\lim_{u\rightarrow \infty}\sum_{m=1}^{+\infty} c(m \mod Np^{u}) z^{m}$ where $m \mod Np^{u}$ is the unique integer in $\{1,\ldots,Np^{u}\}$ congruent to $m$ modulo $Np^{u}$ and the limit is in $L$.
\newline\indent Since $\mathbb{Z}_{p}^{(N)}$ is compact, $c$ is uniformly continuous ; this implies that the sequence $(\sum_{m=1}^{+\infty} c(m \mod Np^{u}) z^{m})_{u \in \mathbb{N}}$ is Cauchy for the norm $||.||_{\infty}$ on the bounded elements of $L^{\mathbb{N}}$. Moreover, $\sum_{m=1}^{+\infty} c(m \mod Np^{u}) z^{m}=\sum_{m=1}^{p^{u}} \frac{c(m) z^{m}}{1- z^{Np^{u}}}$ is a rational function over $\mathbb{P}^{1}$ whose poles are in $\cup_{\xi \in \mu_{N}(K)}B(\xi,1)$.
\newline\indent Given that $A(U^{\an})$ is the space of uniform limits of rational functions over $\mathbb{P}^{1}$ whose poles are in $\cup_{\xi \in \mu_{N}(K)}B(\xi,1)$, and that the uniform norm $||f||_{\infty} = \sup_{z}|f(z)|_{p}$ is equivalent to our chosen norm on $A(U^{\an})$, the sequence $(\sum_{m=1}^{+\infty} c(m \mod p^{u}) z^{m})_{u \in \mathbb{N}}$, sequence converges to an element of $A(U^{\an})\otimes L$. 
\newline\indent Moreover, the limit in question is the element defined by $\sum_{m=1}^{+\infty} (\lim_{u\rightarrow \infty}c(m \mod Np^{u})) z^{m} = \sum_{m=1}^{\infty} c(m) z^{m}$.
\newline\indent This implies that the map defined in (\ref{eq: mahler transform inverse cyclotomic}) is a map $\mathcal{C}(\mathbb{Z}_{p},L)\rightarrow \{f \in A(U^{\an}\otimes_{\mathbb{Q}_{p}} L)\text{ }|\text{ }f(0)=0\}$, inverse of the map defined in (\ref{eq:mahler transform cyclotomic}).
\end{proof}

\begin{Remark} By considering quotients, Proposition \ref{prop Mahler general} gives a statement on Tate algebras in one variable. Further, by gluing affinoid curves, it gives a surprising alternative description of rigid analytic curves over $\mathbb{Q}_{p}$. We will develop this remark as well as its generalization to more general rigid analytic spaces in another paper.
\end{Remark}

\section{An adelic interpretation of the computation}

Our proof of the main theorem also gives, aside bounds on norms, a inductive way to compute overconvergent $p$-adic multiple polylogarithms and $p$-adic cyclotomic multiple zeta values. This appendix answers to the question of characterizing which type of formula we obtain and, at the same time, to the question of finding a $p$-adic analogue of the type of series appearing in equation (\ref{eq:multizetas}).
\newline\indent The following numbers will play a central role in our explicit theory of $p$-adic cyclotomic multiple zeta values.

\begin{Definition} We call prime weighted cyclotomic multiple harmonic sums the numbers 
$$ \har_{p^{\alpha}} \big((n_{i})_{d};(\xi_{i})_{d+1}\big) = (p^{\alpha})^{n_{1}+\ldots+n_{d}} \frak{h}_{p^{\alpha}} \big((n_{i})_{d};(\xi_{i})_{d+1}\big) $$
\end{Definition}

The prime weighted cyclotomic multiple harmonic sums satisfy an inequality "valuation $\geqslant$ weight" :

\begin{Lemma} \label{lemma mhs valuation geq weight}We have $v_{p}\big(\har_{p^{\alpha}}\big((n_{i})_{d};(\xi_{i})_{d+1}\big)\big) \geqslant  \sum_{i=1}^{d} n_{i}$.
\end{Lemma}

\begin{proof} We split the domain of summation $\{ (m_{i})_{d} \in \mathbb{N}^{d}\text{ }|\text{ } 0<m_{1}<\cdots<m_{d}<p^{\alpha}\}$ of multiple harmonic sums (equation (\ref{eq:multiple harmonic sums})) into the disjoint union of two subsets : first, the subset characterized by $p^{\alpha-1}|n_{i}$ for all $i$, whose contribution to the multiple harmonic sum is exactly $\har_{p}\big((n_{i})_{d};(\xi_{i}^{p^{\alpha-1}})_{d+1}\big)$ which has valuation $\geqslant  \sum_{i=1}^{d} n_{i}$, and, secondly, its complement, whose contribution has higher $p$-adic valuation.
\end{proof}

The formulas obtained by our proof can be written in terms of prime weighted multiple harmonic sums. However, since they are only elements of the $N$-th cyclotomic field which plays the role of a field of coefficients, in order to be able to talk about "linear combinations" of them without having them absorbed in the field of coefficients, we consider all values of $p$ and $\alpha$ at the same time. Let us denote by $K=K_{p}$, and let $\mathcal{P}_{N}$ be the set of prime numbers which do not divide $N$.
\newline\indent We note that the constants $\kappa_{d},\kappa'_{d},\kappa''_{d}$ of Corollary \ref{last corollary} can be chosen independent of $p$ (this follows from the fact that $-\frac{d}{\log(p)}$ and $-\frac{d}{\log(p)}\log\big(\frac{d}{\log(p)}\big)$, which appear in Proposition \ref{bound on valuation}, are lower bounded uniformly in $p$). Let $C_{N}$ be the $N$-th cyclotomic field, which we regard as included in all $K_{p}$, $p \in \mathcal{P}_{N}$. We denote by $\widehat{\mathcal{O}_{\text{Bound}(d)}^{\sh,e_{0\cup \mu_{N},d}}}$ the set of formal infinite sums $\sum_{n \in \mathbb{N}} w_{n}$ where $w_{n}$ is a $C_{N}$-linear combination of words of weight $n$ and depth $\leqslant d$ with coefficients in $\{x  \in C_{N} | \forall p \in \mathcal{P}_{N}, v_{p}(x) \geqslant - \kappa_{d} - \kappa'_{d}\log(n+\kappa''_{d}) \}$.

\begin{Definition} \label{def adelic space}We fix a positive integer $d$. Let $\widehat{\Har}_{\mathcal{P}_{N}^{\mathbb{N}^{\ast}},d}$ be the $\mathbb{Z}$-module defined as the image of the map $\displaystyle\widehat{\mathcal{O}_{\text{Bound}(d)}^{\sh,e_{0\cup \mu_{N},d}}} \rightarrow \prod_{(p,\alpha)\in \mathcal{P}_{N}\times \mathbb{N}^{\ast}} K_{p}$ which sends $\sum\limits_{n \geqslant 0} w_{n} \mapsto \big( \sum\limits_{n \geqslant 0} \har_{p^{\alpha}}(w_{n}) \big)_{(p,\alpha) \in\mathcal{P}_{N}\times \mathbb{N}^{\ast}}$.
\end{Definition}

\begin{Lemma}
The series $\sum_{n \geqslant 0} \har_{p^{\alpha}}(w_{n})$ as above are convergent in $K_{p}$, uniformly with respect to $p$, $\alpha$, and $\sum_{n\geqslant 0} w_{n})$ in $\mathcal{P}_{N}$, $\mathbb{N}^{\ast}$ and $\widehat{\mathcal{O}_{\text{Bound}(d)}^{\sh,e_{0\cup \mu_{N},d}}}$ respectively. In particular, for any $n_{0}\in\mathbb{N}$, the reduction of $\widehat{\Har}^{\leqslant d}_{\mathcal{P}_{N}^{\mathbb{N}^{\ast}}}$ modulo $(p^{n_{0}})_{(p,\alpha) \in\mathcal{P}_{N}\times \mathbb{N}^{\ast}}$ is a finitely generated $\mathbb{Z}$-module.
\end{Lemma}

\begin{proof} The first part of the statement follows from Lemma \ref{lemma mhs valuation geq weight}. It implies that, for any $n_{0} \in \mathbb{N}^{\ast}$, we have $\sum_{n \geqslant 0} \har_{p^{\alpha}}(w_{n}) \equiv \sum_{n =0}^{n_{1}} \har_{p^{\alpha}}(w_{n}) \mod p^{n_{0}}$ where $n_{1}$ can be chosen uniformly in $p$, $\alpha$ and $\sum_{n\geqslant 0} w_{n}$.
\end{proof}

We will state a notion formalizing the sequences $(\har_{p^{\alpha}}(w))_{(p,\alpha)}$ in \cite{II-1}. The kernel of the map in Definition \ref{def adelic space} will be studied and formalized in \cite{II-1}, \cite{II-2}, \cite{II-3}.
\newline\indent In the next Proposition, we take the convention that, for $r_{0} \in \mathbb{N}^{\ast}$, $n \in \mathbb{N}^{\ast}$, $\xi \in \mu_{N}(K)$, $\frak{h}_{r_{0}}(n,\xi)= \xi^{r_{0}}{r_{0}^{-n}}$ (this is compatible with equation (\ref{eq:multiple harmonic sums}) which involves only indices of the type $((n_{i})_{d},(\xi_{i})_{d+1})$).

\begin{Proposition} (i) Any sequence $\big( \Ad_{\Phi^{(\xi)}_{p,\alpha}}(e_{\xi})[w]\big)_{(p,\alpha)\in \mathcal{P}_{N} \times \mathbb{N}^{\ast}}$ with $w$ a word over $e_{0\cup \mu_{N}}$ of depth $d$, $\xi \in \mu_{N}(K)$, is an element of $\widehat{\Har}_{\mathcal{P}_{N}^{\mathbb{N}^{\ast}},d}$.
\newline (ii) Let a word $w \in \mathcal{O}^{\sh,e_{0\cup \mu_{N}}}$,  $l \in \mathbb{N}$ and $\xi \in \mu_{N}(K)$. We have $\displaystyle \big((\Cf_{0}\Li^{\dagger}_{p,\alpha}[w])^{(l,\xi)}(0) \big)_{(p,\alpha) \in \mathcal{P}_{N} \times \mathbb{N}^{\ast}} \in \widehat{\Har}_{\mathcal{P}_{N}^{\mathbb{N}^{\ast}},d}$ and there exist words $w'_{i} \in \mathcal{O}^{\sh,e_{0\cup \mu_{N}}}$ of depth $\leq d$ ($i$ in a finite set $I$) and elements $h_{w,l,\xi}$ and $h_{w,l,\xi,i}$ in $\widehat{\Har}_{\mathcal{P}_{N}^{\mathbb{N}^{\ast}},d}$, such that, for all $r_{0}\in \mathbb{N}^{\ast}$, $\displaystyle \big((\Cf_{0}\Li^{\dagger}_{p,\alpha}[w])^{(l,\xi)}(r_{0}) \big)_{(p,\alpha) \in \mathcal{P}_{N} \times \mathbb{N}^{\ast}, p^{\alpha}>r_{0}}$ is the image of $h_{w,l,\xi} + \sum_{i\in I} h_{w,l,\xi,i}.(\har_{r_{0}}(w'_{i}))_{\substack{(p,\alpha) \in \mathcal{P}_{N}\times \mathbb{N}^{\ast}}}$ by the canonical projection $\prod\limits_{(p,\alpha) \in \mathcal{P}_{N} \times \mathbb{N}^{\ast}} K_{p} \twoheadrightarrow \prod\limits_{\substack{(p,\alpha) \in \mathcal{P}_{N} \times \mathbb{N}^{\ast}\\ p^{\alpha}>r_{0}}} K_{p}$.
\end{Proposition}

\begin{proof} The statement (ii) is actually true for the regularized iterated integrals of \S3, by induction on the depth using Proposition \ref{prop comp reg it int}. This implies (i) and (ii) by Proposition \ref{prop characterization in terms of regularized}.
\end{proof}

We note that we have not encountered in our computation the $p$-adic cyclotomic multiple zeta values, which are coefficients of $\Phi_{p,\alpha}$, but, instead, the coefficients of $\Ad_{\Phi^{(\xi)}_{p,\alpha}}(e_{\xi})$, $\xi \in \mu_{N}(K)$. This goes back to Lemma \ref{Frobenius at 0,0}. Actually, by the formal properties of the Frobenius, the coefficients $\Ad_{\Phi^{(\xi)}_{p,\alpha}}(e_{\xi})$, $\xi \in \mu_{N}(K)$ are the versions of $p$-adic cyclotomic multiple zeta values that arise naturally if we want to make an explicit theory. We will formulate all our subsequent papers, including \cite{I-2}, \cite{I-3}, \cite{II-1} \cite{II-2} \cite{II-3}, as an explicit theory of the numbers $\Ad_{\Phi^{(\xi)}_{p,\alpha}}(e_{\xi})$, $\xi \in \mu_{N}(K)$, which will be formalized in \cite{II-1}.
The simple correspondence between the coefficients of $\Phi_{p,\alpha}$ and those of $\Ad_{\Phi^{(\xi)}_{p,\alpha}}(e_{\xi})$, $\xi \in \mu_{N}(K)$ is explained in a separated paper \cite{Assoc}, from the point of view of the theory of associators. Let us just give the simplest example, which corresponds to the ordinary $p$-adic zeta values, of both this correspondence and our computation.

\begin{Example} (Depth one (d=1) and $\mathbb{P}^{1} \setminus \{0,1,\infty\}$ (N=1)).
\newline For $n \in \mathbb{N}^{\ast}$, we have $\dec_{p,\alpha}(e_{0}^{n-1}e_{1})=e_{0}^{n-1}(e_{1}-e_{1^{(p^{\alpha})}})$ : this follows from Definition \ref{definition of reg dec}, Lemma \ref{lemma depth 0} and the fact that $\Phi_{p,\alpha}^{-1}[e_{0}^{l}]=0$ for all $l \in \mathbb{N}^{\ast}$ (this last fact is proved for $\alpha=-1$ in \cite{Unver cyclotomic}, equation (4.1.2), and the proof works for any $\alpha \in \mathbb{N}^{\ast}$). This implies by Proposition \ref{prop characterization in terms of regularized} that, for all $n \in \mathbb{N}^{\ast}$,  $\Li_{p,\alpha}^{\dagger}[e_{0}^{n-1}e_{1}](z) = (p^{\alpha})^{n} \sum_{\substack{m>0\\ p^{\alpha} \nmid m}} \frac{z^{m}}{m^{n}}$, and, by Proposition \ref{prop Mahler general} :
$$ \Li_{p,\alpha}^{\dagger}[e_{0}^{n-1}e_{1}](z) = (p^{\alpha})^{n}\int_{\mathbb{Z}_{p}} \frac{1_{p^{\alpha} \nmid m}}{m^{n}} d\mu_{z}(m) $$
where $1_{p^{\alpha} \nmid m}$ is the characteristic function of the set $\{m \in \mathbb{Z}_{p}\text{ }|v_{p}(m)<\alpha\}$. 
This generalizes Koblitz's formula 
$\displaystyle- \frac{1}{p} \log_{p} \bigg( \frac{(1-z)^{p}}{1-z^{p}} \bigg) = \int_{\mathbb{Z}_{p}^{\ast}} \frac{d\mu_{z}(x)}{x}$ (\cite{Koblitz} \S5, Lemma 1) and Coleman's formula (\cite{Coleman}, Lemma 7.2 p. 202) which is, with our notations, $\displaystyle\Li_{p,1}^{\dagger}[e_{0}^{n-1}e_{1}]=p^{n}\int_{\mathbb{Z}_{p}^{\ast}} \frac{d\mu_{z}(x)}{x^{n}}$.
\newline\indent If $n\geqslant 2$, we also have $(\Phi_{p,\alpha}^{-1}e_{1}\Phi_{p,\alpha})[e_{0}^{n-1}e_{1}]=0$ (this follows from $\Phi_{p,\alpha}^{-1}[e_{0}^{n-1}]=0$).
For $l,n\in \mathbb{N}^{\ast}$, we have 
$\displaystyle(\Phi_{p,\alpha}^{-1}e_{1}\Phi_{p,\alpha})[e_{0}^{l}e_{1}e_{0}^{n-1}e_{1}]= \Phi_{p,\alpha}^{-1}[e_{0}^{l}e_{1}e_{0}^{n-1}] = (-1)^{n+l}\Phi_{p,\alpha}[e_{0}^{n-1}e_{1}e_{0}^{l}] = (-1)^{n}{l+n-1 \choose l}\Phi_{p,\alpha}[e_{0}^{n-1+l}e_{1}]$  (the first equality follows from $\Phi_{p,\alpha}^{-1}[e_{0}^{l}]=0$, the second one follows from the formula for the antipode of shuffle Hopf algebras (\S1.1.1) and the third one from equation (\ref{eq:shuffle equation})). In particular, for $n_{1}\geqslant 2$, $\zeta_{p,\alpha}(n_{1})=-\Phi_{p,\alpha}[e_{0}^{n_{1}-1}e_{1}] = \displaystyle \frac{(-1)^{n}}{n-1}(\Phi_{p,\alpha}^{-1}e_{1}\Phi_{p,\alpha})[e_{0}e_{1}e_{0}^{n-1}e_{1}]$. One can prove that $\dec(e_{1}e_{0}^{n-1}e_{1})=e_{1}e_{0}^{n-1}(e_{1}-e_{1^{(p^{\alpha})}})$ and $\dec(e_{0}e_{1}e_{0}^{n-1})=(-1)^{n-1}ne_{0}^{n}(e_{1}-e_{1^{(p^{\alpha})}})$, and deduce that :
\begin{equation} \label{eq:Kubota} \zeta_{p,\alpha}(n) = \frac{1}{n-1} \sum_{l \in \mathbb{N}} {-n \choose l} B_{l} \text{ } \har_{p^{\alpha}}(n+l-1)
\end{equation}
Let $L_{p}$ be the $p$-adic $L$-function of Kubota-Leopoldt and $\omega$ be the Teichm\"{u}ller character. The equation (4) p. 173 of \cite{Coleman} is equivalent to $\zeta_{p,1}(n) = p^{n} L_{p}(n,\omega^{1-n})$ for all $n \in \mathbb{N}^{\ast}$, via \cite{Furusho 2}, Example 2.10 (1). Replacing $\zeta_{p,1}(n)$ by $p^{n} L_{p}(n,\omega^{1-n})$ in the $\alpha=1$ case of (\ref{eq:Kubota}), we recover the Theorem 6.6 of \cite{R}. We note that the pole at $s=1$ of the $p$-adic $L$-function of Kubota-Leopoldt is visible on that formula.
\end{Example}


\begin{thebibliography}{50}
\bibitem[AHY]{Akagi Hirose Yasuda} K. Akagi, M. Hirose, S. Yasuda - \emph{Integrality of $p$-adic multiple zeta values and application to finite multiple zeta values}, preprint
\bibitem[Bes]{Besser} A. Besser - \emph{Coleman integration using the Tannakian formalism}, Math. Ann. 322 (2002) n°1, 19-48.
\bibitem[Cha]{Chatzis} A. Chatzistamatiou - \emph{On integrality of $p$-adic iterated integrals}, J. Algebra 474 (2017), 240-270.
\bibitem[Che]{Chen} K. T. Chen - \emph{Iterated path integrals}, Bull. Amer. Math. Soc. 83 (1977), n°5, 831-879. 
\bibitem[CL]{CLS} B. Chiarellotto, B. Le Stum - \emph{F-isocristaux unipotents} - Compositio Math. 116 (1999), 81-110
\bibitem[Co]{Coleman} R. Coleman - \emph{Dilogarithms, regulators and $p$-adic $L$-functions} - Invent Math, 69 (1982) n°2, 171-208
\bibitem[Cz]{Colmez} P. Colmez - \emph{Fonctions d'une variable $p$-adique, notes du cours de M2} available at  https://webusers.imj-prg.fr/$\sim$pierre.colmez/fonctionsdunevariable.pdf
\bibitem[D]{Deligne} P. Deligne - \emph{Le groupe fondamental de la droite projective moins trois points}, Galois Groups over $\mathbb{Q}$ (Berkeley, CA, 1987), Math. Sci. Res. Inst. Publ. 16, Springer-Verlag, New York, 1989.
\bibitem[DG]{Deligne Goncharov} P. Deligne, A. B. Goncharov - \emph{Groupes fondamentaux motiviques de Tate mixtes}, Ann. Sci. Ecole Norm. Sup. 38 (2005), n°1, 1-56
\bibitem[F1]{Furusho 1} H. Furusho - \emph{p-adic multiple zeta values I -- p-adic multiple polylogarithms and the p-adic KZ equation}, Invent Math, 155 (2004), n°2, 253-286.
\bibitem[F2]{Furusho 2} H. Furusho - \emph{p-adic multiple zeta values II -- Tannakian interpretations}, Amer. J. Math, 129 (2007), n°4, 1105-1144.
\bibitem[G]{Goncharov} A. B. Goncharov - \emph{Multiple polylogarithms and mixed Tate motives} arXiv:0103059
\bibitem[J]{Note p-adique} D. Jarossay - \emph{Un cadre explicite pour les multiz\^{e}tas $p$-adiques et les polylogarithmes multiples $p$-adiques}, C. R. Math Acad Sci Paris 353 (2015) 871-876
\bibitem[J.I-2]{I-2} D. Jarossay, \emph{Pro-unipotent harmonic actions and a computation of $p$-adic cyclotomic multiple zeta values}, arXiv:1501.04893 (submitted)
\bibitem[J.I-3]{I-3} D. Jarossay, \emph{Pro-unipotent harmonic actions and a few dynamical aspects of the computation of $p$-adic cyclotomic multiple zeta values}, arXiv:1610.09107 (submitted)
\bibitem[J.II-1]{II-1} D. Jarossay, \emph{Adjoint cyclotomic multiple zeta values and cyclotomic multiple harmonic values}, arXiv:1412.5099
\bibitem[J.II-2]{II-2} D. Jarossay, \emph{The adjoint quasi-shuffle relation of $p$-adic cyclotomic multiple zeta values recovered by the explicit formulas}, arXiv:1601.01158
\bibitem[J.II-3]{II-3} D. Jarossay, \emph{Cyclotomic multiple harmonic values regarded as periods}, arXiv:1601.01159
\bibitem[J.Assoc]{Assoc} D. Jarossay, \emph{Associators, adjoint actions and the depth filtration}, arXiv:1601.01161 (submitted)
\bibitem[Ko]{Koblitz} N. Koblitz, \emph{$p$-adic analysis : a short course on recent work}, London Math Soc, Lecture Note Series 46, November 1980
\bibitem[M]{Mahler} K. Mahler - \emph{An interpolation series for continuous functions of a p-adic variable}, J. Reine Angew Math 199 (1958), 23-34.
\bibitem[R]{R} J. Rosen,
\emph{Asymptotic relations for truncated multiple zeta values},
J. London Math Soc, 91, (2015), n°2, 554-572.
\bibitem[S1]{Shiho 1} A. Shiho - \emph{Crystalline fundamental groups. I. Isocrystals on log crystalline site and log convergent site}, J. Math. Soc. Univ Tokyo 7 (2000), n°4, 509-656
\bibitem[S2]{Shiho 2} A. Shiho - \emph{Crystalline fundamental groups. II. Log convergent cohomology and rigid cohomology}, J. Math. Soc. Univ. Tokyo 9, (2002), n°1, 1-163
\bibitem[U1]{Unver MZV} S. \"{U}nver - \emph{$p$-adic multi-zeta values}. J Number Theory, 108 (2004), 111-156
\bibitem[U2]{Unver cyclotomic} S. \"{U}nver - \emph{Cyclotomic p-adic multi-zeta values in depth two}, Manuscripta Math, 149 (2016), n°3-4, 405-441
\bibitem[U3]{U3} S. \"{U}nver - \emph{A note on the algebra of $p$-adic multi-zeta values}, Commun. Number Theory Phys, 9 (2015), n°4, 689-705
\bibitem[U4]{U4} S. \"{U}nver - \emph{Cyclotomic $p$-adic multi-zeta values}, arXiv:1701.05729
\bibitem[V]{Vologodsky} V. Vologodsky, \emph{Hodge structure on the fundamental group and its application to $p$-adic integration}, Moscow Math J 3 (2003), n°1, 205-247.
\bibitem[Y]{Yamashita} G. Yamashita, \emph{Bounds for the dimension of $p$-adic multiple $L$-values spaces}. Documenta Math, Extra Volume Suslin (2010), 687-723
\end{thebibliography}
\end{document}